\documentclass[11pt]{article}%
\usepackage{amsfonts}
\usepackage{amsmath}
\usepackage{amssymb}
\usepackage{graphicx}
\usepackage{geometry}%
\usepackage{enumitem}
\usepackage{color}
\usepackage{bbm}
\usepackage[titletoc,title]{appendix}
\usepackage{subcaption}
\usepackage{multirow}
\usepackage{hyperref}

\providecommand{\U}[1]{\protect\rule{.1in}{.1in}}
\newtheorem{theorem}{Theorem}
\newtheorem{proposition}{Proposition}
\newtheorem{lemma}{Lemma}
\newtheorem{corollary}{Corollary}
\newtheorem{example}{Example}
\newtheorem{remark}{Remark}
\newtheorem{definition}{Definition}

\normalsize
\DeclareRobustCommand{\lcroof}[1]{
  \hbox{\vtop{\vbox{%
      \hrule\kern 1pt\hbox{%
        $\scriptstyle #1$%
        \kern 1pt}}\kern1pt}%
    \vrule\kern1pt}}

\numberwithin{equation}{section}
\allowdisplaybreaks
\begin{document}

\author{Eric C.K. Cheung\thanks{Corresponding author: School of Risk and Actuarial Studies, UNSW Business School, University of New South Wales, Sydney, NSW 2052, Australia. Email: eric.cheung@unsw.edu.au}, Landy Rabehasaina\thanks{Laboratoire de Math\'ematiques, University Bourgogne Franche Comt\'e, 16 route de Gray, 25030 Besan\c con c\'edex, France. Email: lrabehas@univ-fcomte.fr}, Jae-Kyung Woo\thanks{School of Risk and Actuarial Studies, UNSW Business School, University of New South Wales, Sydney, NSW 2052, Australia. Email: j.k.woo@unsw.edu.au}, Ran Xu\thanks{Department of Mathematics and Statistics, Concordia University, 1455 de Maisonneuve Blvd W, Montreal, Quebec H3G 1M8, Canada. Email: ran.xu@concordia.ca}}
\title{Asymptotic Correlation Structure of Discounted Incurred But Not Reported Claims under Fractional Poisson Arrival Process}

\maketitle

\vspace{-.3in}

\begin{abstract}
This paper studies the joint moments of a compound discounted renewal process observed at different times with each arrival removed from the system after a random delay. This process can be used to describe the aggregate (discounted) Incurred But Not Reported claims in insurance and also the total number of customers in an infinite server queue. It is shown that the joint moments can be obtained recursively in terms of the renewal density, from which the covariance and correlation structures are derived. In particular, the fractional Poisson process defined via the renewal approach is also considered. Furthermore, the asymptotic behaviour of covariance and correlation coefficient of the aforementioned quantities is analyzed as the time horizon goes to infinity. Special attention is paid to the cases of exponential and Pareto delays. Some numerical examples in relation to our theoretical results are also presented.
\end{abstract}

\vspace{-.1in} 

	\section{Introduction}\vspace{-.05in}

To model Incurred But Not Reported (IBNR) claims in insurance, in this paper it is first assumed that the claim arrivals follow a renewal process and the reporting delay for each claim is arbitrary. It is known (e.g. Willmot (1990), Mikosch (2009) and Ross (2014)) that unreported claims in actuarial science can be connected to quantities considered in other fields such as the number of customers or particles in a system with an infinite server queue structure (i.e. customers are served immediately upon their arrivals by one of those many servers). There is a vast literature on systems involving infinite server queues. See e.g. Brown and Ross (1969), Keilson and Seidmann (1988), Liu et al. (1990), and Keilson and Servi (1994) for some classical results; and Ridder (2009), Pang and Whitt (2012), Blom and Mandjes (2013), Jansen et al. (2016), Moiseev and Nazarov (2016), and Blom et al. (2017) for some recent development of the subject. As mentioned in Mikosch (2009, Section 8.2.4), the counting process of the IBNR claims can be viewed as a generating model for the activities of processing packets in a large data network. In this case, the number of active sources at time $t$ is the subject of interest of that model. Aggregate IBNR claim process is also related to the cumulative shock model with delayed termination (e.g. Finkelstein and Cha (2013, Section 4.4)), where the reporting time of a claim (i.e. the claim arrival time plus its reporting lag) corresponds to the time of the effective event which terminates fault. Interested readers are also referred to Blanchet and Lam (2013) for the application of infinite serve queues in large insurance portfolios. See also Badescu et al. (2016a,b) for the theoretical properties and parameter estimations of a marked Cox model for IBNR claims.

In the analysis of such stochastic models, the knowledge of correlation structure of the underlying arrival process is often useful to understand the joint behaviour of the process of our interest observed at different times. Indeed, the correlation between two instants $s$ and $t$ when $t\to\infty$ can be used to assess if the process has a long memory, i.e. if the state of the process at a given time has significant impact on its state at a much later time, and can serve to quantify such impact. This kind of information can be of paramount importance for decision making. In particular, the notion of long-range dependence appears when the correlation function exhibits the form of power-law decay. This definition arises in different fields such as finance, data network, and earthquake modelling. From a statistical point of view, the theoretical results allow the possibility of fitting a model to physical phenomena where long-range dependence or short-range dependence is practically observed, as explained in e.g. in Mikosch (2009, Section 8.2.4). In this paper, we shall focus on studying the covariance and correlation structures of compound renewal sums at different times (in the presence of a random delay in realization) as well as their asymptotic behaviours. To this end, the related joint moments are first derived, which may be of interest in its own right as they can in principle be used in moment-based approximations. It should be also noted that the joint moments and covariance structure of the aggregate discounted claims without any reporting delay were studied by L\'{e}veill\'{e} and Ad\'{e}kambi (2010, 2012) and L\'{e}veill\'{e} and Hamel (2013).

To model the long-range dependence property for the renewal process governing the claim arrival process, a fractional generalization of the Poisson process, known as the fractional Poisson process, will be utilized. It is known that there are different ways to define the fractional Poisson process such as using time-changed processes (e.g. Leonenko et al. (2014)). Here we adopt the definition obtained by a renewal-type treatment (e.g. Mainardi et al. (2004), and Meerschaert et al. (2011)), where the fractional Poisson process is viewed as a renewal process with Mittag-Leffler interarrival times. Consequently, the interarrival times are heavy-tailed and have infinite mean (e.g. Repin and Saichev (2000)) as opposed to the light-tailed exponential interarrival times with finite mean in the Poisson process. The fractional Poisson process is an interesting choice in the context of actuarial science, and in particular it can be used to model the arrival of claims caused by rare and extreme events such as storms and high-magnitude earthquakes in line with the discussions in e.g. Benson et al. (2007) and Biard and Saussereau (2014). Indeed, Benson et al. (2007) commented that the use of this type of renewal process is a critical extension from the Poisson assumption because the application of a Poisson model for a geologic process with heavy-tailed waiting times between events will result in a significant misrepresentation of the associated risk. Our analysis will shed light on how this heavy-tail feature impacts the long-term correlation of the IBNR process.


In what follows, we shall define the model in the context of actuarial science. It is assumed that the claim counting process $\{N(t)\}_{t\ge0}$ is a renewal process with the sequence of arrival times $\{T_{n}\}^{\infty}_{n=1}$ where $T_{n} :=\sum_{i=1}^{n} \tau_i$ for $n\in\mathbb{N}$. Here $\mathbb{N}$ is the set of positive integers, and $\{ \tau_i\}_{i=1}^{\infty}$ represents the sequence of interarrival times that are independent and identically distributed (iid) with cumulative distribution function (cdf) $F(t)$ and probability density function (pdf) $f(t)=F'(t)$. We adopt the usual convention that $T_{0}:=0$. Therefore, the claim count is defined by $N(t):= \max\{ i\in \mathbb{N}\cup \{ 0\}:T_{i}\leq t\}$ with $N(0):=0$. The renewal function of the claim counting process shall be denoted by $m(t):=\mathbb{E}[N(t)]=\sum_{i=1}^{\infty}F^{\ast(i)}(t)$, where $F^{\ast(i)}(t)$ is the $i$-th fold convolution of the cdf $F(t)$ with itself and $f^{\ast(i)}(t):=\mathrm{d} F^{\ast(i)}(t)/\mathrm{d}t$ is its pdf. Let us introduce $\{X_{i}\}_{i=1}^{\infty}$ as the sequence of claim amounts which are also assumed to be iid. Denoting $L_i$ as the time-lag (or reporting delay) corresponding to the $i$-th claim arrival, the time-lags $\{L_{i}\}_{i=1}^{\infty}$ are assumed to form an iid sequence with common cdf $W(t):=1-\overline{W}(t)$ and pdf $w(t):=W'(t)$. It is assumed that the claim counting process $\{N(t)\}_{t\ge0}$ (or the interarrival times $\{ \tau_i\}_{i=1}^{\infty}$), the claim amounts $\{X_{i}\}_{i=1}^{\infty}$ and the delays $\{L_{i}\}_{i=1}^{\infty}$ are mutually independent. Such independence assumption is standard in the actuarial literature (see e.g. Karlsson (1974) and Willmot (1990)), and we refer interested readers to Landriault et al. (2017) for some discussions on possible dependence assumptions on the triplet $(\tau_i,X_i,L_i)$. For later use, the generic random variables of $\tau_i$, $X_i$ and $L_i$ are denoted by $\tau$, $X$ and $L$ respectively. 
It will be seen that the moments of the claim amount are sufficient for our analysis, and therefore we define $\mu_k :=\mathbb{E}[X^{k}]$ and $\mu_0:=1$ for convenience.

The total discounted IBNR claim process $\{Z_\delta(t)\}_{t\ge0}$ is now defined by
\begin{equation}\label{defZt}
Z_\delta(t):=\sum_{i=1}^{N(t)} e^{-\delta (T_{i}+L_{i})} I_{\{ T_{i}+L_{i}>t\}} X_{i}=\sum_{i=1}^{\infty} e^{-\delta (T_{i}+L_{i})} I_{\{ T_i \leq t < T_{i}+L_{i}\}} X_{i},
\end{equation}
where $\delta\geq 0$ is a constant force of interest. We shall adopt the convention that $\sum_{j}^{i}:=0$ whenever $i<j$. For notational simplicity, we let $Z(t)=Z_0(t)$ when $\delta=0$. Note that, although the model is explained in an actuarial setting, it can be applied to a queueing context via a switch of terminology by changing `claim arrival' and `reporting delay' to `customer arrival' and `service time' respectively. In particular, if the distribution of $X$ is assumed to be a point mass at one, then $Z(t)$ can be interpreted as the total number of customers at time $t$ in a $G/G/\infty$ queue. (In such case, the assumption of independence between $\{N(t)\}_{t\ge0}$ and $\{L_{i}\}_{i=1}^{\infty}$ means that the service offered to customers does not depend on their arrival times, which is a standard assumption in $G/G/\infty$ queue.) Consequently, our formulation under an actuarial context is more general than a $G/G/\infty$ queue, as (i) insurance companies possibly take into account the effect of interest (allowing for $\delta>0$) when discounting claims; and (ii) actuaries associate a loss amount $X$ to each claim event when calculating or estimating outstanding claims. In this paper, we are interested in the covariance and correlation structures of $Z_\delta(s)$ and $Z_\delta(t)$ for $0<s\leq t$, which are defined as $\mathbb{C}\mathrm{ov}[Z_\delta(s), Z_\delta(t)]=\mathbb{E}[Z_\delta(s)Z_\delta(t)] - \mathbb{E}[Z_\delta(s)] \mathbb{E}[Z_\delta(t)]$ and $\mathrm{Corr}[Z_\delta(s), Z_\delta(t)]=\mathbb{C}\mathrm{ov}[Z_\delta(s),Z_\delta(t)]/ \sqrt{\mathbb{V}\mathrm{ar}[Z_\delta(s)]\mathbb{V}\mathrm{ar}[Z_\delta(t)]}$ respectively.

Concerning the renewal process for claim arrivals, in most of our analysis it is essential to know explicitly the renewal density $m'(t)$ for $t>0$. The class of fractional Poisson processes turns out to be a good candidate thanks to its nice renewal density (see \eqref{FP}). Certainly, some results will be simplified for the Poisson process as $m'(t)=\lambda>0$ does not depend on $t$.\vspace{-.1in}

\section{Outline of results and long-range dependence of IBNR process}\vspace{-.05in}
In order to have a global overview of the paper, our main contributions are presented here. In Section \ref{sec_general_results}, the marginal moments of $Z_\delta(t)$ and the joint moments of $Z_\delta(s)$ and $Z_\delta(t)$ for $0< s \leq t$ are derived in Theorems \ref{EZn} and \ref{EZnZm} using renewal arguments, and consequently the covariance is studied in Corollary \ref{EZtZt+h} and Theorem \ref{CovZsZI}. Section \ref{FPP} is devoted to the asymptotic behaviours (as $t\to\infty$) of the mean, variance and covariance in the case of fractional Poisson arrival process. To obtain simpler formulas, distributional assumptions on the random delays are made, where exact and asymptotic results are provided in Propositions \ref{Expectation_expo_case}-\ref{CovExpPro} for the exponential case, and asymptotics in Propositions \ref{EZParPro}-\ref{CovPareto} for the Pareto case. We remark that the focus of the paper is to analyze the properties of the IBNR process and to provide probabilistic interpretations. Immediate practical applications of the results, including statistical estimation and fitting with a data set, are outside the scope of this paper and will be topics for future research. For legibility purpose, the asymptotic results obtained in the paper are summarized as follows when claims occur according to a fractional Poisson process with index $\alpha\in (0,1)$. The notation $g_1(t)\propto g_2(t)$ means that $g_1(t)\sim C g_2(t)$ as $t\to\infty$ for some constant $C>0$ that possibly depends on $s$ but not on $t$. Moreover, $L~{\buildrel d \over =}~{\cal E}(\beta)$ means that the delay $L$ is exponentially distributed with mean $1/\beta$ whereas $L~{\buildrel d \over =}~\text{Pareto}(\theta,\eta)$ means that the survival function of $L$ is given by (\ref{Ptail}). We remark that the multiplicative constants are omitted in the summary below not only for the sake of brevity but also for the fact that the asymptotic correlation up to a multiplicative factor will be sufficient for determining whether a process possesses short-range or long-range dependence (see Definition \ref{LRD2}). Nonetheless, the exact expressions for these multiplicative constants are given by \eqref{asyEZdFP}, \eqref{VARZdtasy}, \eqref{COVZdstasy}, \eqref{PEZ2}, \eqref{VARZdtasyPareto1} and \eqref{asymp_cov_Pareto} in the aforementioned propositions and can be computed or approximated (see Section \ref{SEC5}).
\begin{eqnarray}
L~{\buildrel d \over =}~{\cal E}(\beta) &: & \mathbb{E}[Z_\delta(t)]\propto e^{-\delta t}t^{\alpha-1},\quad \mathbb{V}\mathrm{ar}[Z_\delta(t)]\propto e^{-2\delta t}t^{\alpha-1},\quad \mathbb{C}\mathrm{ov}[Z_\delta(s), Z_\delta(t)]\propto e^{-\delta t}t^{\alpha-2}\nonumber\\\label{tableau_recap1}\\
L~{\buildrel d \over =}~\text{Pareto}(\theta,\eta) &: & \begin{array}{|c|c|c|c|}
\hline
 &  0<\eta<1 & \eta = 1 & \eta >1 \\
 \hline
 \mathbb{E}[Z(t)] & \propto  t^{\alpha-\eta}& \propto  t^{\alpha-1} \ln t &  \propto t^{\alpha -  1}\\
 \hline
\end{array}\nonumber\\
&& \begin{array}{|c|c|c|c|c|}
\hline
 & 0< \eta<\alpha  &  \alpha\le\eta<1 & \eta=1 & \eta>1\\
 \hline
 \mathbb{V}\mathrm{ar}[Z(t)] & \propto t^{2(\alpha-\eta)} & \propto t^{\alpha-\eta} & \propto t^{\alpha-1} \ln t & \propto t^{\alpha-1}\\
 \hline
\end{array}\nonumber\\
&& \begin{array}{|c|c|c|}
\hline
 &  0<\eta<2-\alpha & \eta\ge 2-\alpha\\
 \hline
 \mathbb{C}\mathrm{ov}[Z(s), Z(t)] & \propto t^{-\eta}& \propto t^{\alpha-2}\\
 \hline
\end{array}\label{tableau_recap2}
\end{eqnarray}

We emphasize that our contributions are not only to actuarial science but also to queueing theory. In the queueing literature, only few papers have dealt with fractional Poisson arrivals because the techniques are usually very specific to fractional calculus. See e.g. Orsingher and Polito (2011) for a particular queue with a birth and death structure modelled by fractional Poisson arrivals, which is one of the rare papers on the subject. However, there are no results in the literature concerning infinite server queues where the arrivals are modelled by a fractional Poisson process, to the best of our knowledge. The result \eqref{tableau_recap1} and the tables in \eqref{tableau_recap2} enable us to give some qualitative insight on whether the IBNR process (or the $G/G/\infty$ queue as a special case) exhibits long-range dependence under the fractional Poisson setting. The behaviour will turn out to be different depending on the distribution of the delays, as displayed in \eqref{comp_corr_expo}-\eqref{comp_corr_Pareto_Poi} thereafter. The choice of the delays, namely exponential and Pareto, is motivated by the fact that these are representatives of light-tailed and heavy-tailed distributions respectively. These distributions were also chosen for technical purpose, and we however believe that this will shed light on whether behaviours resembling those in \eqref{comp_corr_expo}-\eqref{comp_corr_Pareto_Poi} are valid more generally for wider classes of light-tailed and heavy-tailed distributions in future research. Interested readers are also referred to Resnick and Rootzen (2000) for an infinite server queue with heavy-tailed service time and Poisson arrivals, where heavy-tailed service time was argued to be a good candidate for modelling the transfer of huge files across the World Wide Web. In the following, we first state the notion of long-range dependence defined by Maheshwari and Vellaisamy (2016, Section 2.2).
\begin{definition} (Long-range dependence in Maheshwari and Vellaisamy (2016)) \normalfont\label{LRD}
	Let $\{Y(t)\}_{t\ge 0}$ be a stochastic process. Suppose that the correlation function $\mathrm{Corr}[Y(s),Y(t)]$ satisfies, for all $s>0$,
	\[\lim_{t\rightarrow \infty}\frac{\mathrm{Corr}[Y(s),Y(t)]}{t^{-d}} =c(s)
	\]
for some $d>0$ and $c(s)>0$. If $d\in (0,1]$, then $\{Y(t)\}_{t\ge 0}$ is said to possess long-range dependence property. On the other hand, if $d\in (1,2)$, then $\{Y(t)\}_{t\ge 0}$ has short-range dependence. \hfill$\square$
\end{definition}
From Leonenko et al. (2014, p.10), it is known that the fractional Poisson process with index $\alpha\in (0,1)$ has long-range dependence property. However, the above Definition \ref{LRD} only includes processes where the correlation function exhibits power decay (and it excludes the cases where the correlation decreases exponentially or follows other decay law). According to Mikosch (2009, p.283), one also has the following more general definition of long-range dependence, which is adapted to stochastic processes that are not necessarily stationary.
\begin{definition} (General notion of long-range dependence) \normalfont\label{LRD2}
	Let $\{Y(t)\}_{t\ge 0}$ be a stochastic process. Then $\{Y(t)\}_{t\ge 0}$ is said to possess long-range dependence property if one has for all $s>0$ that
\[
\int_0^\infty |\mathrm{Corr}[Y(s),Y(t)]|\mathrm{d}t=\infty.
\]
On the other hand, $\{Y(t)\}_{t\ge 0}$ has short-range dependence if, for all $s>0$,
\[
\int_0^\infty |\mathrm{Corr}[Y(s),Y(t)]|\mathrm{d}t<\infty.
\]
\end{definition}
\vspace{-.4in}
\hfill$\square$\vspace{.1in}

For convenience, from now on long-range dependence and short-range dependence will be abbreviated as LRD and SRD respectively. It is instructive to note that if a process has the LRD (resp. SRD) property under Definition \ref{LRD}, then it also has the LRD (resp. SRD) property under Definition \ref{LRD2}. In the upcoming discussion, Definition \ref{LRD2} will be adopted. In the case where the delays follow exponential distribution, one sees from (\ref{tableau_recap1}) that 
\begin{equation}\label{comp_corr_expo}
L~{\buildrel d \over =}~{\cal E}(\beta)\quad:\quad\mathrm{Corr}[Z_\delta(s),Z_\delta(t)]\propto t^{-\frac{3-\alpha}{2}},
\end{equation}
and therefore $\{Z_\delta(t)\}_{t\ge 0}$ has SRD. When the delays are $\mbox{Pareto}(\theta,\eta)$ distributed, the asymptotic behaviour of the correlation can be readily obtained using (\ref{tableau_recap2}), resulting in the following table.
\begin{align}\label{comp_corr_Pareto}
&L~{\buildrel d \over =}~\text{Pareto}(\theta,\eta)~:\nonumber\\
&\begin{array}{|c|c|c|c|c|c|c|}
\hline
  & 0<\eta<\alpha & \alpha \le \eta <1 & \eta=1 & 1<\eta\le\frac{3-\alpha}{2} & \frac{3-\alpha}{2}<\eta<2-\alpha & \eta\ge 2-\alpha\\
 \hline
 \mathrm{Corr}[Z(s),Z(t)] & \propto t^{-\alpha} & \propto t^{-(\eta+\alpha)/2} & \propto \frac{t^{-(1+\alpha)/2}}{\sqrt{\ln t}} & \propto t^{-\eta-(\alpha-1)/2} & \propto t^{-\eta-(\alpha-1)/2} & \propto t^{-(3-\alpha)/2}\\
 \hline
 \mbox{Dependence} & \mbox{LRD} & \mbox{LRD} & \mbox{LRD} & \mbox{LRD} & \mbox{SRD} & \mbox{SRD}\\
 \hline
\end{array}
\end{align}

The above results can be compared to the situation where arrivals occur according to a Poisson process (i.e. $\alpha=1$). We only look at the discount-free case (i.e. $\delta=0$) as follows.
When $L~{\buildrel d \over =}~{\cal E}(\beta)$, we use the covariance in \eqref{CovPoidzero} along with the variance in \eqref{VARdzeroPoi}
to easily check that $\mathrm{Corr}[Z(s),Z(t)]  \propto e^{-\beta t}$ in the Poisson case,
and therefore the process $\{Z(t)\}_{t\ge 0}$ has SRD. This is consistent with the fractional Poisson case (i.e. $\alpha\in (0,1)$) described in \eqref{comp_corr_expo}. In contrast, when $L~{\buildrel d \over =}~\mbox{Pareto}(\theta,\eta)$, we omit the details and state that \eqref{CovPoidzero} implies $\mathbb{C}\mathrm{ov}[Z(s), Z(t)] \propto t^{-\eta}$ in the Poisson case. Moreover, one can use \eqref{VARdzeroPoi} with the help of L'H\^opital's rule to verify that $\mathbb{V}\mathrm{ar}[Z(t)] \propto t^{-\eta+1}$ if $0<\eta<1$; $\mathbb{V}\mathrm{ar}[Z(t)] \propto \ln t$ if $\eta=1$; and $\mathbb{V}\mathrm{ar}[Z(t)] \propto 1$ if $\eta>1$, and consequently the following table is constructed.
\begin{align}\label{comp_corr_Pareto_Poi}
L~{\buildrel d \over =}~\text{Pareto}(\theta,\eta)~;~\alpha=1~:\quad
&\begin{array}{|c|c|c|c|}
\hline
  & 0<\eta<1 & \eta=1 & \eta>1\\
 \hline
 \mathrm{Corr}[Z(s),Z(t)] & \propto t^{-(\eta+1)/2} & \propto \frac{t^{-1}}{\sqrt{\ln t}} & \propto t^{-\eta} \\
 \hline
 \mbox{Dependence} & \mbox{LRD} & \mbox{LRD} & \mbox{SRD}\\
 \hline
\end{array}
\end{align}
Comparing \eqref{comp_corr_Pareto} and \eqref{comp_corr_Pareto_Poi}, it is interesting to note that $\{Z(t)\}_{t\ge 0}$ is SRD in the Poisson case but LRD in the fractional Poisson case when the Pareto delay parameter $\eta$ is such that $1<\eta\le\frac{3-\alpha}{2}$. For other values of $\eta$ the dependence property in both cases are in agreement.

These correlation behaviours can be further interpreted as follows. When the delay distribution is exponential which is light-tailed, a claim that has occurred but has not been reported at time $s$ is likely to be reported by time $t>s$ when $t$ grows large, regardless of whether the claim arrives according to a Poisson or a fractional Poisson process. This explains the SRD property. On the other hand, switching from exponential to Pareto delay which is heavy-tailed, an unreported claim at time $s$ is more likely to remain unreported at time $t>s$ (for large $t$) because ${\cal E}(\beta)\le_{\small \mbox{st},\infty}\text{Pareto}(\theta,\eta)$. Here $Y_1\le_{\small \mbox{st},\infty}Y_2$ means that $\mathbb{P}(Y_1>x)\le \mathbb{P}(Y_2>x)$ for $x$ large enough, which is a generalization of the classical stochastic order for random variables (see e.g. Shaked and Shanthikumar (2007, Chapter 1)). Consequently, the process $\{Z(t)\}_{t\ge 0}$ may have a tendency to be LRD. Note that increasing the Pareto shape parameter $\eta$ leads to stochastically smaller delay (in the sense that $\text{Pareto}(\theta_1,\eta_1)\le_{\small \mbox{st},\infty}\text{Pareto}(\theta_2,\eta_2)$ for $\eta_1>\eta_2$), which explains the observation from \eqref{comp_corr_Pareto} and \eqref{comp_corr_Pareto_Poi} that the asymptotic correlation between $Z(s)$ and $Z(t)$ decreases faster (as a function of $t$) as $\eta$ increases, and $\{Z(t)\}_{t\ge 0}$ eventually switches from LRD to SRD once $\eta$ exceeds a certain threshold. In particular, for Poisson claim arrivals (see \eqref{comp_corr_Pareto_Poi}), it is noticed that $\{Z(t)\}_{t\ge 0}$ is LRD if and only if $0<\eta\le 1$, i.e. when the Pareto delay has infinite mean. For fractional Poisson claim arrivals with index $\alpha\in (0,1)$, the interarrival times are heavy-tailed (see \eqref{asymptotics_interclaim_frac}) and have infinite mean, and hence it is natural that there is a competition between the index $\alpha$ for the interarrival time and the Pareto parameter $\eta$, and the threshold in this case is given by $\frac{3-\alpha}{2}$ as in \eqref{comp_corr_Pareto}.
\begin{remark} \normalfont
Instead of analyzing the behaviour of $\mathrm{Corr}[Z(s),Z(t)]$ as $t\to\infty$ (for fixed $s>0$), it will also be interesting to study $\mathrm{Corr}[Z(s),Z(t)]$ as both $s$ and $t$ tend to infinity. For example, one may replace $s$ and $t$ by $t$ and $t+h$ respectively and look at $\mathrm{Corr}[Z(t),Z(t+h)]$, and our exact results (e.g. Corollary \ref{EZtZt+h} and Equation \eqref{CovExp}) are still valid, from which asymptotics as $t\to\infty$ (for fixed $h>0$) may be obtained. However, we leave this as future research because it is the asymptotics of $\mathrm{Corr}[Z(s),Z(t)]$ as $t\to\infty$ that help us determine whether the IBNR process is SRD or LRD.\hfill$\square$


\end{remark}\vspace{-.2in}

\section{General results}\label{sec_general_results}\vspace{-.05in}

In this section, general recursive formulas for the joint moments of the total discounted IBNR claims $Z_\delta(s)$ and $Z_\delta(t)$ are derived without any specific distributional assumptions on the interarrival times $\{ \tau_i\}_{i=1}^{\infty}$, the reporting delays $\{L_i\}_{i=1}^{\infty}$ or the claim amounts $\{X_i\}_{i=1}^{\infty}$ (recall that the moments of the claim amounts suffice in all the analysis). Our results can be conveniently expressed in terms of the Dickson-Hipp (D-H) operator $\mathcal{T}$ (see Dickson and Hipp (2001)). For any integrable real function $g$ (with non-negative domain), the D-H operator is defined by $\mathcal{T}_{u}g(v) := \int_{v}^{\infty} e^{-u(t-v)} g(t)\mathrm{d}t$ with both $u$ and $v$ non-negative. We begin by looking at the expectation in the following proposition.
\begin{proposition}\normalfont\label{EZ}
For a renewal claim arrival process, the mean of the total discounted IBNR claims until time $t>0$ is given by
\begin{equation}\label{EZdt}
\mathbb{E}[Z_\delta(t)] =\mu_1 e^{-\delta t} \int_{0}^{t} \mathcal{T}_{\delta}w(t-x) \mathrm{d}m(x).
\end{equation}
In particular, when $\delta=0$ this reduces to
\begin{equation}\label{EZdtdelta0}
\mathbb{E}[Z(t)] =\mu_1 \int_{0}^{t} \overline{W}(t-x) \mathrm{d}m(x).
\end{equation}
\noindent\textbf{Proof:}
By taking expectation on (\ref{defZt}) with the use of independence assumptions, one finds
\begin{align*}
\mathbb{E}[Z_\delta(t)]& =\mathbb{E} \bigg[ \sum_{i=1}^{\infty} e^{-\delta (T_{i}+L_{i})} I_{\{T_i \leq t < T_{i}+L_{i}\}} X_{i}\bigg]
 = \sum_{i=1}^{\infty} \mathbb{E} \Big[ e^{-\delta (T_{i}+L_{i})} I_{\{T_i \leq t < T_{i}+L_{i}\}} X_{i}\Big]\\
& = \sum_{i=1}^{\infty} \mu_1 \int_{0}^{t} e^{-\delta x} \mathbb{E} \Big[ e^{-\delta L_{i}}I_{\{L_{i}>t-x\}} \Big] \mathrm{d}F^{\ast(i)}(x) = \mu_1 \int_{0}^{t} e^{-\delta x} \bigg[ \int_{t-x}^{\infty} e^{-\delta y}\mathrm{d}W(y) \bigg] \mathrm{d}\sum_{i=1}^{\infty}F^{\ast(i)}(x) \\
& = \mu_1 \int_{0}^{t} e^{-\delta x} \bigg[\int^\infty_{t-x} e^{-\delta(y-(t-x))}e^{-\delta (t-x)} w(y) \mathrm{d}y \bigg] \mathrm{d}m(x),
\end{align*}
which simplifies to give \eqref{EZdt}. Then \eqref{EZdtdelta0} follows immediately as $\mathcal{T}_0 w(x)=\overline{W}(x)$.
\hfill$\square$
\end{proposition}
This direct technique applied in Proposition \ref{EZ}, however, is not very helpful to derive the expression for higher-order moments. For higher moments, it is useful to employ the conditioning technique based on the first claim arrival time $T_{1}$. It is thus convenient to define the new notation
\[
Z_{2,\delta}(t):=\sum_{i=2}^{N(t)} e^{-\delta (T_{i}-T_{1}+L_{i})} I_{\{ T_{i}+L_{i}>t\}} X_{i}.
\]
Conditioning on $T_{1}=x\leq t$, we note that $Z_{2,\delta}(t)$ is also a discounted IBNR claim that has the same distribution as $Z_\delta(t-x)$. Then $Z_\delta(t)$ may be expressed in terms of $Z_{2,\delta}(t)$ as
\begin{equation}\label{def1Zt}
Z_\delta(t)=e^{-\delta T_{1}} I_{\{ T_1\leq t\} } \Big[ e^{-\delta L_{1}}I_{\{ T_{1}+L_{1}>t\} } X_1 + Z_{2,\delta}(t) \Big].
\end{equation}
Subsequently, recursive formula for the marginal moments of $Z_\delta(t)$, namely $\mathbb{E}[Z_\delta^{n}(t)]$, can be derived for $n \in \mathbb{N}$ via a binomial expansion in Theorem \ref{EZn}. We remark that the same result can also be obtained by considering the moment generating function of $Z_\delta(t)$ followed by differentiation: such a technique was used by Woo (2016) who considered the joint moments of discounted compound renewal sums for $k$ dependent business lines (or $k$ different types of dependent claims) evaluated at the same time point $t$. However, our technique of using a binomial expansion is easily applicable in the present context for the joint moments of the discounted IBNR process (for a single line of business) evaluated at different time points $s$ and $t$, namely $Z_\delta(s)$ and $Z_\delta(t)$, as in the proof of Theorem \ref{EZnZm}.
\begin{theorem}\normalfont\label{EZn}
For $n \in \mathbb{N}$ and $ t>0$, a recursive formula for evaluating $\mathbb{E}[Z_\delta^{n}(t)]$ is given by
\begin{equation}\label{Theorem1Result}
\mathbb{E}[Z_\delta^{n}(t)]=\sum_{i=0}^{n-1} \mu_{n-i} \binom{n}{i} e^{-(n-i)\delta t} \int_{0}^{t} e^{-i\delta x}
\mathcal{T}_{(n-i)\delta}w(t-x) \mathbb{E}[Z_\delta^{i}(t-x)]\mathrm{d}m(x),
\end{equation}
where the starting value is $\mathbb{E}[Z_\delta^{0}(t)]=1$.

\noindent\textbf{Proof:} See Appendix \ref{PT1}.
\hfill$\square$
\end{theorem}


\begin{theorem}\normalfont\label{EZnZm}
For $0<s\leq t$ and $n, m\in \mathbb{N}$, the joint moments of $Z_\delta(s)$ and $Z_\delta(t)$ can be calculated recursively as
\begin{align}\label{Thm2Result}
&\mathbb{E}[Z_\delta^{n}(s)Z_\delta^{m}(t)]= \sum_{i=0}^{n-1} \mu_{n-i} \binom{n}{i} e^{-(n-i)\delta s}\int_{0}^{s} e^{-(m+i)\delta x} \mathbb{E}[Z_\delta^{i}(s-x)Z_\delta^{m}(t-x)]
\mathcal{T}_{(n-i)\delta} w(s-x) \mathrm{d}m(x) \nonumber\\
&~~+\sum_{j=0}^{m-1} \sum_{i=0}^{n} \mu_{n\!+\!m\!-\!i\!-\!j} \binom{n}{i}\binom{m}{j} e^{-(n\!+\!m\!-\!i\!-\!j)\delta t}
\int_{0}^{s} e^{-(i\!+\!j)\delta x} \mathbb{E}[Z_\delta^{i}(s\!-\!x) Z_\delta^{j}(t\!-\!x)] \mathcal{T}_{(n\!+\!m\!-\!i\!-\!j)\delta} w(t\!-\!x) \mathrm{d}m(x),
\end{align}
and this requires the application of Theorem \ref{EZn} as a starting point.

\noindent\textbf{Proof:} See Appendix \ref{PT2}.
\hfill$\square$
\end{theorem}

Furthermore, from Theorems \ref{EZ} and \ref{EZnZm}, the covariance formula for $Z_\delta(s)$ and $Z_\delta(t)$ is readily available as given below.
\begin{corollary}\normalfont\label{EZtZt+h}
For a renewal claim arrival process, the covariance of the total discounted IBNR claims at the time points $s$ and $t$ (where $0<s\leq t$) is given by
\begin{align}\label{Covth}
\mathbb{C}\mathrm{ov}[Z_\delta(s), Z_\delta(t)]
&= \mu_{1} e^{-\delta s} \int_{0}^{s} e^{-\delta x}\mathbb{E}[Z_\delta(t-x)] \mathcal{T}_{\delta}w(s-x) \mathrm{d}m(x) +\mu_{2} e^{-2\delta t} \int_{0}^{s} \mathcal{T}_{2\delta} w(t-x) \mathrm{d}m(x)\nonumber\\
&~~+\mu_{1} e^{-\delta t} \int_{0}^{s} e^{-\delta x} \mathbb{E}[Z_\delta(s-x)]
\mathcal{T}_{\delta} w(t-x) \mathrm{d}m(x) -\mathbb{E}[Z_\delta(s)] \mathbb{E}[Z_\delta(t)].
\end{align}
If there is no discounting (i.e. $\delta=0$), this reduces to
\begin{align}\label{Covd0}
\mathbb{C}\mathrm{ov}[Z(s), Z(t)] &=\mu_1 \int^s_0 \bigg\{ \overline{W}(s-x) \mathbb{E}[Z(t-x)]+\mathbb{E}[Z(s-x)] \overline{W}(t-x)\bigg\}\mathrm{d}m(x)\nonumber\\
&~~+\mu_2 \int^s_0 \overline{W}(t-x)\mathrm{d}m(x)-\mathbb{E}[Z(s)] \mathbb{E}[Z(t)].
\end{align}
Furthermore, for a Poisson process with rate $\lambda>0$ (i.e. $\mathrm{d}m(x)=\lambda \mathrm{d}x$), it is simplified to
\begin{equation}\label{CovPoidzero}
\mathbb{C}\mathrm{ov}[Z(s), Z(t)]=  \lambda \mu_2 \int^s_0 \overline{W}(t-x)dx,
\end{equation}
which is consistent with the covariance expression in Blom et al. (2014, Section 4).

\noindent\textbf{Proof:} As $\mathbb{C}\mathrm{ov}[Z_\delta(s), Z_\delta(t)]=\mathbb{E}[Z_\delta(s)Z_\delta(t)]-\mathbb{E}[Z_\delta(s)] \mathbb{E}[Z_\delta(t)]$, setting $n=m=1$ in Theorem \ref{EZnZm} with a substitution of the result in Theorem \ref{EZ} for $s$ and $t$ yields (\ref{Covth}). When $\delta=0$, we have $\mathcal{T}_0 w(x)=\overline{W}(x)$ and hence (\ref{Covd0}) is found.

For a Poisson arrival process, by substituting $\mathbb{E}[Z(t)]=\lambda\mu_1 \int^t_0 \overline{W}(x)\mathrm{d}x$ followed by interchanging the order of integral in the second term of (\ref{Covd0}), the first two terms can be rewritten as
\begin{align*}
&\mu_1\int^s_0 \bigg\{ \overline{W}(s-x) \mathbb{E}[Z(t-x)]+\mathbb{E}[Z(s-x)] \overline{W}(t-x)\bigg\}\mathrm{d}m(x)\\
&~~=\lambda^2 \mu_1^2 \bigg\{ \int^s_0 \overline{W}(x) \bigg[\int^{t-s+x}_0 \overline{W}(y)\mathrm{d}y+\int^{t}_{t-s+x} \overline{W}(y)\mathrm{d}y\bigg]\mathrm{d}x\bigg\}=\lambda^2  \mu_1^2\bigg[\int^s_0 \overline{W}(x)\mathrm{d}x\bigg]\bigg[\int^{t}_0 \overline{W}(y)\mathrm{d}y\bigg],
\end{align*}
which is identical to the last term in (\ref{Covd0}). Therefore, \eqref{CovPoidzero} follows and the proof is complete.
\hfill$\square$
\end{corollary}

The following corollary is immediately obtainable by letting $s=t$ in Corollary \ref{EZtZt+h}.
\begin{corollary}\normalfont\label{VarZt}
For a renewal claim arrival process, the variance of the total discounted IBNR claims at time $t>0$ is given by
\begin{equation}\label{Vartth}
\mathbb{V}\mathrm{ar}[Z_\delta(t)] =2 \mu_{1} e^{-\delta t} \int_{0}^{t} e^{-\delta x}\mathbb{E}[Z_\delta(t-x)] \mathcal{T}_{\delta}w(t-x)\mathrm{d}m(x) +\mu_{2} e^{-2\delta t} \int_{0}^{t} \mathcal{T}_{2\delta} w(t-x) \mathrm{d}m(x) -\{\mathbb{E}[Z_\delta(t)]\}^2.
\end{equation}
When $\delta=0$, this reduces to
\begin{equation}\label{VARdzero}
\mathbb{V}\mathrm{ar}[Z(t)]=2 \mu_{1} \int^t_0 \overline{W}(t-x) \mathbb{E}[Z(t-x)] \mathrm{d}m(x)+\mu_2 \int^t_0 \overline{W}(t-x)\mathrm{d}m(x)-\{\mathbb{E}[Z(t)]\}^2.
\end{equation}
Furthermore, for a Poisson arrival process with rate $\lambda>0$, it is simplified to

\begin{equation}\label{VARdzeroPoi}
\mathbb{V}\mathrm{ar}[Z(t)]= \lambda \mu_2 \int^t_0 \overline{W}(x)\mathrm{d}x.
\end{equation}

Note that $\mathbb{V}\mathrm{ar}[Z(t)]\rightarrow \lambda \mu_2 \mathbb{E}[L]$ as $t\rightarrow \infty$, where $\mathbb{E}[L]=\int^\infty_0 \overline{W}(x)\mathrm{d}x$.
\end{corollary}

\begin{theorem}\label{CovZsZI}\normalfont When $\delta=0$ and $\mathbb{E}[L]<\infty$, the limit of $\mathbb{C}\mathrm{ov}[Z(s), Z(t)]$ is given by
\begin{equation}\label{Theorem3result}
	\mathbb{C}\mathrm{ov}[Z(s), Z(t)]
	\longrightarrow 0 \text{~~as~~} t\rightarrow \infty.
\end{equation}
			\noindent \textbf{Proof:} First, we use the result $\mathbb{E} [Z(\infty)]:= \lim_{t\to \infty}\mathbb{E} [Z(t)] =\mu_1 \mathbb{E}[L] / \mathbb{E} [\tau]$ which is obtainable by applying the Smith's Renewal Theorem to \eqref{EZdtdelta0} as
\begin{equation*}
	\mathbb{E} [Z(t)] \longrightarrow \frac{\mu_1}{\mathbb{E}[\tau]} \int^\infty_0 \overline{W}(y)\mathrm{d}y\text{~~as~~}t\rightarrow \infty.
\end{equation*}
By taking the limit $t\rightarrow \infty$ with the help of dominated convergence, the first term on the right-hand side of (\ref{Covd0}) in Corollary \ref{EZtZt+h} becomes
\[\lim\limits_{t\rightarrow \infty} \mu_1 \int^s_0 \bigg\{ \overline{W}(s-x) \mathbb{E}[Z(t-x)]+\mathbb{E}[Z(s-x)] \overline{W}(t-x)\bigg\}\mathrm{d}m(x) = \frac{\mu_1^2 \mathbb{E}[L]}{\mathbb{E}[\tau]}\int_{0}^{s}\overline{W}(s-x)\mathrm{d}m(x),
\]
where the second term in the integration converges to $0$. This can be cancelled out with the third term in (\ref{Covd0}). Then, the desired result \eqref{Theorem3result} follows. 	\hfill$\square$
\end{theorem}

The result in Theorem \ref{CovZsZI} shows that, provided that the delay time $L$ is integrable, the IBNR claims are asymptotically uncorrelated, i.e. on a first approximation, one has $\mathbb{E} [Z(t)Z(s)] \approx \mathbb{E} [Z(t)]\mathbb{E} [Z(s)] $ as $t\to \infty$. More information on the speed of convergence of the covariance towards $0$ will be given in the forthcoming Propositions \ref{CovExpPro} and \ref{CovPareto} when additional distributional assumptions are made. In fact, one of the obstacles in the proofs in those propositions is to be able to provide a second order approximation for the joint moment $\mathbb{E} [Z(t)Z(s)]$ as $t\to\infty$, in order to get some information on the rate of decrease of $\mathbb{C}\mathrm{ov}[Z(s), Z(t)]$.\vspace{-.1in}

\section{Fractional Poisson Process}\label{FPP}\vspace{-.05in}
In this section, we consider the fractional Poisson process for the claim counting process $\{N(t)\}_{t\ge0}$, so that the renewal function and renewal density are given by, for $ t>0$,
\begin{equation}\label{FP}
m(t)=\mathbb{E}[N(t)]= \frac{\lambda t^{\alpha}}{\Gamma ( 1+\alpha )}\text{~~and~~}m'(t)=\frac{\lambda \alpha t^{\alpha -1}}{\Gamma ( 1+\alpha )}=\frac{\lambda t^{\alpha-1}}{\Gamma(\alpha)},
\end{equation}
where $\Gamma(x)$ is the (ordinary) Gamma function defined by $\Gamma(x)=\int_{0}^{\infty}y^{x-1}e^{-y}\mathrm{d}y$ for $\mathrm{Re}(x)>0$, and $\alpha$ and $\lambda$ are parameters such that $0<\alpha\leq 1$ and $\lambda>0$. From Equation (28) of Laskin (2003), its variance is given by
\begin{equation*}
\mathbb{V}\mathrm{ar}[N(t)]=m(t)+[m(t)]^2 \left[\frac{\alpha \mathrm{B}(\alpha, 1/2)}{2^{2\alpha-1}}-1\right],
\end{equation*}
where $\mathrm{B}(a,b)=\Gamma(a)\Gamma(b)/\Gamma(a+b)=\int^1_0 y^{a-1}(1-y)^{b-1}\mathrm{d}y$ is the Beta function. Therefore, the variance is increasing in time according to a power law $t^{2\alpha}$ when $0<\alpha<1$. If $\alpha=1$, this process becomes the ordinary Poisson process with rate $\lambda$. From Repin and Saichev (2000), the asymptotic behaviour for the survival function of the interarrival times satisfies (when $0<\alpha<1$)
\begin{equation}\label{asymptotics_interclaim_frac}
1-F(t) \sim Ct^{-\alpha}\text{~~as~~}t\rightarrow \infty,
\end{equation}
where $C$ is some constant. For the estimation of the parameters $\alpha$ and $\lambda$ of a fractional Poisson process, interested readers are referred to Cahoy et al. (2010), who first derived the estimators by matching the first two moments of the log of the interarrival time $\tau$ and then proved their asymptotic normality.

In the sequel, two different distributional assumptions for the random delay variable $L$, namely exponential and Pareto, will be considered in detail. To proceed with our analysis, we first recall the definitions and properties of some special functions.
\begin{enumerate}
	\item[(i)] The Gaussian (or ordinary) hypergeometric function ${}_{2}\mathrm{F}_1(a,b,c,z)$ is defined as the power series (for $|z|\le 1$)
	\begin{equation}\label{ohf}	 {}_{2}\mathrm{F}_1(a,b,c,z)=\sum_{n=0}^{\infty}\frac{(a)_n ( b)_n}{(c)_n}\frac{z^n}{n!},
	\end{equation}
where $(a)_n=a(a+1)\cdots(a+n-1)$ for $n\in\mathbb{N}$ is the Pochhammer symbol with $(a)_0=1$. Its integral representation is given by
\begin{equation}\label{ohf1}
{}_{2}\mathrm{F}_1(a,b,c,z) = \frac{1}{\mathrm{B}(b,c-b)}\int_{0}^{1}x^{b-1}(1-x)^{c-b-1}(1-zx)^{-a}\mathrm{d}x, \quad \mathrm{Re}(c)>\mathrm{Re}(b)>0,
\end{equation}
provided that the right-hand side converges. We also have the relationship between hypergeometric function and incomplete Beta function given by
\begin{equation}\label{ibf}
\mathrm{B}(a,b,x)=\int^x_0 t^{a-1}(1-t)^{b-1}\mathrm{d}t={\frac {x^{a}}{a}}{}_{2}\mathrm{F}_{1}(a,1-b;a+1;x).
\end{equation}
Another important identity for the Gaussian  hypergeometric function is
\begin{equation}\label{ohf2}
{}_{2}\mathrm{F}_{1}(a,b;c;1)={\frac {\Gamma (c)\Gamma (c-a-b)}{\Gamma (c-a)\Gamma (c-b)}},\qquad \mathrm{Re}(c)>\mathrm{Re}(a+b).
\end{equation}

\item[(ii)] The confluent hypergeometric function of the first kind or the Kummer's function ${}_{1}\mathrm{F}_1(a,b,z)$ introduced by Kummer (1837) (see e.g. Abramowitz and Stegun (1972, Chapter 13)) is defined by
\begin{equation}\label{Ku0}
{}_{1} \mathrm{F}_1(a,b,z) = \sum_{n=0}^{\infty} \frac{(a)_n}{(b)_n}\frac{z^n}{n!},
\end{equation}
and its integral representation is
\begin{equation}\label{Ku}
{}_{1}\mathrm{F}_1(a,b,z) = \frac{1}{\mathrm{B}(a,b-a)}\int_{0}^{1} e^{zx}x^{a-1}(1-x)^{b-a-1}\mathrm{d}x,\qquad \mathrm{Re}(b)>\mathrm{Re}(a)>0.
\end{equation}
We also have the transform relation between two Kummer's functions as
\begin{equation}\label{Kut}
{}_{1}\mathrm{F}_1(a,b,z) = e^{z}{}_{1}\mathrm{F}_1(b-a,b,-z).
\end{equation}
The Kummer's function satisfies
\begin{equation}\label{asymptotics1_Kummer}
{}_1\mathrm{F}_1(a,b,z) =\frac{\Gamma(b)}{\Gamma(a)}\,e^zz^{a-b}\{ 1+O(|z|^{-1})\},\qquad \mathrm{Re}(z)>0.
\end{equation}
See Equation (13.1.4) on p.504 of Abramowitz and Stegun (1972).

\end{enumerate}\vspace{-.2in}
	
\subsection{Exponential delay}
In this section, it is assumed that $\overline{W}(x)=e^{-\beta x}$ for $x\ge0$, where $\beta>0$ is the exponential parameter for the delay $L$ so that $\mathbb{E}[L]=1/\beta$ and $\mathbb{E}[e^{-sL}]=\beta/(\beta+s)$.\\

\noindent{\bf Expectation of total discounted IBNR process.} We first establish an expression for $\mathbb{E}[Z_\delta(t)]$ and provide its asymptotic behaviour as $t\to \infty$.
\begin{proposition}\label{Expectation_expo_case}\normalfont
With exponential delay in a fractional Poisson claim arrival process, the mean of the total discounted IBNR claims $Z_\delta(t)$ when $t>0$ is given by, for $\delta=0$,
\begin{equation}\label{EZdzeroFP}
\mathbb{E}[Z(t)]=\frac{\mu_1\lambda e^{-\beta t}t^\alpha }{\Gamma(1+\alpha)}{}_1 \mathrm{F}_1
(\alpha,1+\alpha,\beta t)=\mu_1 e^{-\beta t} {}_1\mathrm{F}_1
(\alpha,1+\alpha,\beta t) m(t),
\end{equation}
and for $\delta>0$,
\begin{equation}\label{EZdFP}
\mathbb{E}[Z_\delta(t)]=\frac{\mu_1\lambda  \mathbb{E}[e^{-\delta L}] e^{-(\beta+\delta) t}t^\alpha }{\Gamma(1+\alpha)}{}_1 \mathrm{F}_1
(\alpha,1+\alpha,\beta t)=\mu_1 \mathbb{E}[e^{-\delta L}] e^{-(\beta+\delta) t} {}_1 \mathrm{F}_1
(\alpha,1+\alpha,\beta t) m(t),
\end{equation}
\noindent where $m(t)$ is the renewal function given in (\ref{FP}) for the fractional Poisson process. In addition, the corresponding asymptotic behaviours are
\begin{equation}\label{asyEZdzeroFP}
\mathbb{E}[Z(t)] \sim \frac{\mu_1\lambda}
{\beta \Gamma(\alpha)} t^{\alpha-1}=\mu_1 \mathbb{E}[L] m'(t) \text{~~as~~}t\rightarrow \infty,
 \end{equation}
 and
\begin{equation}\label{asyEZdFP}
 \mathbb{E}[Z_\delta(t)] \sim \frac{\mu_1\lambda \mathbb{E}[e^{-\delta L}]}{\beta \Gamma(\alpha)}\,  e^{-\delta t}t^{\alpha-1}=\mu_1 \mathbb{E}[L] \mathbb{E}[e^{-\delta L}] e^{-\delta t} m'(t)\text{~~as~~} t\rightarrow \infty,
 \end{equation}
respectively.

\noindent \textbf{Proof:} First, we note that for exponential delay, one has $\mathcal{T}_\delta w(x)=\frac{\beta}{\beta+\delta}e^{-\beta x}=\mathbb{E}[e^{-\delta L}] e^{-\beta x}=\mathbb{E} [e^{-\delta L}] \overline{W}(x)$. In particular, $\mathcal{T}_0 w(x)=\overline{W}(x)$  when $\delta=0$. Hence, $\mathcal{T}_\delta w(x)= \mathbb{E} [e^{-\delta L} ]\mathcal{T}_0 w(x)$, and as a result of Proposition \ref{EZ} we have the relation
\begin{equation}\label{EZEdZ}
\mathbb{E}[Z_\delta(t)]=e^{-\delta t}\mathbb{E}[e^{-\delta L}]\mathbb{E}[Z(t)].
\end{equation}
In other words, it is immediate to obtain the mean in the presence of discounting from the mean without discounting.

When $\delta=0$, we have from Proposition \ref{EZ} that $\mathbb{E}[Z(t)] =\mu_1 e^{-\beta t}\int^t_0 e^{\beta x}\mathrm{d}m(x)$. Then, with the fractional Poisson process for claim arrivals, substitution of its renewal density (\ref{FP}) followed by a change of variable $y=x/t$ yields
\begin{equation}\label{expression_expectation}
\mathbb{E}[Z(t)]= \frac{\mu_1 \lambda \alpha  e^{-\beta t} }{\Gamma(1+\alpha)}\int^t_0 e^{\beta x}x^{\alpha-1} \mathrm{d}x=
\frac{\mu_1\lambda \alpha e^{- \beta  t} }{\Gamma(1+\alpha)} \bigg[t^\alpha \int^1_0 e^{(\beta t)y} y^{\alpha-1}\mathrm{d}y\bigg].
\end{equation}
Utilizing the Kummer's function in (\ref{Ku}), the above expectation may be neatly expressed as (\ref{EZdzeroFP}).
Consequently, using the relation (\ref{EZEdZ}), we easily obtain (\ref{EZdFP}).

It is indeed useful to have expressed $\mathbb{E}[Z(t)]$ and $\mathbb{E}[Z_\delta(t)]$ in terms of ${}_1\mathrm{F}_1(a,b,z)$ as its asymptotic behaviour is well established in the literature. More precisely, with the help of (\ref{asymptotics1_Kummer}), one finds (\ref{asyEZdzeroFP}) and (\ref{asyEZdFP}). \hfill$\square$
\end{proposition}

\begin{example}\normalfont (Poisson process) For notational convenience, we denote the present value of a $t$-year continuous annuity payable at rate \$1 per year as $\bar{a}_{\lcroof{t}\,\delta} :=\int^t_0 e^{-\delta x}\mathrm{d}x= (1-e^{-\delta t})/\delta$. When the claims arrive according to a Poisson process with rate $\lambda$, from Proposition \ref{EZ} we have $\mathbb{E}[Z(t)]=\mu_1 \lambda \bar{a}_{\lcroof{t}\,\beta}$ and thus $\mathbb{E}[Z_\delta(t)]=\mu_1 \lambda e^{-\delta t} \mathbb{E}[e^{-\delta L}] \bar{a}_{\lcroof{t}\,\beta}$ from (\ref{EZEdZ}). Of course, the same results are obtainable by setting $\alpha=1$ in Proposition \ref{Expectation_expo_case}, since when $\alpha =1 $ the fractional Poisson process is just a Poisson process. In particular, noting that $e^{-\beta t} t ~{}_1 \mathrm{F}_1(1,2,\beta t)=\bar{a}_{\lcroof{t}\,\beta}$ from (\ref{Ku0}), the result follows from (\ref{EZdzeroFP}) for $\alpha=1$. Also, asymptotic behaviours are given as
 \[
\mathbb{E}[Z(t)] \sim \mu_1\lambda  \mathbb{E}[{L}]\text{~~and~~} \mathbb{E}[Z_\delta(t)] \sim \mu_1\lambda \mathbb{E}[L] \mathbb{E}[e^{-\delta L}]\,  e^{-\delta t}\text{~~as~~}t\rightarrow \infty.
\]
\end{example}
\vspace{-.28in}
\hfill$\square$\vspace{.1in}

\noindent{\bf Variance of total discounted IBNR process.} Next, from (\ref{Vartth}), we find the variance of $Z_\delta(t)$ and its asymptotic result in the following proposition.
\begin{proposition}\normalfont\label{P3}
With exponential delay in a fractional Poisson claim arrival process, the variance of the total discounted IBNR claims $Z_\delta(t)$ when $t>0$ is given by
\begin{align}\label{VARZdt}
\mathbb{V}\mathrm{ar}[Z_\delta(t)] &=\frac{2\mu_1^2 \lambda^2  \alpha \{\mathbb{E}[e^{-\delta L}]\}^2 e^{-2(\beta+\delta) t}}{[\Gamma(1+\alpha)]^2} \sum^\infty_{n=0} A_n t^{2\alpha+n}  {}_1 \mathrm{F}_1 (\alpha, 2\alpha+n+1,2\beta t)\nonumber\\
&~~+\frac{\mu_2 \lambda\mathbb{E}[e^{-2\delta L}]  e^{-(\beta+2\delta) t} t^\alpha}{\Gamma(1+\alpha)}{}_1 \mathrm{F}_1(\alpha,1+\alpha,\beta t)
 - \frac{\mu_1^2 \lambda^2  \{\mathbb{E}[e^{-\delta L}]\}^2 e^{-2(\beta+\delta) t} t^{2\alpha}}{[\Gamma(1+\alpha)]^2}[{}_1 \mathrm{F}_1(\alpha,1+\alpha,\beta t)]^2,
\end{align}
where
\begin{equation}\label{An}
A_n=\dfrac{\alpha \beta^n }{n!(\alpha+n) }\mathrm{B}(\alpha,\alpha+n+1).
\end{equation}
In addition, the asymptotic behaviour when $0<\alpha<1$ is
\begin{equation}\label{VARZdtasy}
\mathbb{V}\mathrm{ar}[Z_\delta(t)] \sim \left(\frac{\mu_1^2 \lambda^2 \{\mathbb{E}[e^{-\delta L}]\}^2 }{\beta^{\alpha+1} \Gamma(\alpha)}+\frac{\mu_2 \lambda \mathbb{E}[e^{-2\delta L}] }{\beta \Gamma(\alpha)}\right) e^{-2 \delta t} t^{\alpha-1}\text{~~as~~}t\rightarrow \infty.
\end{equation}

\noindent \textbf{Proof:} The three terms on the right-hand side of (\ref{Vartth}) are expressed as follows. The first term can be calculated by using (\ref{EZdFP}), $\mathcal{T}_\delta w(x)=\frac{\beta}{\beta+\delta}e^{-\beta x}=\mathbb{E}[e^{-\delta L}] e^{-\beta x}$ and (\ref{FP}) as
\begin{align}\label{f1}
&2 \mu_{1} e^{-\delta t} \int_{0}^{t} e^{-\delta x}\mathbb{E}[Z_\delta(t-x)] \mathcal{T}_{\delta}w(t-x)\mathrm{d}m(x)\nonumber\\
&~~=\frac{2\mu_1^2 \lambda^2  \alpha   \{\mathbb{E}[e^{-\delta L}]\}^2 }{[\Gamma(1+\alpha)]^2}  e^{-2\delta t}\int^t_0 e^{-2\beta x}x^\alpha (t-x)^{\alpha-1}
{}_1 \mathrm{F}_1(\alpha,1+\alpha,\beta x)  \mathrm{d}x.
\end{align}
Using the infinite series representation of the Kummer's function in (\ref{Ku0}) with a change of variable from $x/t$ to $y$, the integral term above can be expressed as
\begin{align}\label{ftermVarEP1}
& \int^t_0 e^{-2\beta x}x^\alpha (t-x)^{\alpha-1}
{}_1 \mathrm{F}_1(\alpha,1+\alpha,\beta x)  \mathrm{d}x  =\sum^\infty_{n=0}\frac{\alpha \beta^n}{n!(\alpha+n)} t^{2\alpha+n} \int^1_0  e^{-2\beta t y} y^{\alpha+n} (1-y)^{\alpha-1} \mathrm{d}y\nonumber\\
&~~= \sum^\infty_{n=0} A_n t^{2\alpha+n}{}_1\mathrm{F}_1
(\alpha+n+1,2\alpha+n+1,-2\beta t) =\sum^\infty_{n=0} A_n t^{2\alpha+n} e^{-2\beta t} {}_1\mathrm{F}_1 (\alpha, 2\alpha+n+1,2\beta t),
\end{align}
where $A_n$ is defined in (\ref{An}) and the last equality is due to the Kummer's transformation (\ref{Kut}). From (\ref{EZdt}) together with (\ref{FP}) and (\ref{EZdFP}), the sum of the remaining two terms in (\ref{Vartth}) is given by
			\begin{align}\label{st}
	&		 \frac{\mu_2\mathbb{E}[Z_{2\delta}(t)]}{\mu_1}-\{\mathbb{E}[Z_\delta(t)]\}^2\nonumber\\
	&~~=
	\frac{\mu_2\lambda  \mathbb{E}[e^{-2\delta L}] e^{-(\beta+2\delta) t}t^\alpha }{\Gamma(1+\alpha)}{}_1 \mathrm{F}_1
(\alpha,1+\alpha,\beta t)-\frac{\mu_1^2\lambda^2  (\mathbb{E}[e^{-\delta L}])^2 e^{-2(\beta+\delta) t}t^{2\alpha} }{[\Gamma(1+\alpha)]^2}[{}_1 \mathrm{F}_1
(\alpha,1+\alpha,\beta t)]^2.
			\end{align}
Hence, combining \eqref{f1}-\eqref{st} results in (\ref{VARZdt}).

To obtain the asymptotic formula \eqref{VARZdtasy} as $t\rightarrow \infty$, we focus on the  first term of $\mathbb{V}\mathrm{ar}[Z_\delta(t)]$ in \eqref{VARZdt} (or \eqref{ftermVarEP1}) and would like to show that
\begin{equation}\label{limitsum}
\lim_{t\rightarrow\infty} \frac{\sum^\infty_{n=0} A_n t^{2\alpha+n} e^{-2\beta t} {}_1\mathrm{F}_1 (\alpha,2\alpha+n+1,2\beta t)}{t^{\alpha-1}}=C
\end{equation}
for some constant $C$. Note that \eqref{asymptotics1_Kummer} implies
\begin{equation*}
\lim_{z\rightarrow\infty} \frac{_1\mathrm{F}_1(a,b,z)}{\frac{\Gamma(b)}{\Gamma(a)}\,e^zz^{a-b}} =1.
\end{equation*}
In order to apply such an asymptotic result, we need to check the validity of interchanging limit and infinite summation on the left-hand side of \eqref{limitsum}. It can be proved that the related sequence is uniformly bounded (which gives a sufficient condition). This part of the proof is not trivial and is provided in Part 1 of the `Supplementary materials'. We can now evaluate the limit on the left-hand side of \eqref{limitsum} by interchanging the order of limit and infinite summation followed by substitution of \eqref{An}, and this gives rise to
\begin{align*}
&\lim_{t\rightarrow\infty} \frac{\sum^\infty_{n=0} A_n t^{2\alpha+n} e^{-2\beta t} {}_1\mathrm{F}_1 (\alpha,2\alpha+n+1,2\beta t)}{t^{\alpha-1}} =\sum^\infty_{n=0}A_n \left[\lim_{t\rightarrow\infty} t^{\alpha+n+1} e^{-2\beta t} {}_1\mathrm{F}_1 (\alpha,2\alpha+n+1,2\beta t)\right]\\
=&\sum^\infty_{n=0} \frac{\alpha}{\alpha+n}\,\frac{\beta^n }{n!}\,\frac{\Gamma(\alpha)\Gamma(\alpha+n+1)}{\Gamma(2\alpha+n+1)} \left[\frac{1}{(2\beta)^{\alpha+n+1}}\,\frac{\Gamma(2\alpha+n+1)}{\Gamma(\alpha)} \lim_{t\rightarrow\infty}\frac{_1\mathrm{F}_1 (\alpha,2\alpha+n+1,2\beta t)}{\frac{\Gamma(2\alpha+n+1)}{\Gamma(\alpha)}\,e^{2\beta t}(2\beta t)^{-\alpha-n-1}}\right]\\
=&~\frac{\alpha}{(2\beta)^{\alpha+1}} \sum^\infty_{n=0} \frac{\Gamma(\alpha+n)}{n!\,2^n} =\frac{\alpha}{(2\beta)^{\alpha+1}} \sum^\infty_{n=0}\frac{\int_0^\infty y^{\alpha+n-1}e^{-y}\,\mathrm{d}y}{n!\,2^n} =\frac{\alpha}{(2\beta)^{\alpha+1}} \int_0^\infty y^{\alpha-1} \left[\sum^\infty_{n=0}\frac{(y/2)^n}{n!}\right]e^{-y}\,\mathrm{d}y\\
=&~\frac{\alpha}{(2\beta)^{\alpha+1}} \int_0^\infty y^{\alpha-1}e^{-\frac{y}{2}}\,\mathrm{d}y=\frac{\alpha}{(2\beta)^{\alpha+1}}\,2^{\alpha}\Gamma(\alpha) =\frac{\Gamma(\alpha+1)}{2\beta^{\alpha+1}}.
\end{align*}
Hence, we identify $C:=\frac{\Gamma(\alpha+1)}{2\beta^{\alpha+1}}$ in (\ref{limitsum}) so that (\ref{ftermVarEP1}) asymptotically behaves like $C t^{\alpha-1}$. Combining (\ref{f1}) and (\ref{st}) with asymptotic results of (\ref{ftermVarEP1}) and (\ref{asyEZdFP}) yields (\ref{VARZdtasy}).				 
					\hfill$\square$
\end{proposition}

\noindent{\bf Covariance of total discounted IBNR processes.} We finish this section by studying the covariance of $Z_\delta(s)$ and $Z_\delta(t)$ and its asymptotic behaviour as $t\to \infty$.
\begin{proposition}\normalfont\label{CovExpPro}
With exponential delay in a fractional Poisson claim arrival process, the covariance of the total discounted IBNR claims $Z_\delta(s)$ and $Z_\delta(t)$ for $0<s\leq t$ is given by
\begin{align}\label{CovExp}
&\mathbb{C}\mathrm{ov}[Z_\delta(s),Z_\delta(t)] = \frac{\mu_1^2\lambda^2 \alpha \{\mathbb{E}[e^{-\delta L}]\}^2}{[\Gamma(1+\alpha)]^2} e^{-(\beta+\delta)(s+t)} \sum_{n=0}^{\infty}A_n\Big[  s^{2\alpha+n}  {}_1\mathrm{F}_1 (\alpha, 2\alpha\!+\!n\!+\!1,2\beta s)+s^\alpha \mathrm{W}(2\beta s, n,t) \Big]\nonumber\\
&+ \frac{\mu_2\lambda \mathbb{E}[e^{-2\delta L}]}{\Gamma(1+\alpha)} e^{-(\beta+2\delta) t}s^\alpha {}_1\mathrm{F}_1 (\alpha, 1\!+\!\alpha,\beta s)- \frac{\mu_1^2\lambda^2 \{\mathbb{E}[e^{-\delta L}]\}^2}{[\Gamma(1+\alpha)]^2} e^{-(\beta+\delta) (s+t)}(st)^\alpha  {}_1 \mathrm{F}_1
(\alpha,1\!+\!\alpha,\beta s){}_1 \mathrm{F}_1
(\alpha,1\!+\!\alpha,\beta t),
\end{align}
where
\begin{equation}\label{Wsnt}
\mathrm{W}(2\beta s,n,t) :=\frac{1}{\mathrm{B}(\alpha,\alpha+n+1)}\int_{0}^{1}e^{2\beta s y}y^{\alpha-1}(t- sy)^{\alpha+n}\mathrm{d}y.
\end{equation}
In addition, the asymptotic behaviour when $0<\alpha<1$ is
\begin{equation}\label{COVZdstasy}
\mathbb{C}\mathrm{ov}[Z_\delta(s),Z_\delta(t)] \sim \left(\frac{\mu_1^2 \lambda(1-\alpha) \{\mathbb{E}[e^{-\delta L}]\}^2 }{\beta \Gamma(\alpha)} e^{-(\beta+\delta) s}\left[ \int_{0}^{s} e^{\beta x}x\mathrm{d}m(x)\right]\right) e^{-\delta t} t^{\alpha-2}\text{~~as~~}t\rightarrow \infty.
\end{equation}

\noindent \textbf{Proof:} See Appendix \ref{PCovExpPro}.
\hfill$\square$
\end{proposition}


\subsection{Pareto delay}
This section considers a heavy-tailed distribution for random delay. In particular, it is assumed that $L$ follows a Pareto distribution with tail
\begin{equation}\label{Ptail}
\overline{W}(x) = 1- W(x) = \Big(\dfrac{\theta}{\theta + x}\Big)^\eta,\qquad x\ge0,
\end{equation}
where $\eta,\theta>0$ are the parameters. In this case, $\mathbb{E}[L]=\frac{\theta}{\eta-1}$ when $\eta>1$. However, if $0<\eta\le 1$ then $L$ has infinite mean, and consequently $\overline{W}(x)$ is not a directly Riemann integrable function and the Smith's Renewal Theorem is no longer applicable. In the rest of this section, it is assumed that $\delta=0$. The proofs of all three propositions require the next lemma.
\begin{lemma}\normalfont\label{L3}
Define, for $y\in[0,1)$,
\begin{equation}\label{G1ydef}
G^*(y):=\int_{0}^{y}(1-x)^{-\eta}x^{\alpha -1}\mathrm{d}x,
\end{equation}
\noindent where $\eta\geq 1$ and $\alpha\in (0,1)$. The above integral satisfies, as $y\rightarrow1^-$,
\begin{equation}\label{G1yasym}
G^*(y) =\begin{cases}
\frac{1}{\eta-1}(1-y)^{-\eta+1}+ o\left( (1-y)^{-\eta +1}\right),& \text{if }\eta >1,\\
-\ln(1-y)+C+o(1), & \text{if } \eta=1,
\end{cases}
\end{equation}
for some constant $C$. In particular, when $1<\eta\le 2$, finer asymptotic results can be given by, as $y\rightarrow1^-$,
\begin{equation}\label{G1yasymfine}
G^*(y) =\begin{cases}
\frac{1}{\eta-1}(1-y)^{-\eta+1}+C^*+o(1),& \text{if }1<\eta<2,\\
(1-y)^{-1}-(1-\alpha)\ln(1-y)+C^{**}+o(1), & \text{if } \eta=2,
\end{cases}
\end{equation}
for some constants $C^*$ and $C^{**}$. 

\noindent \textbf{Proof:} See Part 2 of the `Supplementary materials'. 
\hfill$\square$
\end{lemma}

\noindent{\bf Expectation of total IBNR process.} We have the following proposition for the mean of $Z(t)$.
\begin{proposition}\normalfont \label{EZParPro} With Pareto delay in a fractional Poisson claim arrival process, the mean of the total IBNR claims $Z(t)$ when $t>0$ is given by
\begin{equation}\label{expr_EZ_phi}
\mathbb{E}[Z(t)] = \dfrac{\mu_1 \lambda \theta^{\eta}}{\Gamma(\alpha)(\theta+t)^{\eta-\alpha}}  \int_{0}^{\frac{t}{\theta+t}}(1- y)^{-\eta}y^{\alpha-1}\mathrm{d}y.
\end{equation}
Its asymptotic behaviour is, as $t\rightarrow \infty$,
\begin{equation}\label{PEZ2}
\mathbb{E}[Z(t)] \sim
\begin{cases}
\dfrac{\mu_1 \lambda \theta^\eta \Gamma(1-\eta)}{\Gamma(\alpha+1-\eta)} t^{\alpha-\eta}, & \mathrm{if}~~0<\eta<1,\vspace{.07in}\\
\dfrac{\mu_1 \lambda \theta}{\Gamma(\alpha)}t^{\alpha-1}\ln t, & \mathrm{if}~~\eta=1,\vspace{.07in}\\
\dfrac{\mu_1 \lambda \theta}{(\eta-1) \Gamma(\alpha)} t^{\alpha-1}, & \mathrm{if}~~\eta>1.
\end{cases}
\end{equation}

\noindent \textbf{Proof:} Substitution of (\ref{FP}) and (\ref{Ptail}) into (\ref{EZdtdelta0}) yields 
\begin{equation*}
\mathbb{E}[Z(t)] = \mu_1\int_{0}^{t} \Big(\frac{\theta}{\theta + t-x}\Big)^\eta \frac{\lambda x^{\alpha-1}}{\Gamma(\alpha)}\mathrm{d}x = \frac{\mu_1 \lambda}{\Gamma(\alpha)} \Big(\frac{\theta}{\theta+t}\Big)^\eta\int_{0}^{t}\Big(1- \frac{x}{\theta+t} \Big)^{-\eta}x^{\alpha-1}\mathrm{d}x,
\end{equation*}
from which \eqref{expr_EZ_phi} follows by a change of variable from $x/(\theta+t)$ to $y$.

Regarding the asymptotic behaviour, we separate the analysis into three cases as follows.

\noindent\textbf{Case 1.} $0<\eta<1$: The integral in \eqref{expr_EZ_phi} can be expressed in terms of the incomplete Beta function (\ref{ibf}), i.e.
\begin{equation*}
\mathbb{E}[Z(t)]= \dfrac{\mu_1 \lambda \theta^{\eta}}{\Gamma(\alpha)(\theta+t)^{\eta-\alpha}}
\mathrm{B}\bigg(\alpha,1-\eta, \frac{t}{\theta+t} \bigg).
\end{equation*}
As $\lim_{t\rightarrow\infty}\mathrm{B}(\alpha,1-\eta, \frac{t}{\theta+t}) = \mathrm{B}(\alpha,1-\eta) = \frac{\Gamma(\alpha)\Gamma(1-\eta)}{\Gamma(\alpha+1-\eta)}$, one finds the first asymptotic result in (\ref{PEZ2}).

\noindent\textbf{Case 2.} $\eta=1$: Replacing $y$ by $t/(\theta+t)$ in Lemma \ref{L3} and using the second result of \eqref{G1yasym} gives
\[
\int_{0}^{\frac{t}{\theta+t}}(1- y)^{-1}y^{\alpha-1}\mathrm{d}y =-\ln\left(\frac{\theta}{\theta+t}\right)+C+o(1) \sim \ln t \text{~~as~~}t\rightarrow \infty.
\]
Moreover, it is obvious that $1/(\theta+t)^{1-\alpha}\sim t^{\alpha-1}$ as $t\rightarrow \infty$. By substituting the these asymptotics into (\ref{expr_EZ_phi}), we arrive at the second result in (\ref{PEZ2}).

\noindent\textbf{Case 3.} $\eta>1$: Similar to Case 2, applying the first result of \eqref{G1yasym} yields
\[
\int_{0}^{\frac{t}{\theta+t}}(1- y)^{-\eta}y^{\alpha-1}\mathrm{d}y =\frac{1}{\eta-1}\left(\frac{\theta}{\theta+t}\right)^{-\eta+1}+ o\bigg( \left(\frac{\theta}{\theta+t}\right)^{-\eta +1}\bigg) \sim\frac{\theta^{-\eta+1}}{\eta-1} t^{\eta-1}\text{~~as~~}t\rightarrow \infty.
\]
\noindent Combining with the fact that $1/(\theta+t)^{\eta-\alpha}\sim t^{\alpha-\eta}$ as $t\rightarrow \infty$ gives the third asymptotic formula of (\ref{PEZ2}) thanks to (\ref{expr_EZ_phi}).\hfill$\square$


\end{proposition}

\begin{remark}\label{RemarkParetomean}\normalfont 
Consider the fractional Poisson claim arrival process with $0<\alpha<1$. One particular consequence of the asymptotic results (\ref{PEZ2}) in Proposition \ref{EZParPro} is that
\[
\lim_{t\rightarrow \infty}\mathbb{E}[Z(t)] =
\begin{cases}
 0, &\mathrm{if}~~0<\alpha<\eta,
\\
\mu_1\lambda \theta^\eta \Gamma(1-\eta), &\mathrm{if}~~0<\eta=\alpha<1,\\
\infty, & \mathrm{if}~~0<\eta<\alpha<1.
\end{cases}
\]
The above result can be intuitively understood in the following way. Since the interarrival time $\tau$ and the delay $L$ exhibit power decay according to (\ref{asymptotics_interclaim_frac}) and (\ref{Ptail}), one checks easily that
$$
\tau \le_{\small \mbox{st},\infty} L\ (\mbox{resp. }L \le_{\small \mbox{st},\infty} \tau ) \mbox{~ when~~}\alpha> \eta\ (\mbox{resp. when }\eta>\alpha),
$$
Hence, if $ \alpha>\eta$ then delays are too long and cannot make up for the faster arriving claims, which justifies why the expected value of IBNR claims tends to infinity as time goes by. \hfill$\square$
\end{remark}

\noindent{\bf Variance of total IBNR process.} The upcoming Proposition \ref{P6} is concerned with the asymptotic behaviour of $\mathbb{V}\mathrm{ar}[Z(t)]$. We first introduce the following lemma in preparation for its proof.
\begin{lemma}\normalfont\label{L2}
Let $G: [0,1] \rightarrow \mathbb{R}$ be a continuous function, and define the integral, for $s>0$,
\begin{equation}\label{Js}
J(s):= \int_{0}^{s}(v+s-x)^{-\gamma}x^{-\xi}G\bigg(\frac{s-x}{v+s-x}\bigg)\mathrm{d}x,\qquad s>0
\end{equation}
where $v>0$ is fixed. It is assumed that $\xi<1$ and $\gamma \ge 1$ so that the integral must converge. We have the asymptotic result, as $s\rightarrow \infty$,
\begin{equation}\label{Jsasym}
J(s) \sim\begin{cases}
s^{-\xi} v^{1-\gamma}\int_{0}^{1}(1-z)^{\gamma-2}G(z)\mathrm{d}z,& \text{if }\gamma>1, \\
s^{-\xi}(\ln s)G(1),& \text{if }\gamma =1.
\end{cases}
\end{equation}
\noindent \textbf{Proof:} See Part 3 of the `Supplementary materials'.
\hfill$\square$
\end{lemma}

\begin{proposition}\normalfont \label{P6} With Pareto delay in a fractional Poisson claim arrival process, the asymptotic behaviour of the variance of the total IBNR claims $Z(t)$ when $0<\alpha<1$ is, as $t\rightarrow \infty$,
{
	\small
\begin{equation}\label{VARZdtasyPareto1}
\mathbb{V}\mathrm{ar}[Z(t)] \sim
\begin{cases}
\left\{\dfrac{2\mu_1^2\lambda^2 \theta^{2\eta}\Gamma(\alpha+1-2\eta)\Gamma(1-\eta)}{\Gamma(2\alpha+1-2\eta)\Gamma(\alpha+1-\eta)} -\dfrac{\mu_1^2 \lambda^2 \theta^{2 \eta} [\Gamma(1-\eta)]^2}{[\Gamma(\alpha+1-\eta)]^2}\right\} t^{2(\alpha-\eta)}, &\mathrm{if}~~0<\eta<\alpha,\vspace{.07in}\\
\mu_1^2 \lambda^2 \theta^{2\alpha} [\Gamma(1-\alpha)]^2+\mu_2 \lambda \theta^\alpha \Gamma(1-\alpha), & \mathrm{if}~~\eta=\alpha,\vspace{.07in}\\
\dfrac{\mu_2 \lambda \theta^\eta \Gamma(1-\eta)}{\Gamma(\alpha+1-\eta)} t^{\alpha-\eta}, & \mathrm{if}~~\alpha<\eta<1,\vspace{.07in}\\
\dfrac{\mu_2 \lambda \theta}{\Gamma(\alpha)} t^{\alpha-1}\ln t, & \mathrm{if}~~ \eta=1,\vspace{.07in}\\
\left\{\dfrac{2\mu_1^{2}\lambda^2\theta^{\alpha+1}}{[\Gamma(\alpha)]^2}\displaystyle{\int_{0}^{1}(1-y)^{2\eta-\alpha -2}G^*(y)\mathrm{d}y} +\dfrac{\mu_2 \lambda \theta}{(\eta-1) \Gamma(\alpha)}\right\} t^{\alpha-1}, & \mathrm{if}~~\eta>1,
\\
\end{cases}
\end{equation}}
where $G^*(y)$ is given by (\ref{G1ydef}).

\noindent \textbf{Proof:} Recall that the variance of $Z(t)$ is given by (\ref{VARdzero}). Application of (\ref{FP}), (\ref{Ptail}) and (\ref{expr_EZ_phi}) yields that its first term equals
\begin{align}\label{PVZ}
&2 \mu_{1} \int^t_0 \overline{W}(t-x) \mathbb{E}[Z(t-x)] \mathrm{d}m(x)\nonumber\\
&~~=
2\mu_1\int_{0}^{t}\Big(\frac{\theta}{\theta+t-x}\Big)^{\eta}\frac{\mu_1\lambda\theta^{\eta}}{\Gamma(\alpha)(\theta+t-x)^{\eta-\alpha}}\bigg[\int_{0}^{\frac{t-x}{\theta+t-x}}(1- y)^{-\eta}y^{\alpha-1}\mathrm{d}y\bigg]\frac{\lambda x^{\alpha-1}}{\Gamma(\alpha)}\mathrm{d}x\nonumber\\
&~~=\frac{2\mu_1^{2}\lambda^2\theta^{2\eta}}{[\Gamma(\alpha)]^2}\int_{0}^{t}\frac{x^{\alpha-1}}{(\theta+t-x)^{2\eta-\alpha}}\bigg[\int_{0}^{\frac{t-x}{\theta+t-x}}(1- y)^{-\eta}y^{\alpha-1}\mathrm{d}y\bigg]\mathrm{d}x.
\end{align}

To analyze the asymptotic behaviour of $\mathbb{V}\mathrm{ar}[Z(t)]$, we mainly focus on the term (\ref{PVZ}). 
The analysis is separated into three cases as follows.

\noindent\textbf{Case 1.} $0<\eta<1$:  With the help of (\ref{ibf}) and (\ref{ohf}), \eqref{PVZ} can be rewritten as
\begin{align}
& \frac{2\mu_1^{2}\lambda^2\theta^{2\eta}}{[\Gamma(\alpha)]^2}\int_{0}^{t}\frac{x^{\alpha-1}}{(\theta+t-x)^{2\eta-\alpha}}\mathrm{B}\left(\alpha,1-\eta,\frac{t-x}{\theta+t-x}\right)\mathrm{d}x\label{PVZ1}\\
& ~~=\frac{2\mu_1^{2}\lambda^2\theta^{2\eta}}{[\Gamma(\alpha)]^2}\int_{0}^{t}\frac{x^{\alpha-1}}{(\theta+t-x)^{2\eta-\alpha}}\frac{(\tfrac{t-x}{\theta+t-x})^\alpha}{\alpha} {}_{2}\mathrm{F}_1\left(\alpha,\eta,\alpha+1, \frac{t-x}{\theta+t-x}\right)\mathrm{d}x\nonumber\\
&~~ = \frac{2\mu_1^{2}\lambda^2\theta^{2\eta} }{[\Gamma(\alpha)]^2} \int_{0}^{t}\frac{x^{\alpha-1}(t-x)^{\alpha}}{(\theta+t-x)^{2\eta}} \sum_{n=0}^{\infty}\frac{(\eta)_n}{n!(\alpha+n)}\Big(\frac{t-x}{\theta+t-x}\Big)^{n}\mathrm{d}x\nonumber\\
&~~ = \frac{2\mu_1^{2}\lambda^2\theta^{2\eta} }{[\Gamma(\alpha)]^2} \sum_{n=0}^{\infty}\frac{(\eta)_n}{n!(\alpha+n)}\int_{0}^{1}\frac{t^{2\alpha+n}}{(\theta+t)^{2\eta+n}}y^{\alpha-1}(1-y)^{\alpha+n}\Big(1-\frac{t}{\theta+t}y\Big)^{-(2\eta+n)}
\mathrm{d}y\nonumber\\
&~~=\frac{2\mu_1^2\lambda^2}{[\Gamma(\alpha)]^2}\Big(\frac{\theta}{\theta+t}\Big)^{2\eta}t^{2\alpha}\sum_{n=0}^{\infty} \Big(\frac{t}{\theta+t}\Big)^{n}\frac{(\eta)_n}{n!(\alpha+n)}\mathrm{B}(\alpha, \alpha+n+1){}_{2}\mathrm{F}_1\left(2\eta+n,\alpha, 2\alpha+n+1,\frac{t}{\theta+t}\right),\label{PVZ12}
\end{align}
where the last two lines are due to a change of variable from $x/t$ to $y$ together with (\ref{ohf1}). Note that the second and the third terms in (\ref{VARdzero}) are equivalent to $\mu_2 \mathbb{E}[Z(t)]/\mu_1$ and $\{\mathbb{E}[Z(t)]\}^2$ respectively, and their asymptotics are given by (\ref{PEZ2}) as
\begin{equation}\label{two}
\frac{\mu_2 \mathbb{E}[Z(t)]}{\mu_1} \sim \frac{\mu_2 \lambda \theta^{\eta} \Gamma(1-\eta)} {\Gamma(\alpha+1-\eta)} t^{\alpha-\eta} \text{~~and~~} \{\mathbb{E}[Z(t)]\}^2\sim
\frac{\mu_1^2 \lambda^2 \theta^{2 \eta} [\Gamma(1-\eta)]^2}{[\Gamma(\alpha+1-\eta)]^2} t^{2(\alpha-\eta)} \text{~~as~~}t\rightarrow \infty.
\end{equation}
This case can be further subdivided into two cases.


\noindent\textbf{Case 1a.} $0<\eta<\frac{\alpha+1}{2}$: Now, we shall proceed to prove that \eqref{PVZ12} is asymptotically proportional to $t^{2(\alpha-\eta)}$ as $t\rightarrow\infty$. To see this, we let
\begin{equation}\label{Ant}
\chi^*_n(t):=\Big(\frac{t}{\theta+t}\Big)^{n}\frac{(\eta)_n}{n!(\alpha+n)}\mathrm{B}(\alpha, \alpha+n+1){}_{2}\mathrm{F}_1\left(2\eta+n,\alpha, 2\alpha+n+1,\frac{t}{\theta+t}\right)
\end{equation}
and look at
\begin{align}
\lim_{t\rightarrow \infty} \frac{\frac{2\mu_1^2\lambda^2}{[\Gamma(\alpha)]^2}\big(\frac{\theta}{\theta+t}\big)^{2\eta}t^{2\alpha}\sum_{n=0}^{\infty} \chi^*_n(t)}{t^{2(\alpha-\eta)}}
&= \frac{2\mu_1^2\lambda^2 \theta^{2\eta}}{[\Gamma(\alpha)]^2}\lim_{t\rightarrow \infty} \Big(\frac{t}{\theta+t}\Big)^{2\eta} \sum_{n=0}^{\infty} \chi^*_n(t)\nonumber\\
&= \frac{2\mu_1^2\lambda^2 \theta^{2\eta}}{[\Gamma(\alpha)]^2}\lim_{t\rightarrow \infty} \sum_{n=0}^{\infty} \chi^*_n(t)\label{one0}\\
&= \frac{2\mu_1^2\lambda^2 \theta^{2\eta}}{[\Gamma(\alpha)]^2}\mathrm{B}(\alpha,\alpha+1-2\eta)\mathrm{B}(\alpha,1-\eta).\label{one}
\end{align}
See \eqref{prop6step} in Appendix \ref{PP6} for the calculation of the limit $\lim_{t\rightarrow \infty}\sum_{n=0}^{\infty} \chi^*_n(t)$ using the uniform convergence of the series $\sum^\infty_{n=0} \chi^*_n(t)$ on $t\in [0,\infty)$. Combining (\ref{two}) and (\ref{one}), the first two results of (\ref{VARZdtasyPareto1}) and part of the third result when $\alpha<\eta<\frac{\alpha+1}{2}$ follow.


\noindent\textbf{Case 1b.} $\frac{\alpha+1}{2}\leq \eta<1$: The integral in (\ref{PVZ1}) of the same form as (\ref{Js}) and thus, a direct application of (\ref{Jsasym}) with $v=\theta$, $s=t$, $\gamma=2\eta-\alpha$, $\xi=1-\alpha$, and $G(x)=\mathrm{B}(\alpha, 1-\eta, x)$ yields, as $t\rightarrow \infty$,
\[
\int_{0}^{t}\frac{x^{\alpha-1}}{(\theta+t-x)^{2\eta-\alpha}}\mathrm{B}\left(\alpha,1-\eta,\frac{t-x}{\theta+t-x}\right)\mathrm{d}x
\sim
\begin{cases}
C t^{\alpha-1}, &\mathrm{if}~~\frac{\alpha+1}{2}<\eta<1
\\
t^{\alpha-1}(\ln t)\mathrm{B}(\alpha,1-\eta), & \mathrm{if}~~ \eta=\frac{\alpha+1}{2},
\end{cases}
\]
where $C:=\theta^{\alpha-2\eta+1} \int^1_0 (1-z)^{2\eta-\alpha-2} \mathrm{B}(\alpha,1-\eta,z)\mathrm{d}z$. Because the assumption $\frac{\alpha+1}{2}\leq \eta$ implies $\alpha<\eta$ as $0<\alpha<1$, combining with the results in (\ref{two}) we observe that the second term in (\ref{VARdzero}) dominates the asymptotic behaviour which is proportional to $t^{\alpha-\eta}$. Therefore, the part of the third result of (\ref{VARZdtasyPareto1}) when $\frac{\alpha+1}{2}\leq \eta<1$ follows.

\noindent\textbf{Case 2.} $\eta\ge1$: In this case, we can re-express the integral in (\ref{PVZ}) in terms of (\ref{G1ydef}) as
  \begin{equation}\label{Atdef}
	 A^*(t):=\int_{0}^{t}\frac{x^{\alpha-1}}{(\theta+t-x)^{2\eta-\alpha}}\bigg[\int_{0}^{\frac{t-x}{\theta+t-x}}(1- y)^{-\eta}y^{\alpha-1}\mathrm{d}y\bigg]\mathrm{d}x= \int_{0}^{t}\frac{x^{\alpha-1}}{(\theta+t-x)^{2\eta-\alpha}}G^*\bigg(\frac{t-x}{\theta+t-x} \bigg)\mathrm{d}x.
	\end{equation}
It is shown in Appendix \ref{PP6} that
\begin{equation}\label{Atasym}
A^*(t) \sim t^{\alpha-1} \theta^{\alpha-2\eta +1} \int_{0}^{1}(1-y)^{2\eta-\alpha -2}G^*(y)\mathrm{d}y \text{~~as~~}t\rightarrow \infty.
\end{equation}

\noindent\textbf{Case 2a.} $\eta=1$: According to (\ref{PEZ2}), the asymptotics of the second and the third terms in (\ref{VARdzero}) are
\[
\frac{\mu_2 \mathbb{E}[Z(t)]}{\mu_1} \sim \frac{\mu_2 \lambda \theta}{\Gamma(\alpha)}t^{\alpha-1}\ln t \text{~~and~~} \{\mathbb{E}[Z(t)]\}^2\sim
\frac{\mu_1^2 \lambda^2 \theta^2}{[\Gamma(\alpha)]^2}t^{2\alpha-2}(\ln t)^2  \text{~~as~~}t\rightarrow \infty.
\]
Comparing with \eqref{Atasym}, it is clear that the second term is the dominant one, and one obtains the fourth result of (\ref{VARZdtasyPareto1}).

\noindent\textbf{Case 2b.} $\eta>1$: Again, (\ref{PEZ2}) implies that the second and the third terms in (\ref{VARdzero}) asymptotically behave like
\[
\frac{\mu_2 \mathbb{E}[Z(t)]}{\mu_1} \sim \frac{\mu_2 \lambda \theta}{(\eta-1) \Gamma(\alpha)} t^{\alpha-1}\text{~~and~~} \{\mathbb{E}[Z(t)]\}^2\sim \frac{\mu_1^2 \lambda^2 \theta^2}{(\eta-1)^2 [\Gamma(\alpha)]^2} t^{2\alpha-2} \text{~~as~~}t\rightarrow \infty.
\]
Comparison with \eqref{Atasym} reveals that the dominant terms in (\ref{VARdzero}) are the first two, giving rise to the final result in (\ref{VARZdtasyPareto1}).\hfill$\square$
\end{proposition}

\noindent{\bf Covariance of total IBNR processes.} While the proof of the previous proposition requires the asymptotic results \eqref{G1yasym}, the finer asymptotics \eqref{G1yasymfine} are needed in the proof of the upcoming proposition.
\begin{proposition}\label{CovPareto}\normalfont
Let $s>0$. With Pareto delay in a fractional Poisson claim arrival process, the asymptotic behaviour of the covariance of the total IBNR claims $Z(t)$ and $Z(s)$ when $0<\alpha<1$ is, as $t\rightarrow \infty$,
\begin{equation}\label{asymp_cov_Pareto}
\mathbb{C}\mathrm{ov}[Z(s),Z(t)]\sim 
\begin{cases}
D_1 t^{-\eta}, &\mathrm{if}~~0<\eta<2-\alpha,\\
D_2 t^{\alpha-2}, &\mathrm{if}~~\eta =2-\alpha,\\
D_3 t^{\alpha-2} , &\mathrm{if}~~\eta >2-\alpha ,
\end{cases}
\end{equation}
where
\begin{align*}
D_1:=&~\frac{\mu_2\lambda \theta^\eta s^\alpha}{\Gamma(\alpha+1)} +\frac{\mu_1\lambda \theta^\eta}{\Gamma(\alpha)} \int_0^s \mathbb{E}[Z(s-x)]x^{\alpha-1}\mathrm{d} x,\\
D_2:=&~\frac{\mu_2\lambda \theta^{2-\alpha} s^\alpha}{\Gamma(\alpha+1)} +\frac{\mu_1\lambda \theta^{2-\alpha}}{\Gamma(\alpha)}\int_0^s \mathbb{E}[Z(s-x)]x^{\alpha-1}\mathrm{d} x +\frac{\mu_1^2\lambda\theta}{\Gamma(\alpha)} \int^s_0 \overline{W}(s-x)x \mathrm{d}m(x),\\
D_3:=&~\frac{\mu_1^2\lambda\theta (1-\alpha)}{(\eta -1)\Gamma(\alpha)} \int^s_0 \overline{W}(s-x)x \mathrm{d}m(x).
\end{align*}
The expectation $\mathbb{E}[Z(s-x)]$ appearing in the above constants is given by (\ref{expr_EZ_phi}).

\noindent \textbf{Proof:} See Appendix \ref{appendix_proof_cov_Pareto}.\hfill$\square$
\end{proposition}\vspace{-.2in}

\section{Numerical Examples}\label{SEC5}\vspace{-.05in}

\begin{figure}[!h]
	\centering
\begin{minipage}{0.3\textwidth}
\includegraphics[width=\linewidth]{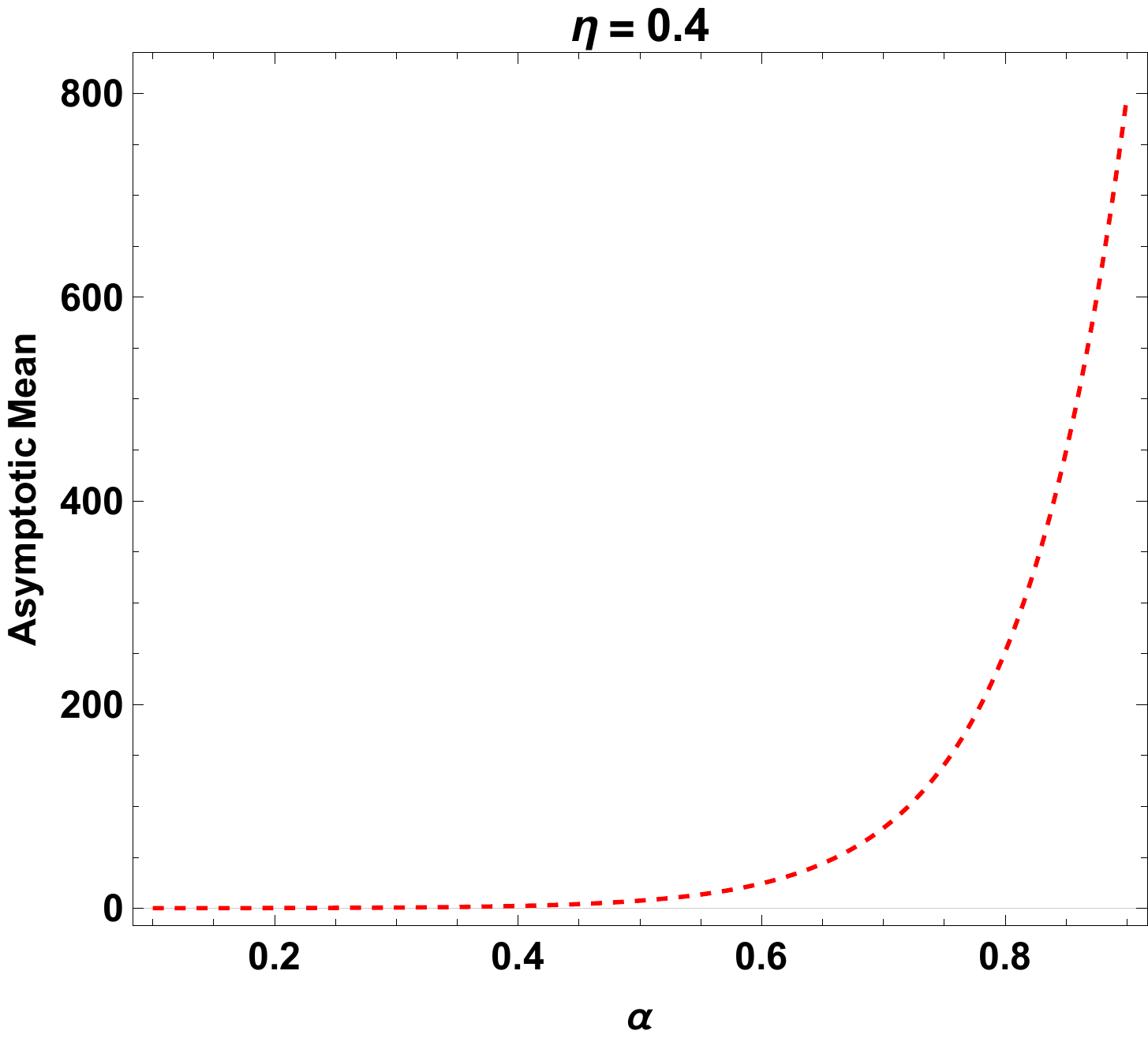}
\end{minipage}
\begin{minipage}{0.3\textwidth}
	\includegraphics[width=\linewidth]{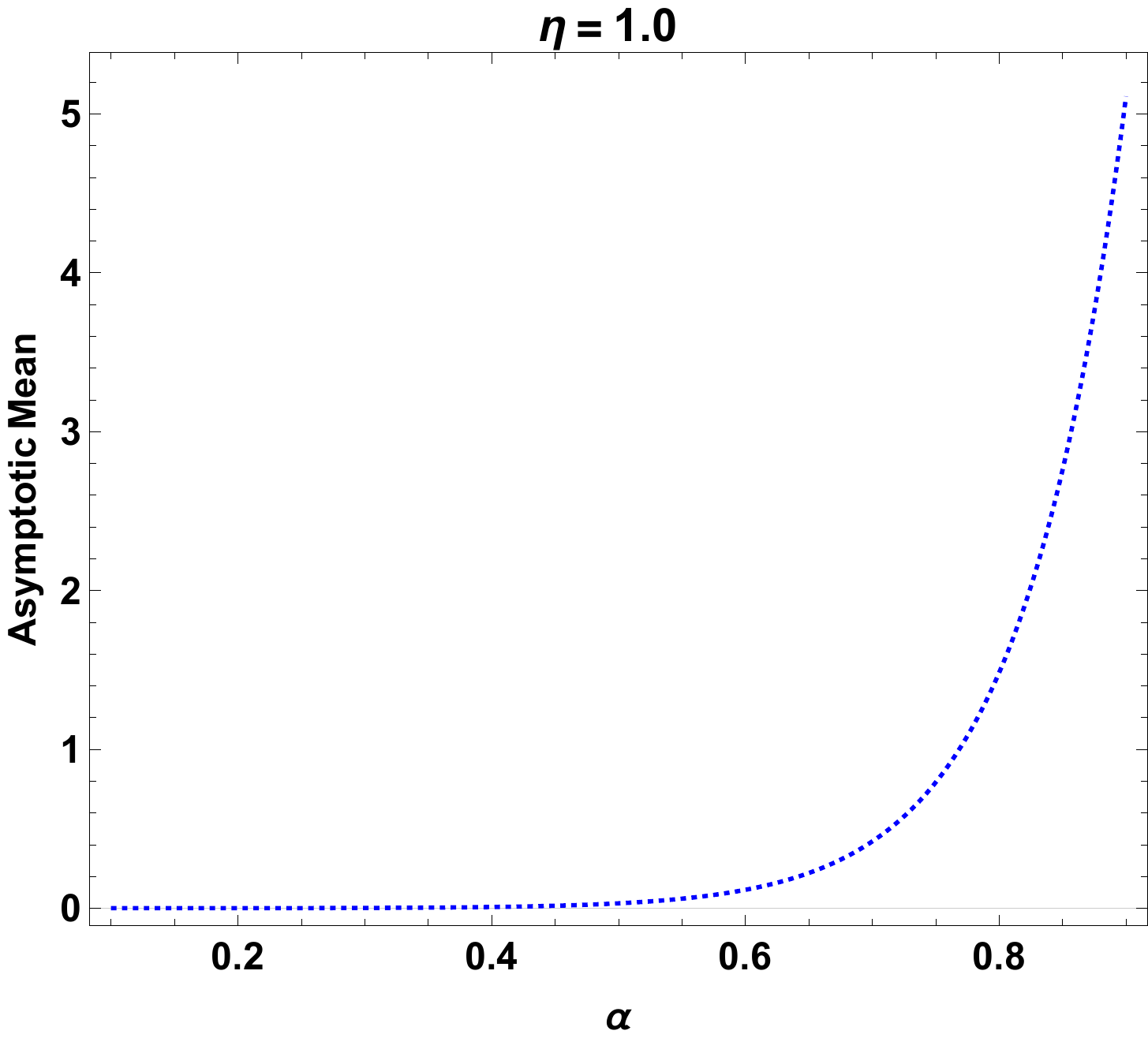}
\end{minipage}
\begin{minipage}{0.3\textwidth}
	\includegraphics[width=\linewidth]{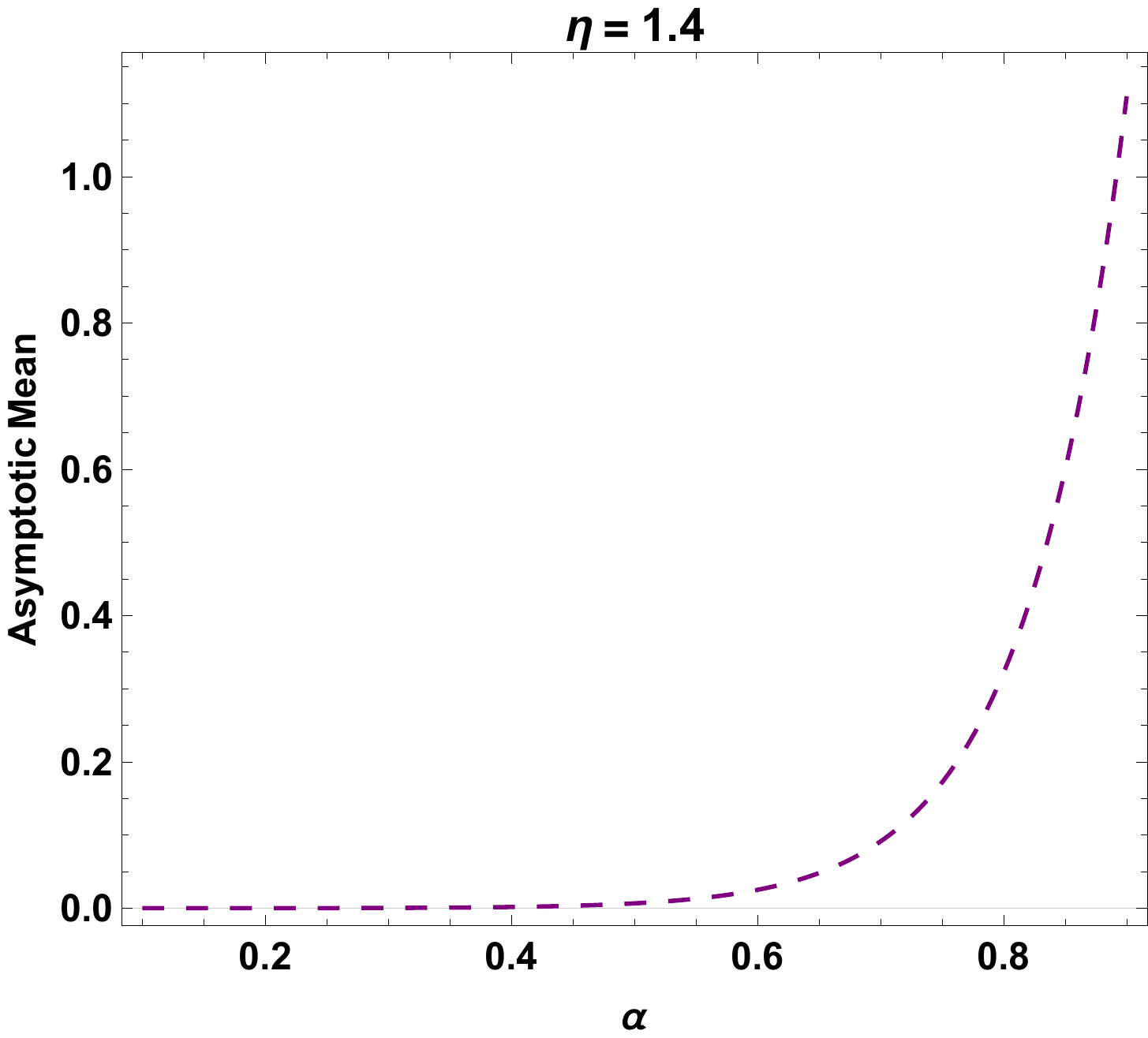}
\end{minipage}
\caption{Asymptotic mean of $Z(t)$ against $\alpha$}
\label{figue:mean}
\end{figure}

\begin{figure}[!h]
 	\centering
 	\begin{minipage}{0.3\textwidth}
 		\includegraphics[width=\linewidth]{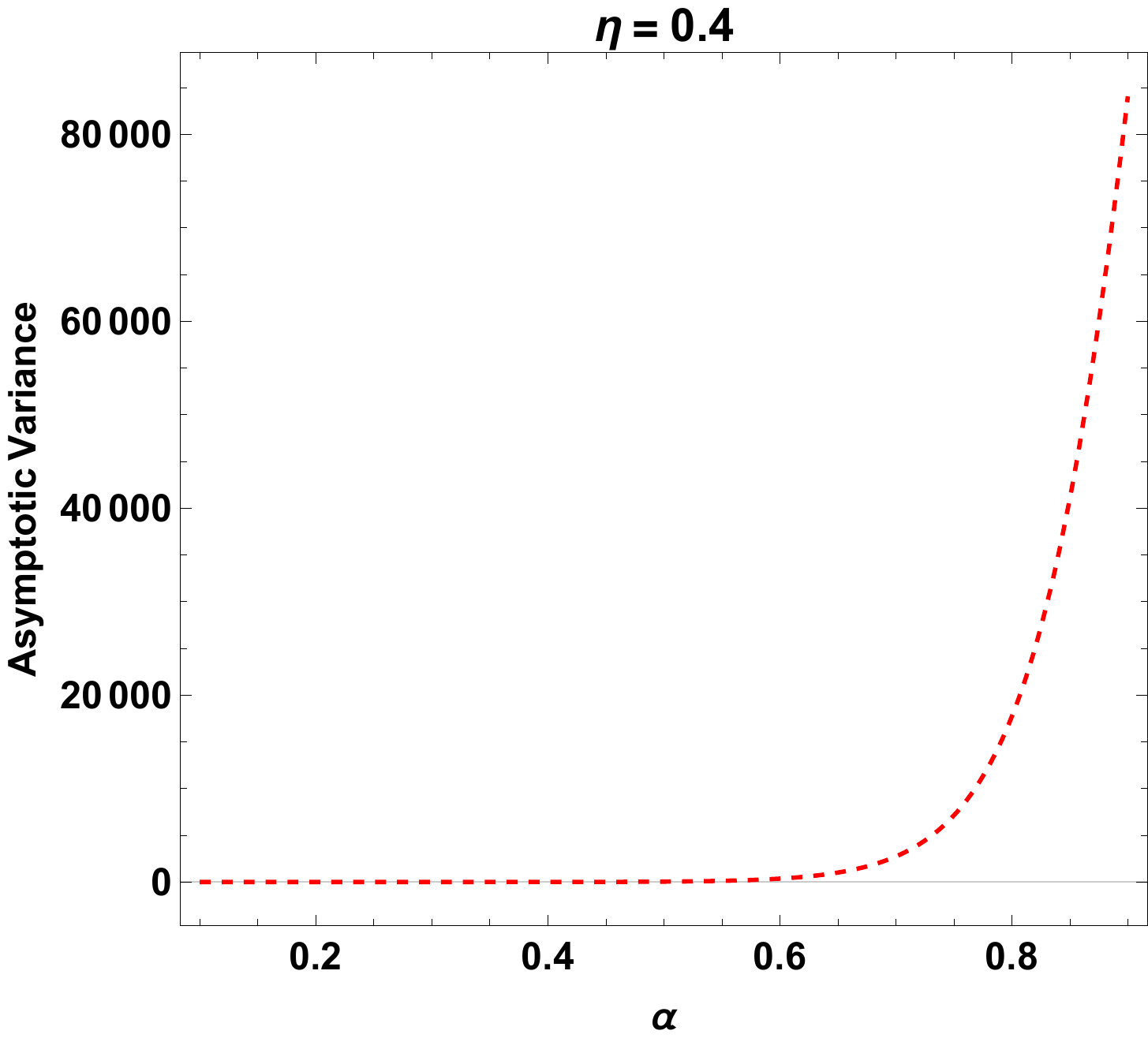}
 	\end{minipage}
 	\begin{minipage}{0.3\textwidth}
 		\includegraphics[width=\linewidth]{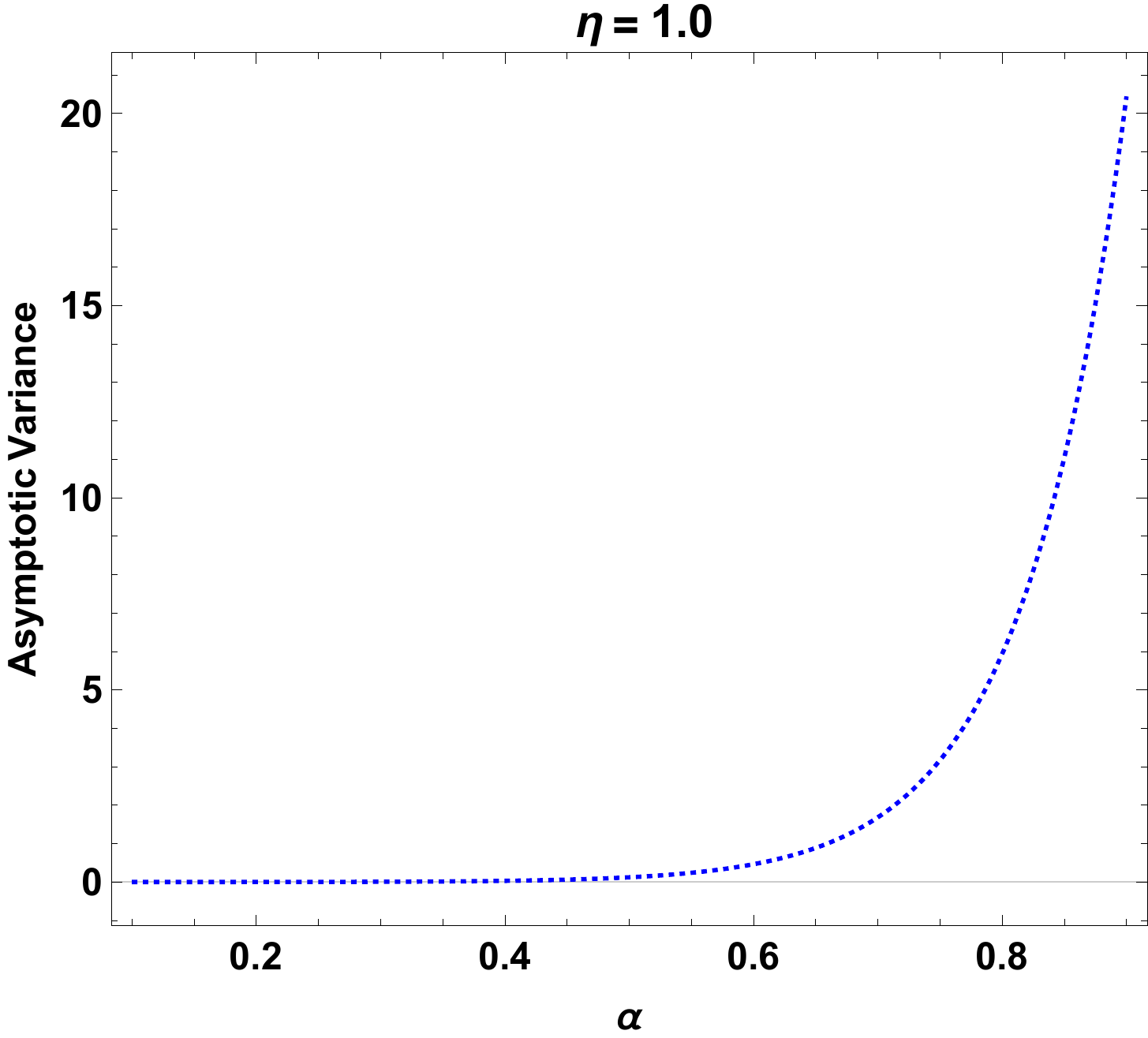}
 	\end{minipage}
 	\begin{minipage}{0.3\textwidth}
 		\includegraphics[width=\linewidth]{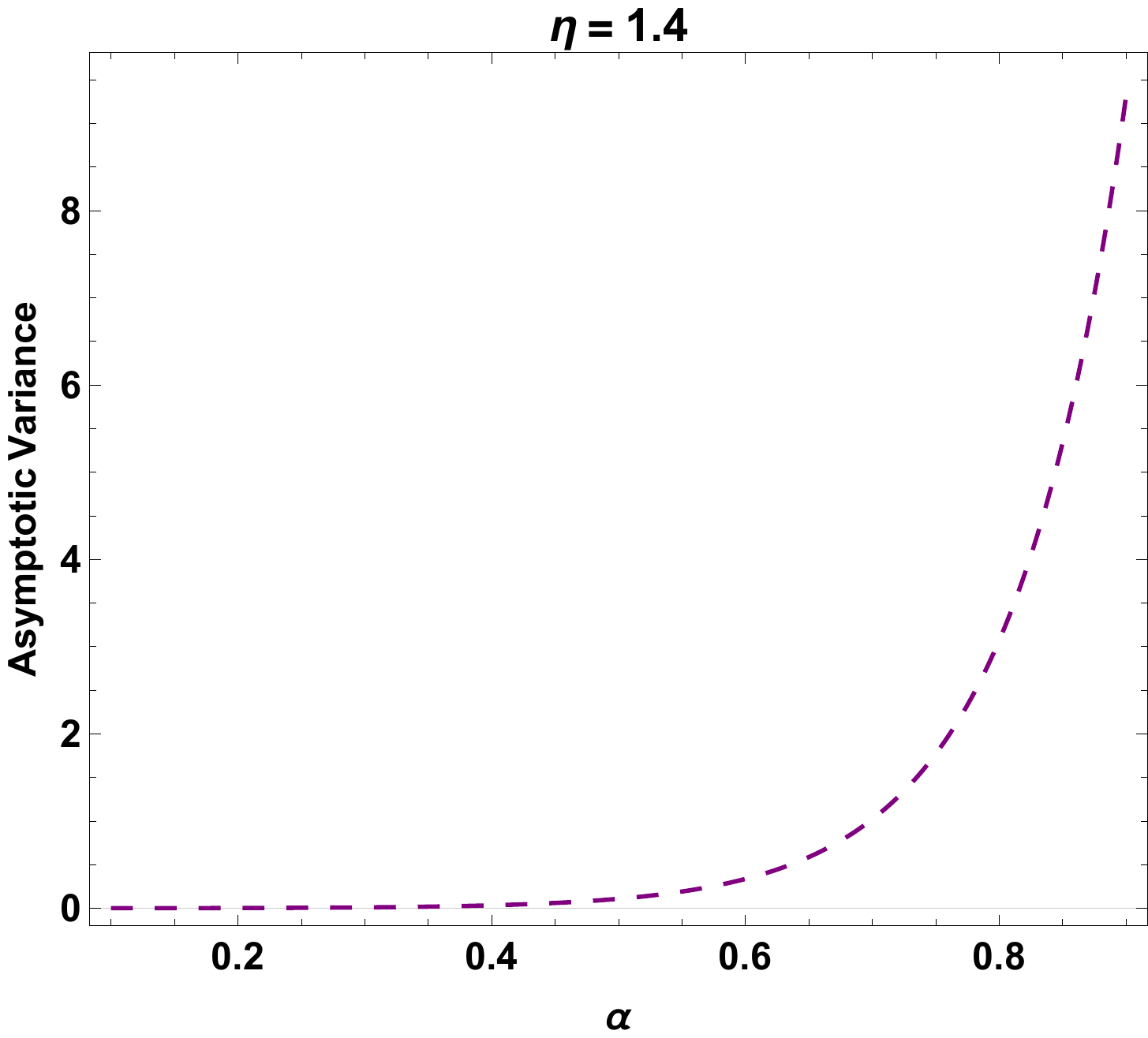}
 	\end{minipage}
 	\caption{Asymptotic variance of $Z(t)$ against $\alpha$}
 	\label{figure:variance}
 \end{figure}

\begin{figure}[!h]
	\centering
	\begin{minipage}{0.3\textwidth}
		\includegraphics[width=\linewidth]{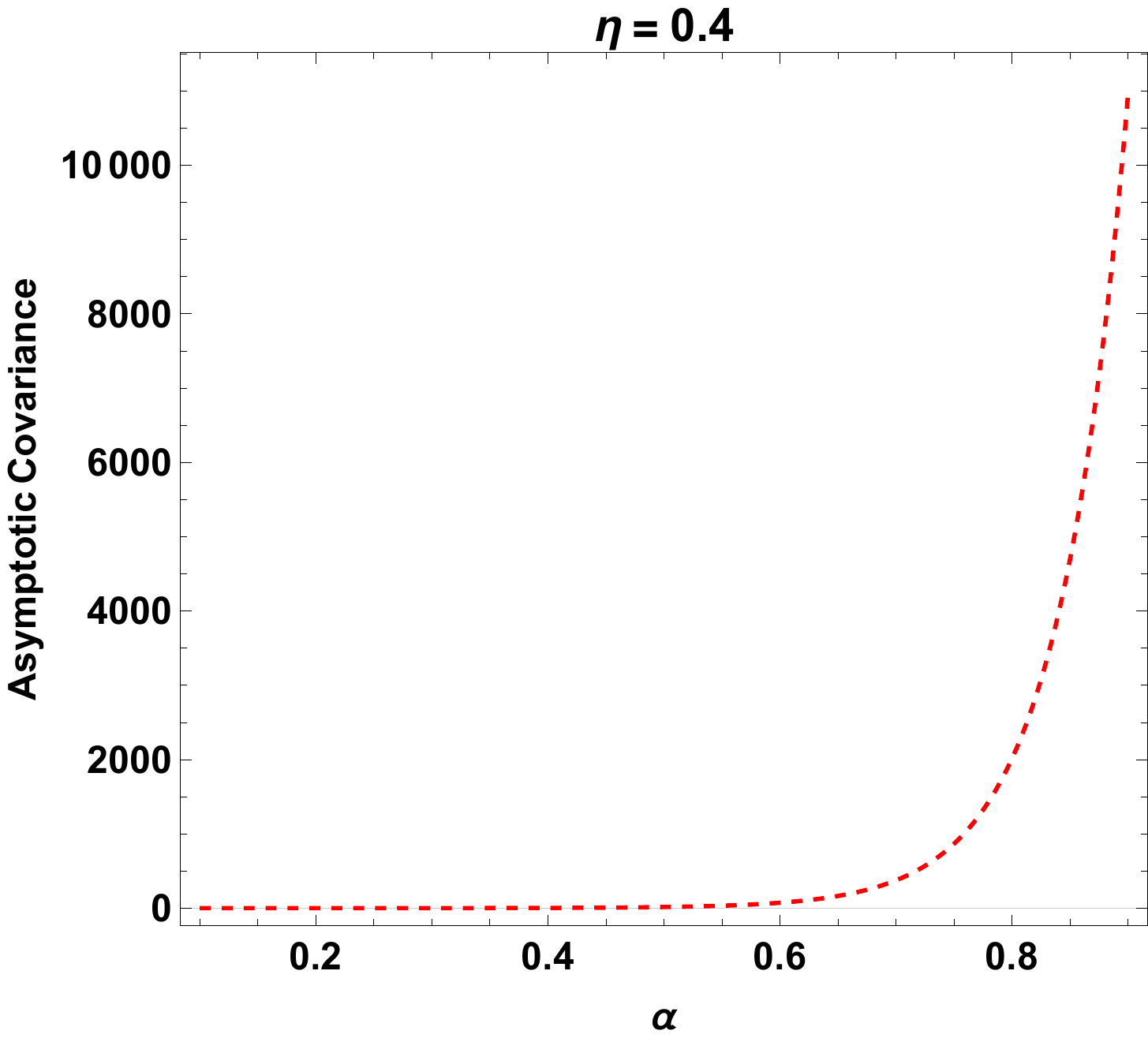}
	\end{minipage}
	\begin{minipage}{0.3\textwidth}
		\includegraphics[width=\linewidth]{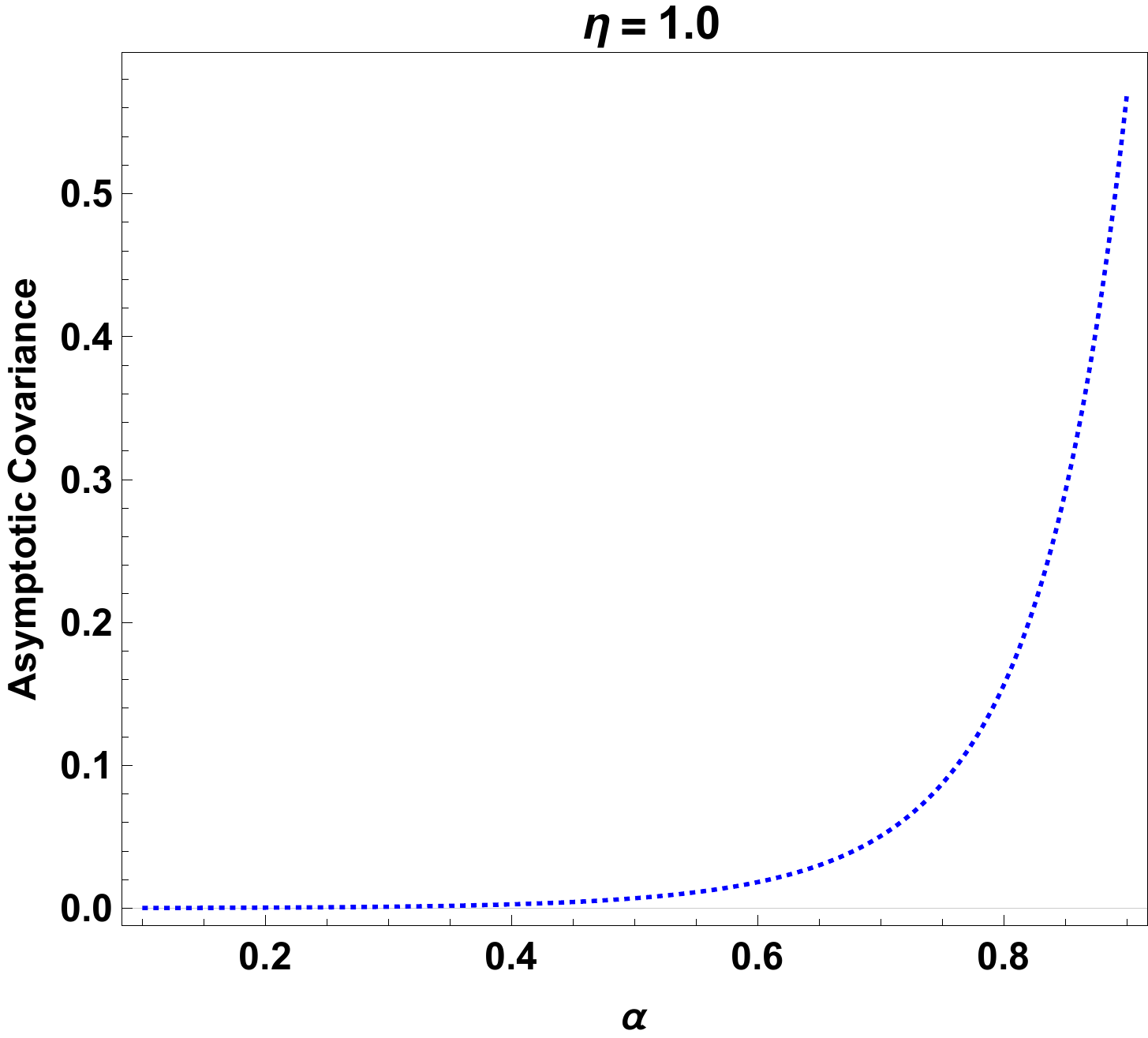}
	\end{minipage}
	\begin{minipage}{0.3\textwidth}
		\includegraphics[width=\linewidth]{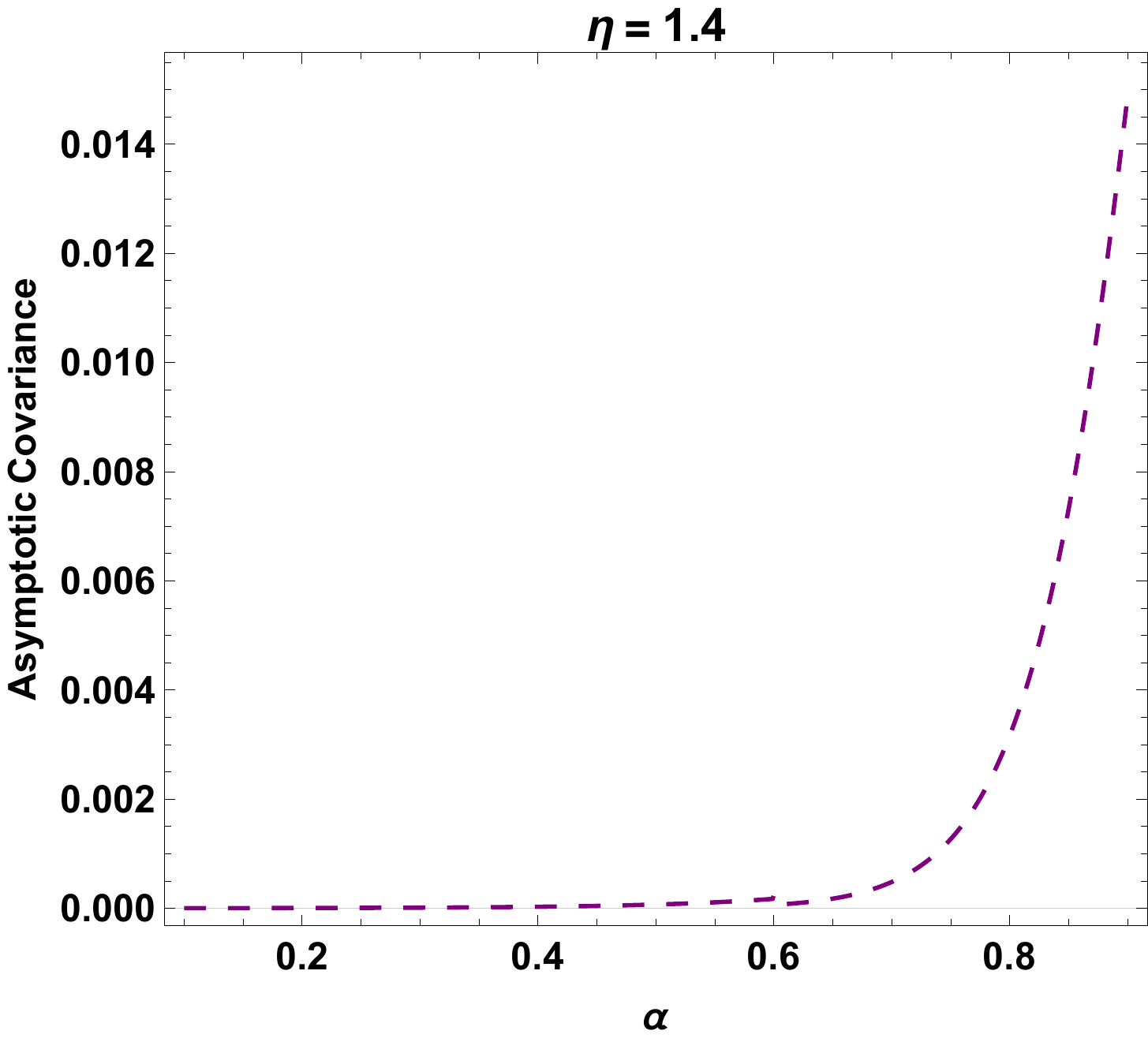}
	\end{minipage}
	\caption{Asymptotic covariance between $Z(s)$ and $Z(t)$ against $\alpha$}
	\label{figure:covariance}
\end{figure}

\begin{figure}[h]
	\centering
	\begin{minipage}{0.3\textwidth}
		\includegraphics[width=\linewidth]{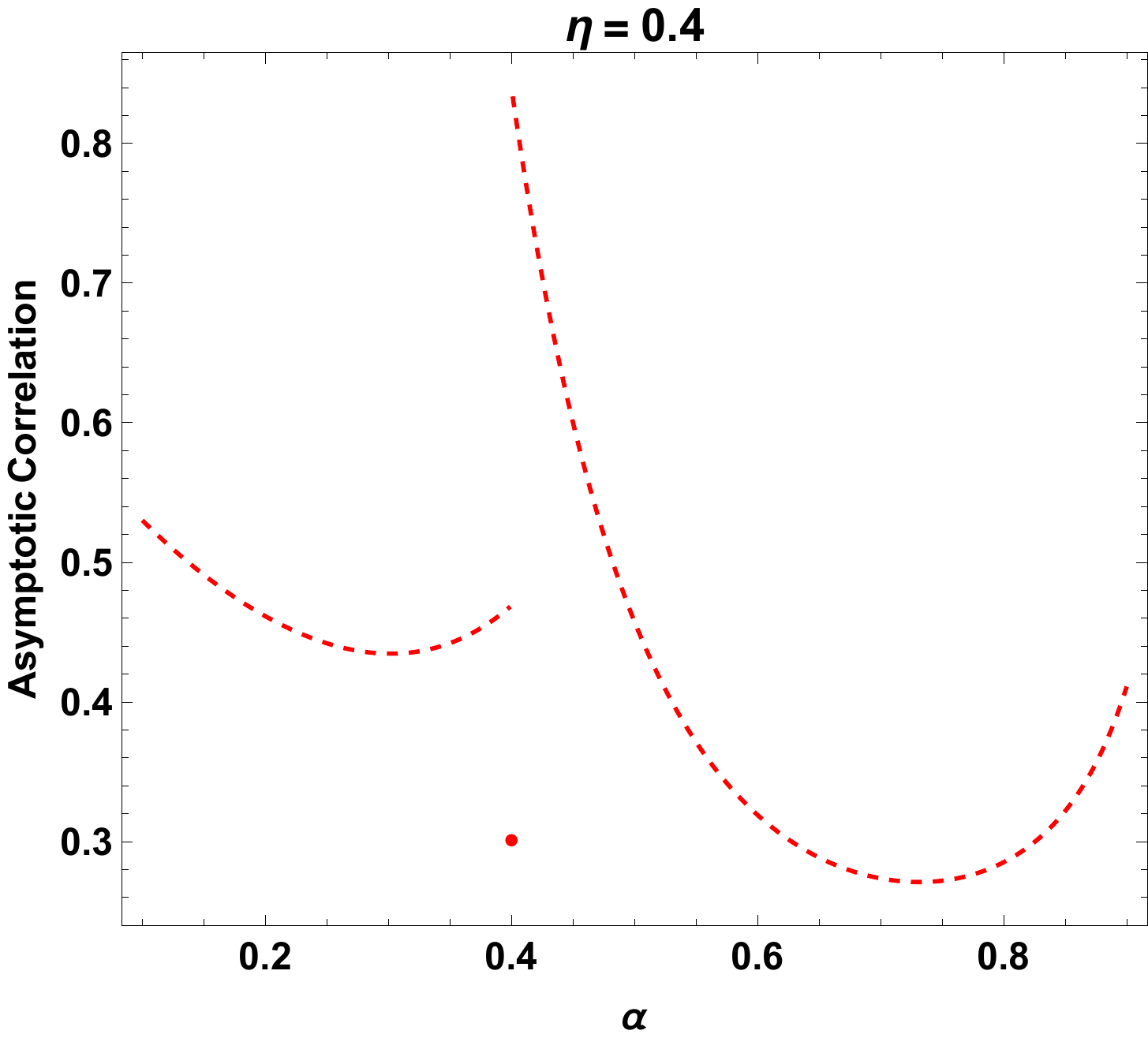}
	\end{minipage}
	\begin{minipage}{0.3\textwidth}
		\includegraphics[width=\linewidth]{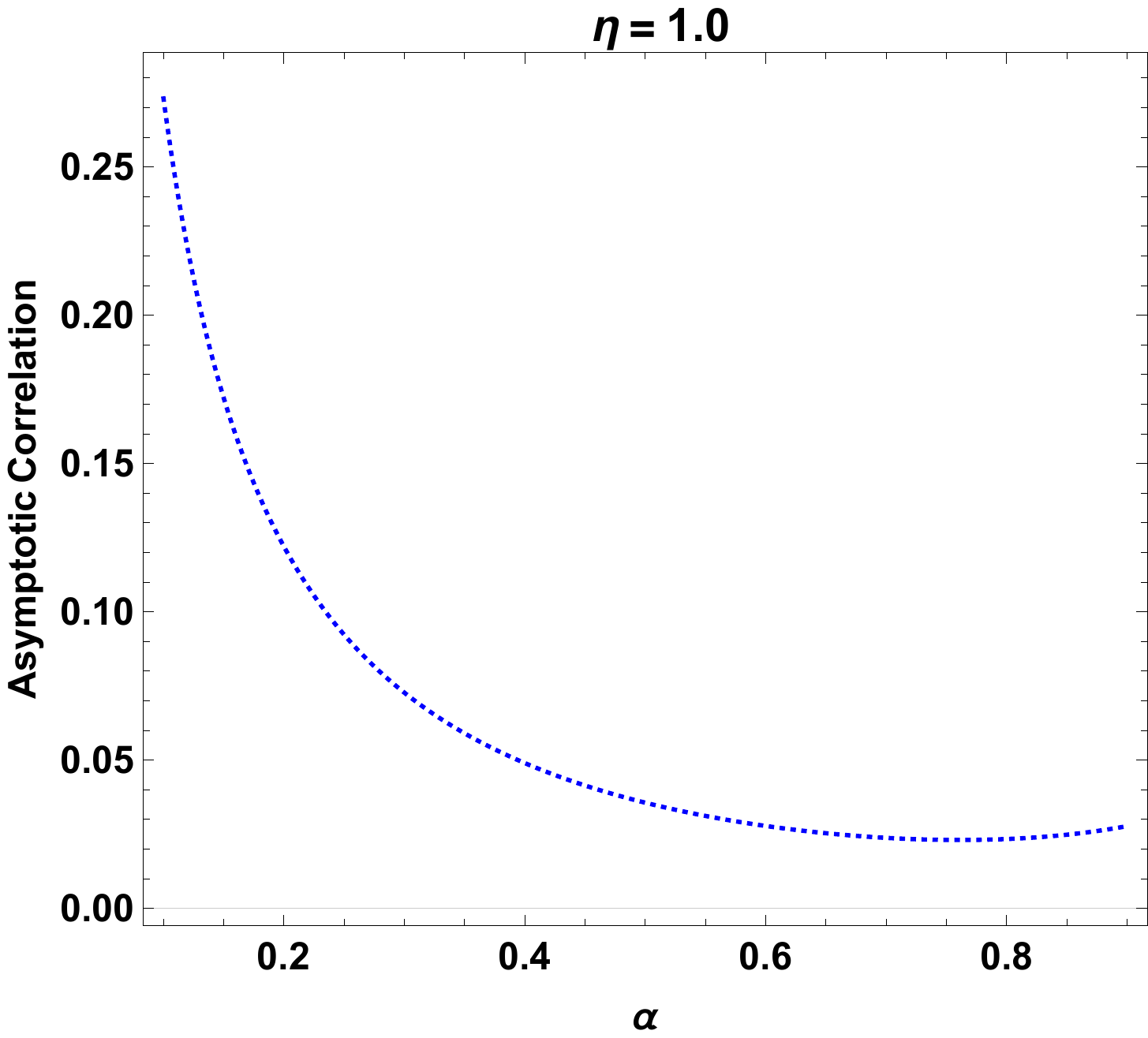}
	\end{minipage}
	\begin{minipage}{0.3\textwidth}
		\includegraphics[width=\linewidth]{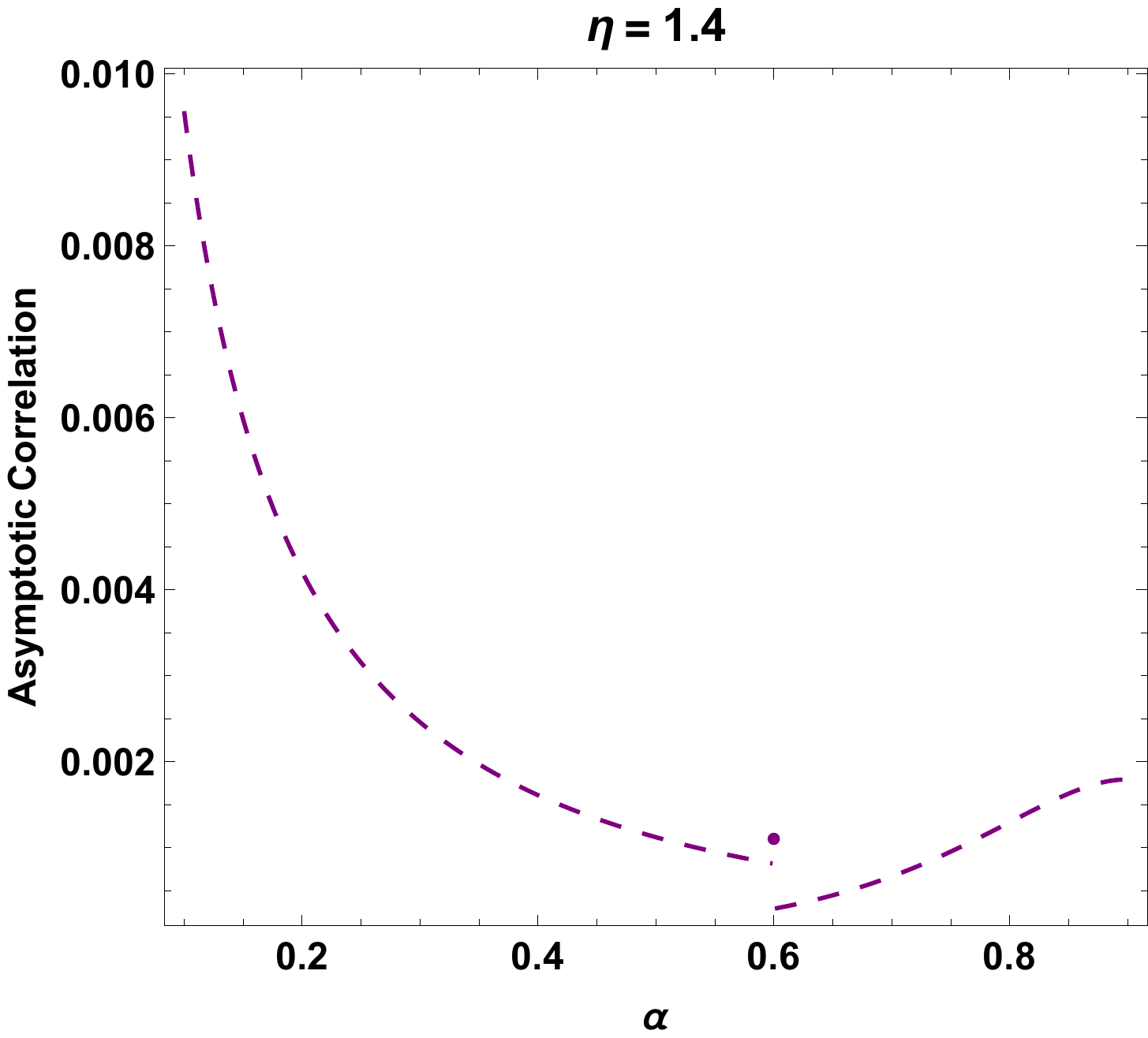}
	\end{minipage}
	\caption{Asymptotic correlation between $Z(s)$ and $Z(t)$ against $\alpha$}
	\label{figure:correlation}
\end{figure}

This section is dedicated to numerical illustrations of the theoretical results regarding some exact and asymptotic expressions of the mean, variance, covariance and correlation of IBNR processes for fractional Poisson claim arrivals with exponential and Pareto delay distributions. In all examples, it is assumed that the first two moments of the claim amount $X$ are $\mu_{1}=1$ and $\mu_{2}=4$ respectively (so that $\mathbb{V}\mathrm{ar}[X]=3$ and the coefficient of variation of $\sqrt{3}$ indicates claim with high variability). In addition, the parameter $\lambda$ of the fractional Poisson process is assumed to be $\lambda=1.5$, and except for Example \ref{example:interest} we assume a force of interest of $\delta=0$.
\begin{example}\normalfont\label{example:mean_var} \emph{(Impact of parameters of delay distribution)}
In this example, the underlying claim arrival process is assumed to have fractional Poisson parameter $\alpha=0.6$, and both ${\cal E}(\beta)$ and $\text{Pareto}(\theta,\eta)$ delays will be considered under various choices of parameters. We first calculate the asymptotic mean and variance of the IBNR claims $Z(t)$ by plugging in a large $t$ (where $t=100,000$) into our asymptotic expressions. The results are listed in Table \ref{table:asymean} and Table \ref{table:asyvar}. The corresponding exact values, whenever they can be calculated explicitly using our theoretical results, are quoted in parenthesis. The asymptotic and the exact values are very close in each case. Table \ref{table:asymean} shows that $\mathbb{E}[Z(t)]$ decreases as $\beta$ increases for exponential delays. Intuitively, this is due to the fact that ${\cal E}(\beta_1)\le_{\small \mbox{st}}{\cal E}(\beta_2)$ for $\beta_1>\beta_2$ (where $Y_1\le_{\small \mbox{st}}Y_2$ means that usual stochastic ordering $\mathbb{P}(Y_1>x)\le \mathbb{P}(Y_2>x)$ for all $x\ge0$), as stochastically smaller delays imply that incurred claims are reported quicker and there are less IBNR claims. Similarly, when one switches to Pareto delay, the monotonicity of $\mathbb{E}[Z(t)]$ in $\theta$ and $\eta$ can be explained by the stochastic orderings $\text{Pareto}(\theta,\eta_1)\le_{\small \mbox{st}}\text{Pareto}(\theta,\eta_2)$ for $\eta_1>\eta_2$ and $\text{Pareto}(\theta_1,\eta)\le_{\small \mbox{st}}\text{Pareto}(\theta_2,\eta)$ for $\theta_1<\theta_2$. Note that the values of $\mathbb{V}\mathrm{ar}[Z(t)]$ in Table \ref{table:asyvar} show the same pattern as those in Table \ref{table:asymean}. Moreover, the values of $\mathbb{E}[Z(t)]$ and $\mathbb{V}\mathrm{ar}[Z(t)]$ in Table \ref{table:asymean} are of small magnitude for exponential delays, and this is consistent with the implication of the summary in \eqref{tableau_recap1} that $\mathbb{E}[Z(t)]$ and $\mathbb{V}\mathrm{ar}[Z(t)]$ tend to zero as $t$ tends to infinity. However, for Pareto delays, the magnitude of $\mathbb{E}[Z(t)]$ and $\mathbb{V}\mathrm{ar}[Z(t)]$ is large when $\eta=0.2,0.4$ but is small when $\eta=1.0,1.2,1.4$. Such an observation is consistent with \eqref{tableau_recap2} which indicates that both $\mathbb{E}[Z(t)]$ and $\mathbb{V}\mathrm{ar}[Z(t)]$ tend to infinity or zero depending on whether $\eta<\alpha$ or $\eta>\alpha$ (recall that $\alpha$ is fixed at $0.6$). See also Remark \ref{RemarkParetomean}.

Next, the asymptotic covariance and correlation between $Z(s)$ and $Z(t)$ are provided in Table \ref{table:covariance} and Table \ref{table:correlation} respectively. In each table, we consider two scenarios with different initial time $s$ but the same time frame $t-s$. Specifically, in Case 1 we assume $s=10,000$ and $t=100,000$ while in Case 2 we let $s=20,000$ and $t=110,000$. Although the asymptotic covariances in Table \ref{table:covariance} show similar pattern to those in Table \ref{table:asymean} and Table \ref{table:asyvar} concerning means and variances, it is important to note that the covariance should always tend to zero as $t$ tends to infinity according to \eqref{tableau_recap1} and \eqref{tableau_recap2}. The high values of the covariance for Pareto delays when $\eta=0.2,0.4$ are attributed to the fact that the proportionality constant $D_1$ is large in these cases (see \eqref{asymp_cov_Pareto}). In calculating the asymptotic correlation in Table \ref{table:correlation}, we are dividing the results in Table \ref{table:covariance} by the square root of the product of the asymptotic variances $\mathbb{V}\mathrm{ar}[Z(s)]$ and $\mathbb{V}\mathrm{ar}[Z(t)]$ because both $s$ and $t$ are large. It is observed that correlation decreases when $\beta$ increases for ${\cal E}(\beta)$ delay or when $\eta$ increases for $\text{Pareto}(\theta,\eta)$ delay. However, the behaviour of the correlation in response to a change in $\theta$ for $\text{Pareto}(\theta,\eta)$ delay depends on the fixed value of $\eta$. In particular, the correlation decreases in $\theta$ when $\eta=0.2,0.4$ but increases in $\theta$ when $\eta=1.0,1.2,1.4$. For both exponential and Pareto delays, the correlation increases as the initial time $s$ increases (while keeping the time frame $t-s$ fixed).
\hfill$\square$
\end{example}\vspace{-.1in}
\begin{table}[h]
		\centering
	\parbox{0.45\linewidth}{
			\scriptsize \addtolength{\tabcolsep}{-3pt}
			\begin{tabular}{|c|c|c|c|c|c|}
			\hline
			\multicolumn{2}{|c|}{\bf{Exponential}} & \multicolumn{4}{|c|}{\bf{Pareto}}\\
			\hline
			\multirow{2}{*}{$\beta$}&\multirow{2}{*}{$\mathbb{E}[Z(t)]$} & \multirow{2}{*}{$\eta$}&\multicolumn{3}{|c|}{$\mathbb{E}[Z(t)]$}\\
			\cline{4-6}
			&& & $\theta=0.5$&$\theta=1.0$&$\theta=2.0$\\
			\hline
			\multirow{2}{*}{$0.1$}&0.100726& \multirow{2}{*}{0.2}& 171.345 & 196.824 & 226.091\\
			&(0.100730)& & (171.339)& (196.812)& (226.068)\\
			\hline
			\multirow{2}{*}{0.2}&0.050363&\multirow{2}{*}{0.4}&18.4377 & 24.3287 & 32.1020 \\
			& (0.050364)& &(18.4294)&(24.3120)&(32.0685)\\
			\hline
			\multirow{2}{*}{0.5}&0.020145& \multirow{2}{*}{1.0}&0.057982 &0.115965&0.231930\\
			& (0.020145)&&(0.066325)&(0.125668)&(0.237372)\\
			\hline
			\multirow{2}{*}{1.0}&0.010073& \multirow{2}{*}{1.2}&0.025181 &0.050363 & 0.100726\\
			& (0.010073)&&(0.023468)&(0.046426)&(0.091681)\\
			\hline
			\multirow{2}{*}{2.0}&0.005036& \multirow{2}{*}{1.4}& 0.012591 &0.025181 &0.050363 \\
			&(0.005036)&&(0.012545)&(0.025060)&(0.050041)\\
			\hline
		\end{tabular}
	\caption{Asymptotic mean of $Z(t)$ for $t=100,000$}
	\label{table:asymean}
}
\hfil
\parbox{0.45\linewidth}{
	\scriptsize\addtolength{\tabcolsep}{-3pt}
		\begin{tabular}{|c|c|c|c|c|c|}
		\hline
		\multicolumn{2}{|c|}{\bf{Exponential}} & \multicolumn{4}{|c|}{\bf{Pareto}}\\
		\hline
		\multirow{2}{*}{$\beta$}& \multirow{2}{*}{$\mathbb{V}\mathrm{ar}[Z(t)]$} & \multirow{2}{*}{$\eta$}&\multicolumn{3}{|c|}{$\mathbb{V}\mathrm{ar}[Z(t)]$}\\
		\cline{4-6}
		&& & $\theta=0.5$&$\theta=1.0$&$\theta=2.0$\\
		\hline
		\multirow{2}{*}{0.1}& 1.004398&\multirow{2}{*}{0.2}& 14,755.3 & 19,469.8 & 25,690.5\\
		& (0.994294)& &&&\\
		\hline
		\multirow{2}{*}{0.2}&0.399871& \multirow{2}{*}{0.4}&210.101 & 365.806 & 636.906 \\
		&(0.397343)&&&&\\
		\hline
		\multirow{2}{*}{0.5}&0.126382&\multirow{2}{*}{1.0}&0.231930 &0.463859 & 0.927718\\
		&(0.125977)&&&&\\
		\hline
		\multirow{2}{*}{1.0}&0.055399&\multirow{2}{*}{1.2}&0.120935&0.262715 &0.588618\\
		&(0.055298)&&&&\\
		\hline
		\multirow{2}{*}{2.0}&0.025129&\multirow{2}{*}{1.4}& 0.061263& 0.133768 &0.301616 \\
		&(0.025104)&&&&\\
		\hline
	\end{tabular}
	\caption{Asymptotic variance of $Z(t)$ for $t=100,000$}
	\label{table:asyvar}
}
	\end{table}
\vspace{-.2in}
\begin{table}[h]
		\centering
		\scriptsize
		\begin{tabular}{|c|c|c|c|c|c|c|c|c|c|}
			\hline
		\multicolumn{3}{|c|}{\bf{Exponential}}&\multicolumn{7}{|c|}{\textbf{Pareto}}\\
		\hline
		\multirow{2}{*}{$\beta$}& Case 1: $s=10,000$& Case 1: $s=20,000$& \multirow{2}{*}{$\eta$}&\multicolumn{3}{|c|}{Case 1: $s=10,000~~ t=100,000$}&\multicolumn{3}{|c|}{Case 2: $s=20,000~~ t=110,000$}\\
		\cline{5-10}
		&$t=100,000$&$t=110,000$&&$\theta=0.5$ &$\theta=1.0$&$\theta=2.0$&$\theta=0.5$&$\theta=1.0$&$\theta=2.0$\\
		\hline
		0.1&$1.02\mathrm{e}{-3}$&$1.35\mathrm{e}{-3}$ &0.2&2,131.32&2,786.52&3,646.48&4,112.99&5,389.03&7,066.23\\
		0.2&$2.55\mathrm{e}{-4}$&$3.38\mathrm{e}{-4}$&0.4&45.4345&73.6313&120.865&73.4282&119.888&198.112\\
		0.5&$4.08\mathrm{e}{-5}$&$5.41\mathrm{e}{-5}$&1.0&$8.81\mathrm{e}{-3}$&$1.82\mathrm{e}{-2}$&$3.88\mathrm{e}{-2}$&$1.20\mathrm{e}{-2}$&$2.48\mathrm{e}{-2}$&$5.22\mathrm{e}{-2}$\\
		1.0&$1.02\mathrm{e}{-5}$&$1.35\mathrm{e}{-5}$&1.2&$7.49\mathrm{e}{-4}$&$1.75 \mathrm{e}{-3}$&$ 4.16\mathrm{e}{-3}$&$1.01\mathrm{e}{-3}$&$2.35\mathrm{e}{-3}$& $ 5.54\mathrm{e}{-3}$\\
		2.0&$2.55\mathrm{e}{-6}$&$3.38\mathrm{e}{-5}$&1.4&$8.01\mathrm{e}{-5} $&$2.34\mathrm{e}{-4}$&$7.07\mathrm{e}{-4}$&$1.06\mathrm{e}{-4}$&$3.10\mathrm{e}{-4}$& $9.36\mathrm{e}{-4}$\\
		\hline
		\end{tabular}
	\caption{Asymptotic covariance between $Z(s)$ and $Z(t)$}
	\label{table:covariance}
	\end{table}
\vspace{-.2in}
\begin{table}[h]
		\centering
		\scriptsize
		\begin{tabular}{|c|c|c|c|c|c|c|c|c|c|}
			\hline
		\multicolumn{3}{|c|}{\bf{Exponential}}&\multicolumn{7}{|c|}{\textbf{Pareto}}\\
		\hline
		\multirow{2}{*}{$\beta$}& Case 1: $s=10,000$& Case 2: $s=20,000$& \multirow{2}{*}{$\eta$}&\multicolumn{3}{|c|}{Case 1: $s=10,000~~ t=100,000$}&\multicolumn{3}{|c|}{Case 2: $s=20,000~~ t=110,000$}\\
		\cline{5-10}
		&$t=100,000$&$t=110,000$&&$\theta=0.5$ &$\theta=1.0$&$\theta=2.0$&$\theta=0.5$&$\theta=1.0$&$\theta=2.0$\\
		\hline
		0.1&$6.39\mathrm{e}{-4}$&$9.94\mathrm{e}{-4}$ &0.2&$3.63\mathrm{e}{-1}$&$3.60\mathrm{e}{-1}$&$3.57\mathrm{e}{-1}$&$5.11\mathrm{e}{-1}$&$5.07\mathrm{e}{-1}$&$5.04\mathrm{e}{-1}$\\
		0.2&$4.02\mathrm{e}{-4}$&$6.24\mathrm{e}{-4}$&0.4&$3.43\mathrm{e}{-1}$&$3.19\mathrm{e}{-1}$&$3.01\mathrm{e}{-1}$&$4.73\mathrm{e}{-1}$&$4.44\mathrm{e}{-1}$&$4.21\mathrm{e}{-1}$\\
		0.5&$2.04\mathrm{e}{-4}$&$3.16\mathrm{e}{-4}$&1.0&$2.68\mathrm{e}{-2}$&$2.77\mathrm{e}{-2}$&$2.95\mathrm{e}{-2}$&$4.12\mathrm{e}{-2}$&$4.24\mathrm{e}{-2}$&$4.47\mathrm{e}{-2}$\\
		1.0&$1.16\mathrm{e}{-4}$&$1.80\mathrm{e}{-4}$&1.2&$3.91\mathrm{e}{-3}$&$4.21\mathrm{e}{-3}$&$ 4.46\mathrm{e}{-3}$&$6.15\mathrm{e}{-3}$&$6.60\mathrm{e}{-3}$& $6.95\mathrm{e}{-3}$\\
		2.0&$6.40\mathrm{e}{-5}$&$9.94\mathrm{e}{-5}$&1.4&$8.25\mathrm{e}{-4} $&$1.10\mathrm{e}{-3}$&$1.48\mathrm{e}{-3}$&$1.28\mathrm{e}{-3}$&$1.71\mathrm{e}{-3}$& $2.29\mathrm{e}{-3}$\\
		\hline
		\end{tabular}
	\caption{Asymptotic correlation between $Z(s)$ and $Z(t)$}
	\label{table:correlation}
	\end{table}
\vspace{-.2in}
\begin{example}\normalfont\label{example:alpha} \emph{(Impact of fractional Poisson parameter $\alpha$)}
In this example, the effect of $\alpha$ (where $0.1<\alpha<0.9$) on the asymptotics of the IBNR claims is examined. We focus on the case of Pareto delays with $t=100,000$ (and additionally $s=10,000$ for covariance and correlation). Moreover, we let $\theta=1.0$ and consider three scenarios where $\eta=0.4,1.0,1.4$. The results are given in Figures \ref{figue:mean}-\ref{figure:correlation}. For each fixed $\eta$, the asymptotic mean in Figure \ref{figue:mean} is increasing in $\alpha$. Intuitively, as $\alpha$ increases, the interarrival times of the claims becomes less heavy-tailed (see \eqref{asymptotics_interclaim_frac}) so that claims arrive more frequently, thereby increasing the IBNR claims. It is also noted that the asymptotic mean is continuous in $\alpha$, as evident from each of the three pieces in \eqref{PEZ2}. The asymptotic variance in Figure \ref{figure:variance} and the asymptotic covariance in Figure \ref{figure:covariance} show similar trend at a first glance. However, we remark that for the case $\eta=0.4$ (which corresponds to the first subfigure of Figure \ref{figure:variance}) there is indeed a discontinuity at $\alpha=\eta$, which can hardly be seen due to the scale of the plot. Such a curve is plotted using the first three pieces of \eqref{VARZdtasyPareto1}. In contrast, the second and third subfigures of Figure \ref{figure:variance} for the cases $\eta=1.0$ and $\eta=1.4$ are continuous in $\alpha$ and are plotted using the fourth and fifth pieces of \eqref{VARZdtasyPareto1}. Similarly, the first two subfigures of Figure \ref{figure:covariance} concerning the covariance are continuous in $\alpha$ and are both plotted using the first piece of \eqref{asymp_cov_Pareto}, but the third subfigure is discontinuous at $\alpha=0.6$ and is plotted using all three pieces of \eqref{asymp_cov_Pareto}. The aforementioned discontinuities explain those observed in the first and third subfigures of Figure \ref{figure:correlation}.\hfill$\square$
\end{example}\vspace{-.1in}
\vspace{-.2in}
\begin{example}\normalfont\label{example:changetime} \emph{(Change of correlation over time $t$)}
In this example, we fix the initial time $s=10,000$, we examine how the asymtotic correlation between $Z(s)$ and $Z(t)$ changes when $t$ increases from $100,000$ to $200,000$. The results for $\alpha=0.3,0.5,0.8$ are plotted in Figure \ref{figure:timeexppar} for both ${\cal E}(\beta=1)$ and $\text{Pareto}(\theta=1,\eta=1)$ delays. The correlation always decreases in $t$ because an IBNR claim at time $s$ is more likely to be reported by time $t$ as $t$ increases. For exponential delays (where the IBNR process is always SRD), the correlation increases as the fractional Poisson parameter $\alpha$ increases. This can be explained by \eqref{comp_corr_expo} where the correlation is always asymptotically proportional to $t^{-\frac{3-\alpha}{2}}$ which is increasing in $\alpha$. Note that the proportionality constant also depends on $\alpha$, but one expects that the effect of a change in $\alpha$ on $t^{-\frac{3-\alpha}{2}}$ dominates for large $t$. On the other hand, the correlation under Pareto delays decreases in $\alpha$ when $\alpha$ increases from 0.3 through 0.5 to 0.8 in this example. However, one should note from Example \ref{example:alpha} that in general the asymptotic correlation is not necessarily monotone in $\alpha$ (see Figure \ref{figure:correlation}).\hfill$\square$
\end{example}\vspace{-.1in}
\begin{figure}[!h]
	\centering
	\begin{minipage}{0.4\textwidth}
		\includegraphics[width=5.5cm]{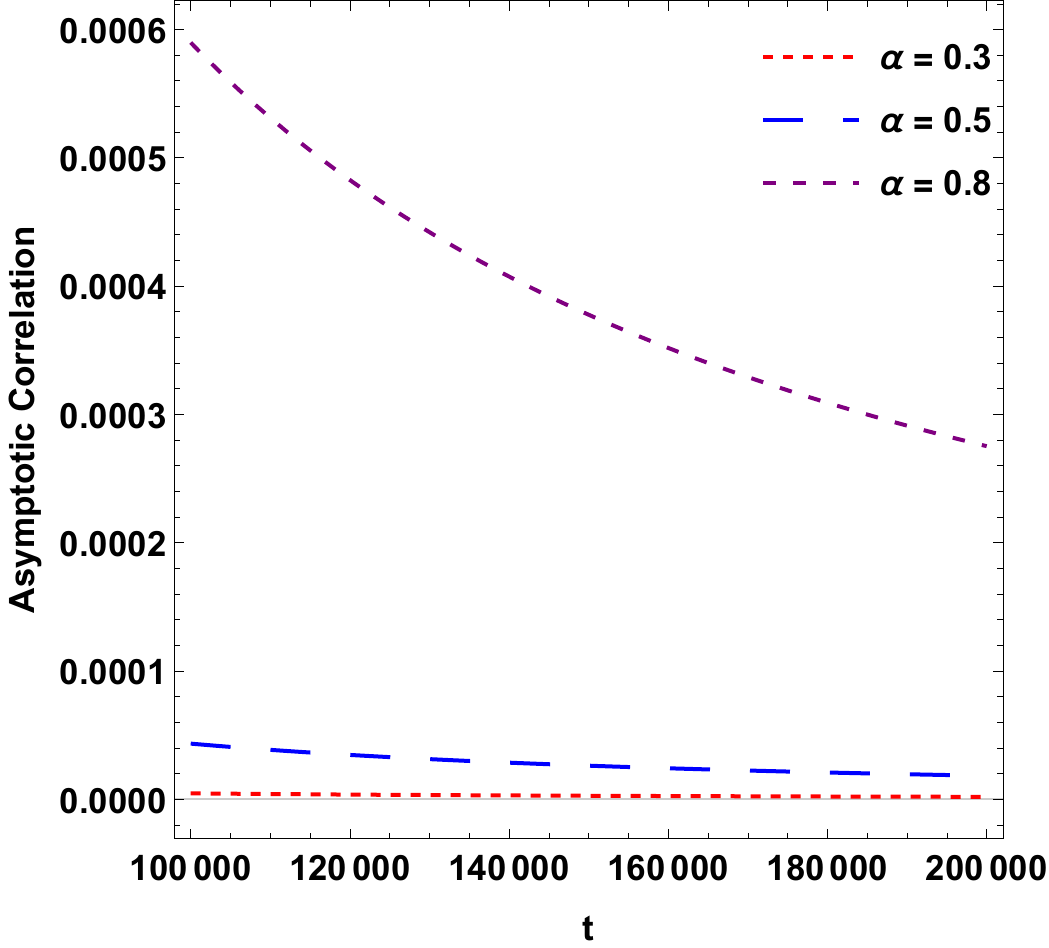}
	\end{minipage}
	\begin{minipage}{0.4\textwidth}
		\includegraphics[width=5.5cm]{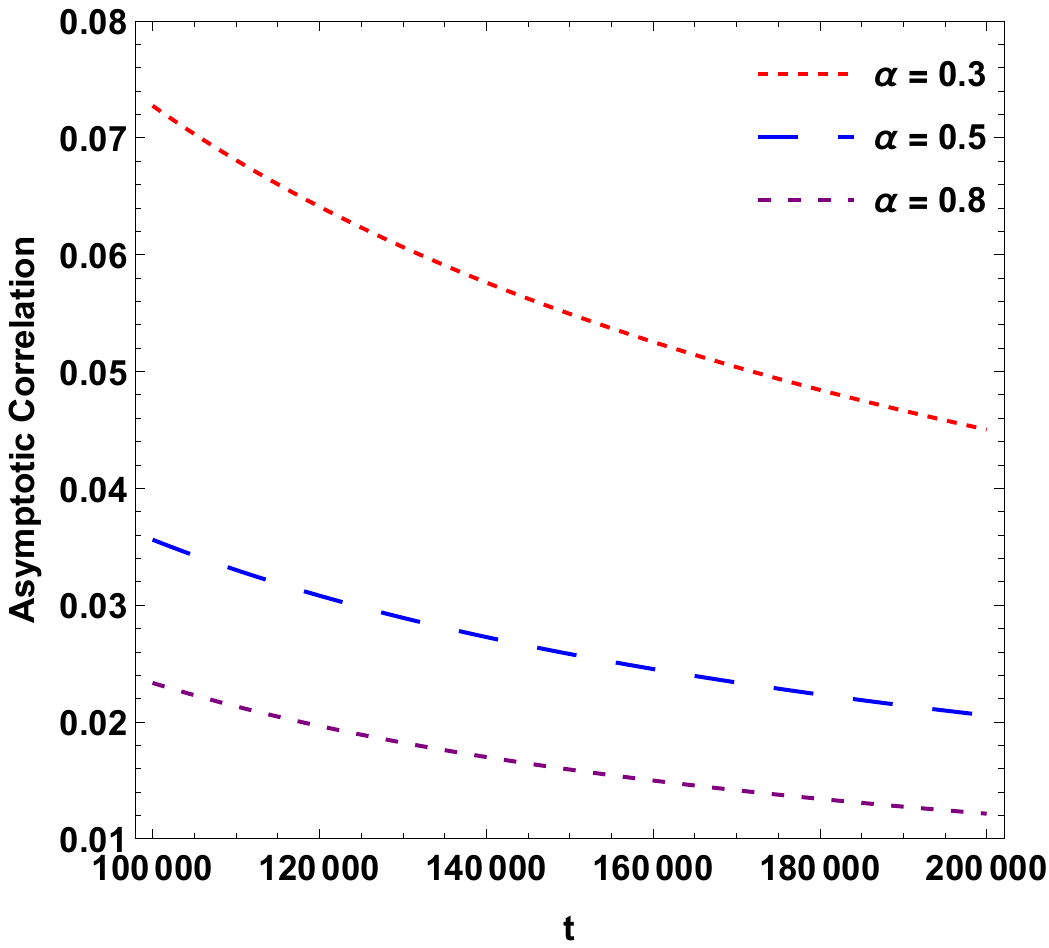}
	\end{minipage}
	\caption{Asymptotic correlation over time $t$ for exponential (left) and Pareto (right) delays}
	\label{figure:timeexppar}
\end{figure}\vspace{-.1in}
\begin{example}\normalfont\label{example:interest} \emph{(Impact of force of interest $\delta$)}
This final example briefly illustrates the effect of $\delta$ on the exact values of mean and variance of the IBNR claims for exponential delays with parameter $\beta=1.0$. We let $\alpha=0.6$ and consider four time values $t=20,50,100,200$ in Figure \ref{figure:interestmeanvar}. It can be seen that both $\mathbb{E}[Z_\delta(t)]$ and $\mathbb{V}\mathrm{ar}[Z_\delta(t)]$ quickly approaches zero as $\delta$ increases due to heavier discounting. Such an effect is more pronounced for a larger value of $t$, which is expected because of the presence of $e^{-\delta t}$ and $e^{-2\delta t}$ in the asymptotic formulas for the mean and variance respectively (see \eqref{tableau_recap1}).\hfill$\square$
\end{example}
\begin{figure}[!h]
	\centering
	\begin{minipage}{0.4\textwidth}
		\includegraphics[width=5.5cm]{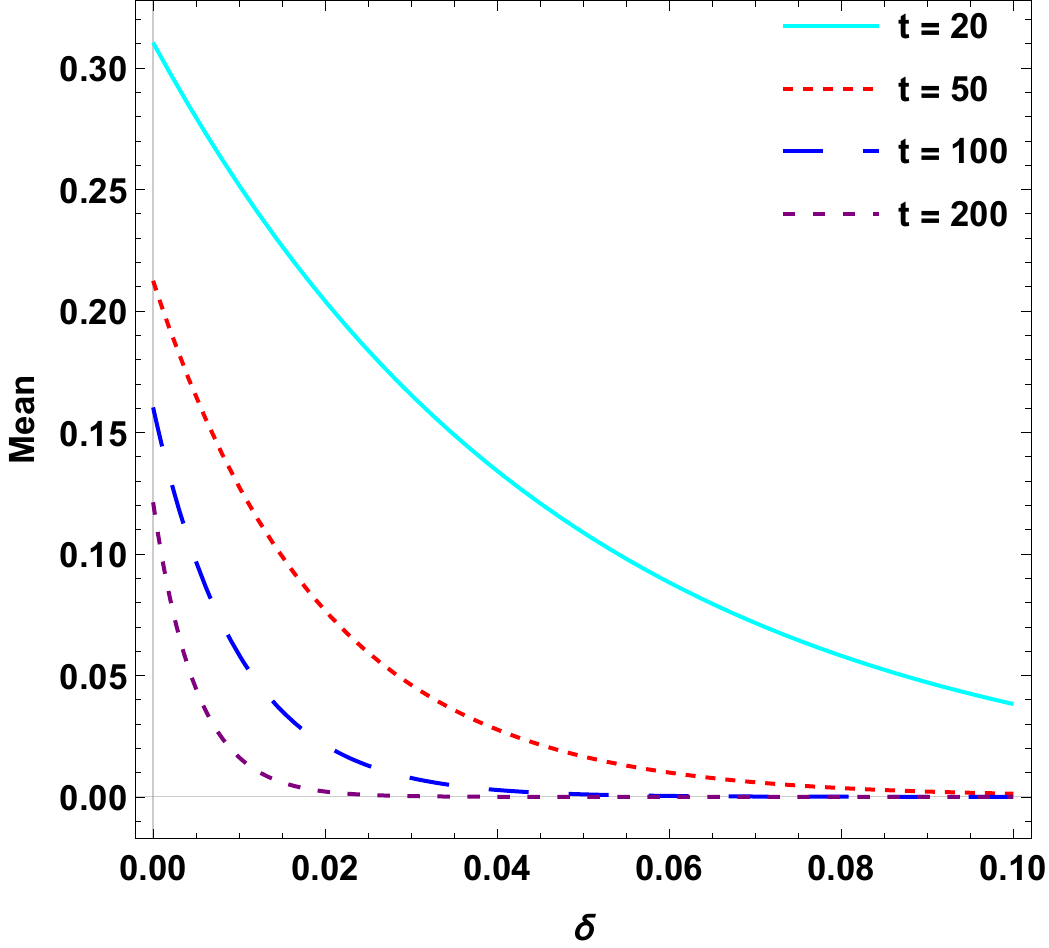}
	\end{minipage}
	\begin{minipage}{0.4\textwidth}
		\includegraphics[width=5.5cm]{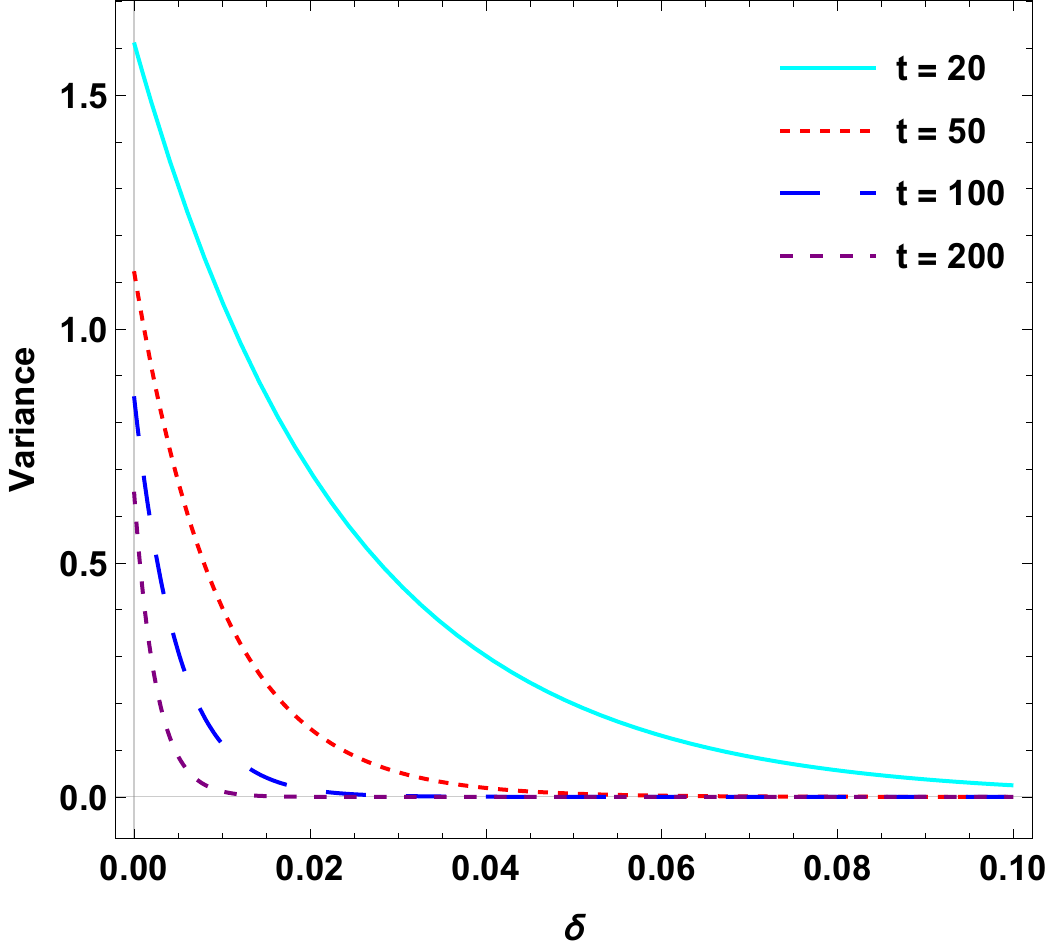}
	\end{minipage}
	\caption{Mean (left) and variance (right) against $\delta$}
	\label{figure:interestmeanvar}
\end{figure}	


\noindent\textbf{Acknowledgements}

The authors would like to thank the anonymous reviewers for helpful comments and suggestions which greatly improved an earlier version of the paper as well as Dr. Olde Daalhuis for insightful discussions concerning Part 1 of the supplementary materials in relation to the proof of Proposition \ref{P3}. Landy Rabehasaina and Jae-Kyung Woo gratefully acknowledge the support from the Joint Research Scheme France/Hong Kong Procore Hubert Curien grant No 35296 and F-HKU710/15T and the UNSW Business School 2018 International Research Collaboration Travel Funds.\vspace{-.1in}

\begin{footnotesize}

\vspace{-.2in}
 
 \end{footnotesize}

\begin{appendices}
\section{ }

\subsection{Proof of Theorem \ref{EZn}}\label{PT1}

We begin by utilizing (\ref{def1Zt}) to expand $Z_\delta^{n}(t)$ (for $n \in \mathbb{N}$ and $ t>0$) into a binomial sum as
\begin{align*}
Z_\delta^{n}(t) &= e^{-n\delta T_{1}} I_{\{ T_{1}\leq t\}} \Big[ e^{-\delta L_{1}}I_{\{ T_{1}+L_{1}>t\} } X_1 + Z_{2,\delta}(t) \Big]^{n}\\
&= e^{-n\delta T_{1}} I_{\{ T_{1}\leq t\}} \bigg[ Z_{2,\delta}^{n}(t) + \sum_{i=0}^{n-1} \binom{n}{i} e^{-(n-i)\delta L_{1}}I_{\{ T_{1}+L_{1}>t\}}X_{1}^{n-i} Z_{2,\delta}^{i}(t) \bigg].
\end{align*}
\noindent Conditioning on $T_{1}=x$, its expectation is given by
\begin{align}\label{abc}
\mathbb{E}[Z_\delta^{n}(t)] &=\int_{0}^{t} e^{-n\delta x} \bigg\{ \mathbb{E}[Z_\delta^{n} (t-x)] + \sum_{i=0}^{n-1} \mu_{n-i} \binom{n}{i} \mathbb{E} [Z_\delta^{i}(t-x)] \int_{t-x}^{\infty} e^{-(n-i)\delta y}\mathrm{d}W(y) \bigg\} \mathrm{d}F(x) \nonumber \\
&= \int_{0}^{t} e^{-n\delta x} \bigg\{ \mathbb{E}[Z_\delta^{n} (t-x)] + \sum_{i=0}^{n-1} \mu_{n-i} \binom{n}{i} \mathbb{E} [Z_\delta^{i}(t-x)] e^{-(n-i)\delta (t-x)} \mathcal{T}_{(n-i)\delta} w(t-x) \bigg\} \mathrm{d}F(x).
\end{align}
For convenience, for $t>0$ we define $\widehat{g}_n(t):=e^{n\delta t}\mathbb{E}[Z_\delta^{n}(t)]$ and $\widehat{a}_n(t):=e^{n\delta t}a_n(t)$,
where
\[
a_n(t):= \sum_{i=0}^{n-1} \mu_{n-i} \binom{n}{i} \mathbb{E} [Z_\delta^{i}(t)] e^{-(n-i)\delta t} \mathcal{T}_{(n-i)\delta} w(t).
\]
Further define the convolution operator $\ast$ as $b_1\ast b_2(t):=\int_0^t b_1(t-x)\mathrm{d}b_2(x)$ for $t>0$ where $b_1(\cdot)$ and $b_2(\cdot)$ have positive domain. Via multiplication of \eqref{abc} by $e^{n\delta t}$, one obtains
\begin{equation*}
\widehat{g}_n(t) =\widehat{g}_n\ast F(t) +\widehat{a}_n\ast F(t)
\end{equation*}
for $t>0$, which is a renewal equation satisfied by $\widehat{g}_n(\cdot)$. Therefore, its solution is given by
\begin{equation}\label{gRenewalSol}
\widehat{g}_n(t) =\widehat{a}_n\ast F(t) +\int_0^t \widehat{a}_n\ast F(t-x)\mathrm{d}m(x) =\int_0^t \widehat{a}_n(t-x)\mathrm{d}m(x),
\end{equation}
where the last line follows from the representation $m(x)=\sum_{i=1}^{\infty}F^{\ast(i)}(x)$ of the renewal function. Dividing both sides by $e^{n\delta t}$ yields the desired result \eqref{Theorem1Result}.\vspace{-.1in}


\subsection{Proof of Theorem \ref{EZnZm}}\label{PT2}

First, from (\ref{def1Zt}) with the application of binomial expansion it follows that, for $0<s\leq t$ and $n, m\in \mathbb{N}$,
\begin{align*}
&Z_\delta^{n}(s)Z_\delta^{m}(t) \\
&~= e^{-(n+m)\delta T_{1}} I_{\{ T_{1}\leq s\}} I_{\{ T_{1}\leq t\}} \Big[ e^{-\delta L_{1}}I_{\{ T_{1}+L_{1}>s\} } X_1 + Z_{2,\delta}(s) \Big]^{n} \Big[ e^{-\delta L_{1}} I_{\{ T_{1}+L_{1}>t\} } X_1 + Z_{2,\delta}(t) \Big]^{m} \nonumber \\
&~= e^{-(n+m)\delta T_{1}} I_{\{ T_{1}\leq s\}} \bigg[ Z_{2,\delta}^{n}(s) + \sum_{i=0}^{n-1} \binom{n}{i} e^{-(n-i)\delta L_{1}} I_{\{ T_{1}+L_{1}>s\} } X_{1}^{n-i} Z_{2,\delta}^{i} (s) \bigg] \\
&~~~ \times \bigg[ Z_{2,\delta}^{m}(t) + \sum_{j=0}^{m-1}\binom{m}{j} e^{- (m-j)\delta L_{1}} I_{\{ T_{1}+L_{1}>t\} } X_{1}^{m-j} Z_{2,\delta}^{j}(t) \bigg]\\
&~= e^{-(n\!+\!m)\delta T_{1}} I_{\{ T_{1}\leq s\}} \bigg[ Z_{2,\delta}^{n}(s)Z_{2,\delta}^{m}(t) + I_{\{ T_{1}
\!+\!L_{1}>s\} } \sum_{i=0}^{n-1}\binom{n}{i} e^{-(n-i)\delta L_{1}}X_{1}^{n\!-\!i} Z_{2,\delta}^{i}(s)Z_{2,\delta}^{m}(t) \\
&~~+I_{\{ T_{1}+L_{1}>t\} }\sum_{j=0}^{m-1} \sum_{i=0}^{n} \binom{n}{i}\binom{m}{j} e^{-(n+m-i-j)\delta L_{1}} X_{1}^{n+m-i-j} Z_{2,\delta}^{i}(s)Z_{2,\delta}^{j}(t)\bigg].
\end{align*}
Taking expectations on the both sides of above equation via conditioning on $T_{1}=x$ results in
\begin{align}\label{Hjoinm}
&\mathbb{E}[Z_\delta^{n}(s)Z_\delta^{m}(t)]\nonumber\\
&~= \int_{0}^{s} e^{-(n+m)\delta x} \bigg\{ \mathbb{E}[Z_\delta^{n}(s\!-\!x)Z_\delta^{m}(t\!-\!x)] + \sum_{i=0}^{n-1} \mu_{n-i} \binom{n}{i} \mathbb{E}[Z_\delta^{i}(s\!-\!x)Z_\delta^{m}(t\!-\!x)]  \int_{s-x}^{\infty} e^{-(n-i)\delta y} \mathrm{d}W(y)   \nonumber\\
&~~~+ \sum_{j=0}^{m-1} \sum_{i=0}^{n} \mu_{n\!+\!m\!-\!i\!-\!j} \binom{n}{i}\binom{m}{j} \mathbb{E}[Z_\delta^{i}(s\!-\!x) Z_\delta^{j}(t\!-\!x)]  \int_{t-x}^{\infty} e^{-(n\!+\!m\!-\!i\!-\!j)\delta y} \mathrm{d}W(y) \bigg\} \mathrm{d}F(x)\nonumber\\
&~= \int_{0}^{s} e^{-(n+m)\delta x} \bigg\{ \mathbb{E}[Z_\delta^{n}(s\!-\!x)Z_\delta^{m}(t\!-\!x)] + \sum_{i=0}^{n-1} \mu_{n-i} \binom{n}{i} \mathbb{E}[Z_\delta^{i}(s\!-\!x)Z_\delta^{m}(t\!-\!x)] e^{-(n\!-\!i)\delta (s\!-\!x)} \mathcal{T}_{(n\!-\!i)\delta} w(s\!-\!x)\nonumber\\
&~~~+ \sum_{j=0}^{m-1} \sum_{i=0}^{n} \mu_{n\!+\!m\!-\!i\!-\!j} \binom{n}{i}\binom{m}{j} \mathbb{E}[Z_\delta^{i}(s\!-\!x) Z_\delta^{j}(t\!-\!x)]e^{-(n\!+\!m\!-\! i\!-\!j)\delta (t\!-\!x)}\mathcal{T}_{(n\!+\!m\!-\!
i\!-\!j)\delta} w(t\!-\!x) \bigg\} \mathrm{d}F(x),
\end{align}
where the notation of D-H operator has been used.

Similar to the arguments used in the proof of Theorem \ref{EZn}, we proceed by defining $g_{n,m,h}(s):=\mathbb{E}[Z_\delta^{n}(s)Z_\delta^{m}(s+h)]$ and $\widehat{g}_{n,m,h}(s):=e^{(n+m)\delta s}g_{n,m,h}(s)$ for $s>0$ and fixed $h\ge 0$. Replacing $t$ by $s+h$ in \eqref{Hjoinm} and multiplying both sides by $e^{(n+m)\delta s}$, we arrive at the renewal equation
\begin{equation}\label{gnmRenewalEq}
\widehat{g}_{n,m,h}(s) =\widehat{g}_{n,m,h}\ast F(s) +\widehat{a}_{n,m,h}\ast F(s)
\end{equation}
for $s>0$, where $\widehat{a}_{n,m,h}(s):=e^{(n+m)\delta s}a_{n,m,h}(s)$ with
\begin{align*}
a_{n,m,h}(s):=&\sum_{i=0}^{n-1} \mu_{n-i} \binom{n}{i} \mathbb{E}[Z_\delta^{i}(s)Z_\delta^{m}(s\!+\!h)] e^{-(n\!-\!i)\delta s} \mathcal{T}_{(n\!-\!i)\delta} w(s)\\
&+ \sum_{j=0}^{m-1} \sum_{i=0}^{n} \mu_{n\!+\!m\!-\!i\!-\!j} \binom{n}{i}\binom{m}{j} \mathbb{E}[Z_\delta^{i}(s) Z_\delta^{j}(s\!+\!h)]e^{-(n\!+\!m\!-\! i\!-\!j)\delta (s\!+\!h)}\mathcal{T}_{(n\!+\!m\!-\!i\!-\!j)\delta} w(s\!+\!h).
\end{align*}
Analogous to \eqref{gRenewalSol}, the solution of \eqref{gnmRenewalEq} is
\[
\widehat{g}_{n,m,h}(s) =\int_0^s \widehat{a}_{n,m,h}(s-x)\mathrm{d}m(x),
\]
which is identical to \eqref{Thm2Result} upon dividing by $e^{(n+m)\delta s}$ and then replacing $h$ with $t-s$.\vspace{-.1in}

\subsection{Proof of Proposition \ref{CovExpPro}}\label{PCovExpPro}

{\bf 1. Proving the exact formula \eqref{CovExp}.} We start by proving \eqref{CovExp}. First, plugging (\ref{EZdFP}), $\mathcal{T}_\delta w(x)=\frac{\beta}{\beta+\delta}e^{-\beta x}=\mathbb{E}[e^{-\delta L}] e^{-\beta x}$ and (\ref{FP}) into the first term on the right-hand side of (\ref{Covth}) yields
\begin{align}\label{First_bis}
&\mu_{1} e^{-\delta s} \int_{0}^{s} e^{-\delta x}\mathbb{E}[Z_\delta(t-x)] \mathcal{T}_{\delta}w(s-x) \mathrm{d}m(x)\nonumber\\
&~~=\frac{\mu_1^2 \lambda^2  \alpha   \{\mathbb{E}[e^{-\delta L}]\}^2 }{[\Gamma(1+\alpha)]^2}  e^{-(\beta+\delta)(s+t)}\int^s_0 e^{2\beta x}x^{\alpha-1} (t-x)^{\alpha}
{}_1 \mathrm{F}_1(\alpha,1+\alpha,\beta (t-x))  \mathrm{d}x.
\end{align}
Again, with the application of (\ref{Ku0}) followed by changing a variable $x/s$ to $y$, we find that the integral in (\ref{First_bis}) can be expressed as
\begin{equation}\label{First1}
\sum_{n=0}^{\infty}\frac{\alpha\beta^n}{n!(\alpha+n)} s^{\alpha} \int_{0}^{1}e^{2\beta s y}y^{\alpha-1}(t- sy)^{\alpha+n}\mathrm{d}y= \sum_{n=0}^{\infty} A_n s^\alpha \mathrm{W}(2\beta s, n,t),
\end{equation}
where $A_n$ and $\mathrm{W}(2\beta s,n,t)$ are given by (\ref{An}) and (\ref{Wsnt}) respectively. Similarly, the third term of (\ref{Covth}) can be represented as
\begin{align}\label{Third_bis}
&\mu_{1} e^{-\delta t} \int_{0}^{s} e^{-\delta x} \mathbb{E}[Z_\delta(s-x)]
\mathcal{T}_{\delta} w(t-x) \mathrm{d}m(x) \nonumber\\
&~~=\frac{\mu_1^2 \lambda^2  \alpha   \{\mathbb{E}[e^{-\delta L}]\}^2 }{[\Gamma(1+\alpha)]^2}e^{-(\beta+\delta)(s+t)}
\int^s_0 e^{2\beta x}x^{\alpha-1} (s-x)^{\alpha}
{}_1 \mathrm{F}_1(\alpha,1+\alpha,\beta (s-x))  \mathrm{d}x\nonumber\\
&~~=\frac{\mu_1^2 \lambda^2  \alpha   \{\mathbb{E}[e^{-\delta L}]\}^2 }{[\Gamma(1+\alpha)]^2}e^{-(\beta+\delta)(s+t)} \sum^\infty_{n=0}A_n  s^{2\alpha+n}
{}_1 \mathrm{F}_1(\alpha,2\alpha+n+1,2\beta s),
\end{align}
where the last equality is due to (\ref{ftermVarEP1}). Lastly, using (\ref{EZdt}) and the fact that $\mathcal{T}_{\delta} w(t-x)=e^{-\beta(t-s)}\mathcal{T}_{\delta} w(s-x)$, we observe that the second term in (\ref{Covth}) is simply
\begin{equation}\label{Rest_bis1}
\mu_{2} e^{-2\delta t} \int_{0}^{s} \mathcal{T}_{2\delta} w(t-x) \mathrm{d}m(x) =\frac{\mu_2}{\mu_1}e^{-(\beta+2\delta)(t-s)}\mathbb{E}[Z_{2\delta}(s)],
\end{equation}
which corresponds to the second term on the right-hand side of \eqref{CovExp} thanks to (\ref{EZdFP}). It is also clear that the last terms in (\ref{Covth}) and \eqref{CovExp} are identical again because of (\ref{EZdFP}). Therefore, combining this with (\ref{First_bis})-(\ref{Third_bis}) yields (\ref{CovExp}).

\noindent{\bf 2. Proving the asymptotic formula \eqref{COVZdstasy}.} We now proceed to prove the asymptotic result \eqref{COVZdstasy}. From \eqref{Third_bis} and \eqref{Rest_bis1}, it is noted that the third and the second terms of \eqref{Covth} are $O(e^{-(\beta+\delta) t})$ and $O(e^{-(\beta+2\delta) t})$ respectively.
In what follows, it will be shown that these two terms are asymptotically dominated by the first and the last terms of (\ref{Covth}). Using \eqref{EZdt}, \eqref{EZEdZ} and $\mathcal{T}_\delta w(s-x)=\mathbb{E}[e^{-\delta L}] e^{-\beta(s-x)}$, the sum of the first and the last terms of (\ref{Covth}) can be expressed as
\begin{align}\label{first_last_cov_expo}
&\mu_{1} e^{-\delta s} \int_{0}^{s} e^{-\delta x}\mathbb{E}[Z_\delta(t-x)] \mathcal{T}_{\delta}w(s-x) \mathrm{d}m(x)  -\mathbb{E}[Z_\delta(s)] \mathbb{E}[Z_\delta(t)]\nonumber\\
&~~=\mu_{1} e^{-\delta s} \int_{0}^{s} \big\{e^{-\delta x}\mathbb{E}[Z_\delta(t-x)]  - \mathbb{E}[Z_\delta(t)]\big\} \mathcal{T}_{\delta}w(s-x) \mathrm{d}m(x)\nonumber\\
&~~=\mu_{1} \{\mathbb{E}[e^{-\delta L}]\}^2 e^{-(\beta+\delta) s} e^{-\delta t} \int_{0}^{s} e^{\beta x} \big\{\mathbb{E}[Z(t-x)]  - \mathbb{E}[Z(t)]\big\} \mathrm{d}m(x).
\end{align}
In view of the integrand above, one needs to first study the asymptotics of $\mathbb{E}[Z(t-x)] - \mathbb{E}[Z(t)]$ as $t\to\infty$ for all $x\in[0,s]$ (where $s\in(0,t]$ is fixed). By applying (\ref{expression_expectation}), we rewrite $\mathbb{E}[Z(t-x)] - \mathbb{E}[Z(t)]$ as, for $x\in[0,s]$,
\begin{equation}\label{first_last_cov_expo_decompose}
\mathbb{E}[Z(t-x)] - \mathbb{E}[Z(t)] = \frac{\mu_1\lambda}{\Gamma(\alpha)}[t^{\alpha-2} \Xi(t,x) - \Pi(t,x)],
\end{equation}
where
\begin{align}
\Xi(t,x):=& \int_0^1  [e^{\beta x (1-y)}-1]t^2 e^{-\beta t (1-y)}y^{\alpha-1}\mathrm{d} y,\label{XiDef}\\
\Pi(t,x):=&~[t^\alpha-(t-x)^\alpha] \int_0^1  e^{-\beta (t-x) (1-y)}y^{\alpha-1}\mathrm{d} y,\label{PiDef}
\end{align}
will be analyzed separately.

{\bf Term $\Xi(t,x)$.} Performing the change of variable $v=t(1-y)$ to the term $\Xi(t,x)$ in \eqref{XiDef} yields
\begin{equation}\label{XiDef2}
\Xi(t,x) =\int_0^t (e^{\beta x v/t}-1) t e^{-\beta v}\left( 1-\frac{v}{t}\right)^{\alpha -1} \mathrm{d} v =\int_0^{t/2} + \int_{t/2}^t:= \Xi_1(t,x)+ \Xi_2(t,x)
\end{equation}
with the obvious definitions of $\Xi_1(t,x)$ and $\Xi_2(t,x)$. We first look at $\Xi_1(t,x)$. Note that $(e^{\beta x v/t}-1) t =\sum_{j=1}^\infty (\beta x v)^j/(j! t^{j-1})$ is decreasing in $t>0$. Therefore, for $t$ large enough such that $t\ge 2s$, one has that, for $x\in[0,s]$,
\begin{align}\label{Xi1integrandbound}
I_{\{0\le v\le t/2\}}(e^{\beta x v/t}-1) t e^{-\beta v}\left( 1-\frac{v}{t}\right)^{\alpha -1} &\le  I_{\{0\le v\le t/2\}}(e^{\beta x v/2s}-1) 2s e^{-\beta v}\left( 1-\frac{v}{t}\right)^{\alpha -1}\nonumber\\
&\le I_{\{0\le v\le t/2\}}(e^{\beta v/2}-1) 2s e^{-\beta v} \frac{1}{2^{\alpha-1}}\nonumber\\
&\le \frac{s}{2^{\alpha-2}}(e^{-\beta v/2}-e^{-\beta v}),
\end{align}
which is integrable with respect to $v$ on $[0,\infty)$. For all $v\ge0$, $(e^{\beta x v/t}-1) t$ and $(1-\frac{v}{t})^{\alpha-1}$ tend to $\beta x v$ and $1$ respectively as $t\rightarrow\infty$. Then, by dominated convergence we arrive at
\begin{align}\label{Xi1asymp}
\lim_{t\rightarrow\infty}\Xi_1(t,x) =&\int_0^\infty \left[\lim_{t\rightarrow\infty}I_{\{0\le v\le t/2\}}(e^{\beta x v/t}-1) t e^{-\beta v}\left( 1-\frac{v}{t}\right)^{\alpha -1}\right]\mathrm{d}v
=\int_0^\infty \beta x v e^{-\beta v}\mathrm{d}v\nonumber\\
=&~\frac{x}{\beta}~~\text{for all}~~x\in[0,s].
\end{align}
As for $\Xi_2(t,x)$, changing variable $z=v/t$ leads to
\[
\Xi_2(t,x)=\int_{1/2}^1 (e^{\beta xz}-1)t^2 e^{-\beta t z} (1-z)^{\alpha-1}\mathrm{d}z,
\]
which can be upper bounded as
\begin{equation}\label{ineg_TCD_Xi2}
|\Xi_2(t,x)|\le  (e^{\beta s}-1)t^2 e^{-\beta t /2} \int_{1/2}^1(1-z)^{\alpha-1}\mathrm{d}z~~\text{for all}~~x\in[0,s]~~\text{and}~~t\ge s,
\end{equation}
which tends to $0$ as $t\rightarrow\infty$.
Combining with \eqref{Xi1asymp}, we conclude that the limit of \eqref{XiDef2} is given by
\begin{equation}\label{conv_Xi}
\lim_{t\to \infty}\Xi(t,x)=\frac{x}{\beta}~~\text{for all}~~x\in[0,s].
\end{equation}

Note that the inequality \eqref{Xi1integrandbound} implies the uniform upper bound $|\Xi_1(t,x)|\le C_{\Xi_1}$ where $C_{\Xi_1}:=\frac{s}{2^{\alpha-2}\beta}$ for all $x\in[0,s]$ and $t\ge 2s$. Moreover, because $t^2 e^{-\beta t /2}$ is bounded on $t\ge0$, \eqref{ineg_TCD_Xi2} also means that $|\Xi_2(t,x)|\le C_{\Xi_2}$ for all $x\in[0,s]$ and $t\ge s$ where $C_{\Xi_2}$ is some constant. As a result, we arrive at the upper bound
\begin{equation}\label{upper_bound_Xi}
|\Xi(t,x)| \le C_\Xi ~~\text{for all}~~x\in[0,s]~~\text{and}~~t\ge 2s,
\end{equation}
where $C_\Xi:=C_{\Xi_1}+C_{\Xi_2}$.

{\bf Term $\Pi(t,x)$.} Concerning $\Pi(t,x)$ in \eqref{PiDef}, we first consider, by a change of variable $v=t(1-y)$,
\begin{equation}\label{Pi_convergence}
t\int_0^1  e^{-\beta (t-x) (1-y)}y^{\alpha-1}\mathrm{d} y=\int_0^t e^{-\beta v}e^{\beta xv/t}\left( 1-\frac{v}{t}\right)^{\alpha -1}\mathrm{d} v = \int_0^{t/2} + \int_{t/2}^t := \Pi_1^*(t,x)+ \Pi_2^*(t,x)
\end{equation}
where $\Pi_1^*(t,x)$ and $\Pi_2^*(t,x)$ have obvious definitions. Note that $\Pi(t,x)$ can be rewritten as
\begin{equation}\label{PiDef2}
\Pi(t,x)= \bigg[\frac{1-\left(1-\frac{x}{t}\right)^\alpha}{\frac{1}{t}}\bigg][\Pi_1^*(t,x)+ \Pi_2^*(t,x)]t^{\alpha-2}.
\end{equation}
Applying similar arguments used for analyzing $\Xi_1(t,x)$ and $\Xi_2(t,x)$ gives rise to
\[
\lim_{t\to \infty}\Pi_1^*(t,x) =\int_0^\infty e^{-\beta v} \mathrm{d} v=\frac{1}{\beta}~~\text{for all}~~x\in[0,s]
\]
by dominated convergence, and
\begin{equation}\label{upper_bound_Pi2}
|\Pi_2^*(t,x)| = t\int_{1/2}^1 e^{-\beta(t-x)z} (1-z)^{\alpha-1}\mathrm{d}z \le e^{\beta s} t e^{-\beta t}\int_{1/2}^1 (1-z)^{\alpha-1}\mathrm{d}z~~\text{for all}~~x\in[0,s]~~\text{and}~~t\ge s,
\end{equation}
\noindent which tends to $0$ as $t\rightarrow\infty$. Consequently, (\ref{Pi_convergence}) converges to $1/\beta$ as $t\rightarrow\infty$, and then one uses \eqref{PiDef2} to evaluate
\begin{equation}\label{equiv_Pi}
\lim_{t\to \infty}\frac{\Pi(t,x)}{t^{\alpha-2}} =\bigg[\lim_{t\to \infty}\frac{1-\left(1-\frac{x}{t}\right)^\alpha}{\frac{1}{t}}\bigg]\frac{1}{\beta}=\frac{\alpha x}{\beta}~~\text{for all}~~x\in[0,s].
\end{equation}

In addition, for all $x\in[0,s]$ and $t\ge s$, one may check that $|\Pi_1^*(t,x)|\le C_{\Pi_1^*}$ where $C_{\Pi_1^*}:=\frac{e^{\beta s/2}}{2^{\alpha-1}\beta}$ while \eqref{upper_bound_Pi2} implies that $|\Pi_2^*(t,x)|\le C_{\Pi_2^*}$ for some constant $C_{\Pi_2^*}$. Coupled with the fact that $|1-(1-x/t)^\alpha |\le C^* x/t$ for $x\in [0,s]$ and $t$ large enough (say $t\ge t^*$) where $C^*$ is some constant, we can upper bound \eqref{PiDef2} as
 \begin{equation}\label{upper_bound_Pi}
|\Pi(t,x)| \le \left(C^* \frac{x}{t}\right) (C_{\Pi_1^*}+C_{\Pi_2^*}) t^{\alpha-1}\le C_\Pi t^{\alpha-2}~~\text{for all}~~x\in[0,s]~~\text{and}~~t\ge\max\{s,t^*\},
\end{equation}
where $C_\Pi:=C^*s(C_{\Pi_1^*}+C_{\Pi_2^*})$.\\

{\bf The result.} Plugging (\ref{conv_Xi}) and (\ref{equiv_Pi}) into \eqref{first_last_cov_expo_decompose} yields the asymptotic equivalent, for all $x\in [0,s]$,
\begin{equation}\label{equivalent_partial_first_last_term}
\mathbb{E}[Z(t-x)]  - \mathbb{E}[Z(t)]\sim \frac{\mu_1\lambda }{\Gamma(\alpha)}\frac{1-\alpha}{\beta} x t^{\alpha-2}\text{~~as~~}t\rightarrow \infty.
\end{equation}
We shall now establish an equivalent for the integral term in \eqref{first_last_cov_expo} and consider
\begin{equation}\label{prop4limitstep}
\frac{\int_{0}^{s} e^{\beta x} \big\{\mathbb{E}[Z(t-x)]  - \mathbb{E}[Z(t)]\big\} \mathrm{d}m(x)}{t^{\alpha-2}} =\int_{0}^{s} e^{\beta x} \left\{ \frac{\mathbb{E}[Z(t-x)]  - \mathbb{E}[Z(t)]}{t^{\alpha-2}}\right\}\mathrm{d}m(x).
\end{equation}
Utilizing the definition \eqref{first_last_cov_expo_decompose} together with the bounds \eqref{upper_bound_Xi} and \eqref{upper_bound_Pi}, it is clear that 
\[
\left| \frac{\mathbb{E}[Z(t-x)]  - \mathbb{E}[Z(t)]}{t^{\alpha-2}}\right| \le  \frac{\mu_1\lambda}{\Gamma(\alpha)}(C_\Xi+C_\Pi)~~\text{for all}~~x\in[0,s]~~\text{and}~~t\ge\max\{2s,t^*\}.
\]
Taking limit as $t\to \infty$ on both sides of \eqref{prop4limitstep} along with the use of dominated convergence and \eqref{equivalent_partial_first_last_term}, it is easy to see that
\begin{equation}\label{prop4limitstep2}
\lim_{t\to \infty}\frac{\int_{0}^{s} e^{\beta x} \big\{\mathbb{E}[Z(t-x)]  - \mathbb{E}[Z(t)]\big\} \mathrm{d}m(x)}{t^{\alpha-2}} = \frac{\mu_1\lambda }{\Gamma(\alpha)}\frac{1-\alpha}{\beta} \int_{0}^{s} e^{\beta x} x\mathrm{d}m(x).
\end{equation}

It now remains to conclude how the asymptotic formula \eqref{COVZdstasy} is obtained. Recall that the second and the third terms of \eqref{Covth} are respectively $O(e^{-(\beta+2\delta) t})$ and $O(e^{-(\beta+\delta) t})$. These are negligible compared to the sum of the first and the last terms, namely \eqref{first_last_cov_expo}, which is $O(e^{-\delta t}t^{\alpha-2})$ because of \eqref{prop4limitstep2}. Applying \eqref{first_last_cov_expo} and \eqref{prop4limitstep2} again yields the result \eqref{COVZdstasy}.\vspace{-.1in}

\subsection{Proof of auxiliary results for Proposition \ref{P6}}\label{PP6}

{\bf 1. Proving (\ref{one}) when $0<\eta<\frac{\alpha+1}{2}$.} To begin, we shall show that $\sum^\infty_{n=0} \chi^*_n(t)$ with $\chi^*_n(t)$ given by (\ref{Ant}) is uniformly convergent. It is first observed that
\begin{equation}\label{Cn}
|\chi^*_n(t)| \leq \frac{(\eta)_n}{n!(\alpha+n)}\mathrm{B}(\alpha, \alpha+n+1){}_{2}\mathrm{F}_1(2\eta+n,\alpha, 2\alpha+n+1,1):=C_n^*.
\end{equation}
Utilizing the identity (\ref{ohf2}) which requires $2\alpha+n+1>2\eta+n+\alpha$ (i.e. $\eta<\frac{\alpha+1}{2}$ as assumed), one can write
\[
C_n^* =\frac{(\eta)_n}{n!(\alpha+n)} \frac{\mathrm{\Gamma}(\alpha)\mathrm{\Gamma}(\alpha+1-2\eta)}{\mathrm{\Gamma}(2\alpha+1-2\eta)} =\frac{(\eta)_n}{n!(\alpha+n)}\mathrm{B}(\alpha,\alpha+1-2\eta).
\]
Then, via application of (\ref{ohf}) and (\ref{ohf2}) (requiring $\alpha+1>\alpha+\eta$, i.e. $\eta<1$ which is satisfied), we can sum $C_n^*$'s as
\begin{align}\label{sumCn}
\sum^\infty_{n=0}C_n^* &=\mathrm{B}(\alpha,\alpha+1-2\eta) \sum^\infty_{n=0} \frac{(\eta)_n}{n!(\alpha+n)} =\mathrm{B}(\alpha,\alpha+1-2\eta) \frac{1}{\alpha}
{}_{2}\mathrm{F}_1(\alpha,\eta,\alpha+1,1)\nonumber\\
&=\mathrm{B}(\alpha,\alpha+1-2\eta) \frac{\mathrm{\Gamma}(\alpha)\mathrm{\Gamma}(1-\eta)}{\mathrm{\Gamma}(\alpha+1-\eta)} =\mathrm{B}(\alpha,\alpha+1-2\eta) \mathrm{B}(\alpha,1-\eta) <\infty.
\end{align}
Hence, uniform convergence is proved. Then, further using the definitions of $\chi^*_n(t)$ and $C_n^*$ in (\ref{Ant}) and (\ref{Cn}), we can now evaluate the limit on the right-hand side of \eqref{one0} as
\begin{align}\label{prop6step}
\lim_{t \rightarrow \infty} \sum^\infty_{n=0} \chi^*_n(t) &=\sum^\infty_{n=0}\lim_{t \rightarrow \infty} \chi^*_n(t) \nonumber\\
&=\sum^\infty_{n=0} \frac{(\eta)_n}{n!(\alpha+n)}\mathrm{B}(\alpha, \alpha+n+1)\lim_{t \rightarrow \infty} \bigg[ \Big(\frac{t}{\theta+t}\Big)^{n}{}_{2}\mathrm{F}_1\left(2\eta+n,\alpha, 2\alpha+n+1,\frac{t}{\theta+t}\right)\bigg]\nonumber\\
&=\sum^\infty_{n=0} \frac{(\eta)_n}{n!(\alpha+n)}\mathrm{B}(\alpha, \alpha+n+1){}_{2}\mathrm{F}_1(2\eta+n,\alpha, 2\alpha+n+1,1) =\sum^\infty_{n=0}C_n^*\nonumber\\
&=\mathrm{B}(\alpha,\alpha+1-2\eta) \mathrm{B}(\alpha,1-\eta),
\end{align}
from which \eqref{one} follows. Note that ${}_{2}\mathrm{F}_1(a,b,c,z)$ is continuous at $z=1$ when $\mathrm{Re}(c-a-b)>0$ (i.e. $\eta<\frac{\alpha+1}{2}$ in our case), and \eqref{sumCn} has been applied. 

\noindent{\bf 2. Proving (\ref{Atasym}) when $\eta\ge1$.} Via a change of variable with $y = \frac{t+\theta}{t}\frac{t-x}{\theta+t-x}$ in (\ref{Atdef}), i.e. $x=(\theta +t)\big[ 1 - \frac{\theta}{\theta + t(1-y)}\big]$, so that $x\in[0,t]$ implies $y\in [0,1]$, and $\mathrm{d}x = -(\theta +t)\frac{\theta t}{[\theta+t(1-y)]^2}\mathrm{d}y$, we obtain
\begin{equation}\label{Astarproofstep}
A^*(t) = \theta^{\alpha-2\eta+1}t^{\alpha- 1} \left[\int_{0}^{1}(1-y)^{\alpha -1} \Big(1-\frac{t}{\theta+t} y\Big)^{2\eta -2\alpha -1} G^*\bigg(\frac{t}{\theta+t}y \bigg) \mathrm{d}y\right] \frac{t}{\theta+t}.
\end{equation}	
Based on the fact that $\frac{t}{\theta+t}\rightarrow 1^-$ as $t\rightarrow\infty$, we shall establish the limit of $\frac{A^*(t)}{t^{\alpha-1}}$ as $t\rightarrow\infty$. It is first noted that for all $y\in[0,1)$, $G^*(\frac{t}{\theta+t}y)$ increases towards $G^*(y)$ as $t\rightarrow \infty$. In order to determine the limit of $\frac{A^*(t)}{t^{\alpha-1}}$, we shall use the asymptotic result of $G^*(y)$ as $y\rightarrow1^-$ in Lemma \ref{L3}. In what follows, we focus on the case $\eta>1$ as the argument will be similar when $\eta=1$. We shall need the following lemma.
\begin{lemma}\normalfont\label{L4}
Given $\theta>0$. For every $z\in(0,1)$ one has the limiting result
\begin{equation*}
\lim_{t\rightarrow\infty}\int_{z}^{1}(1-y)^{\alpha -1}\Big(1-\frac{t}{\theta+t}y \Big)^{\eta-2\alpha} \mathrm{d}y =\int_{z}^{1}(1-y)^{\eta-\alpha -1}\mathrm{d}y,
\end{equation*}
which is a convergent integral for $\eta>1$ and $\alpha\in (0,1)$.

\noindent \textbf{Proof:} We distinguish two cases according to the sign of $\eta-2\alpha$.

\noindent\textbf{Case 1.} $\eta\le 2\alpha$: Because $\big(1-\frac{t}{\theta+t}y \big)^{\eta-2\alpha}$ increases towards $(1-y)^{\eta-2\alpha}$ as $t\rightarrow \infty$, the result follows by monotone convergence.

\noindent\textbf{Case 2.} $\eta>2\alpha$: For all $t\ge0$ and $y\in(z,1)$, it is clear that
\[
(1-y)^{\alpha -1}\Big(1-\frac{t}{\theta+t}y \Big)^{\eta-2\alpha} \le (1-y)^{\alpha -1}.
\]
Since $\int_{z}^{1}(1-y)^{\alpha-1}\mathrm{d}y$ is finite, application of the dominated convergence theorem yields the result.
\hfill$\square$
\end{lemma}

From (\ref{G1yasym}), one has that there exists $x_0\in(0,1)$ such that
\[
|G^*(x)|\le \frac{3}{2}\frac{1}{\eta-1}|1-x|^{-\eta+1} ~~\text{for all}~~ x \in(x_0,1).
\]
This implies that we can first fix a $z\in(x_0,1)$, and then there exists $t_0$ large enough such that
\begin{equation}\label{G1absol}
\bigg\vert G^* \bigg(\frac{t}{\theta+t}y\bigg)\bigg\vert \le \frac{3}{2}\frac{1}{\eta-1}\bigg\vert 1-\frac{t}{\theta+t}y \bigg\vert^{-\eta+1}~~\text{for all}~~ y\in(z,1)~~\text{and}~~t\ge t_0.
\end{equation}
For example, one may pick $t_0$ such that $\frac{t_0}{\theta+t_0}z = x_0$, i.e. $t_0=\frac{\theta x_0}{z-x_0}$.

Next, due to \eqref{Astarproofstep}, we turn our attention the asymptotic behaviour of
\begin{equation}\label{J1J2}
\frac{A^*(t)}{t^{\alpha-1}}\frac{1}{\theta^{\alpha-2\eta+1}}\frac{\theta+t}{t} = \int_{0}^{1}(1-y)^{\alpha-1}\Big(1-\frac{t}{\theta+t}y \Big)^{2\eta-2\alpha-1}G^*\bigg(\frac{t}{\theta+t}y \bigg) \mathrm{d}y=\int_0^z + \int_z^1:=J_1(z,t) + J_2(z,t),
\end{equation}
where $J_1(z,t)$ and $J_2(z,t)$ have obvious definitions, and $z\in(x_0,1)$ is fixed.
Starting with the term $J_1(z,t)$, it is noted that $|[1-\frac{t}{\theta+t}y]^{2\eta-2\alpha -1}G^*(\frac{t}{\theta+t}y) |$ is upper bounded by $\max\{1,(1-z)^{2\eta-2\alpha-1}\}G^*(z)$ for all $t\ge 0$ and $y\in[0,z]$ and the integral $\int_{0}^{z}(1-y)^{\alpha-1}\mathrm{d}y$ is finite. Application of the dominated convergence theorem gives rise to
\begin{equation}\label{J1ztasy}
\lim_{t \rightarrow \infty} J_1(z,t) =\int_{0}^{z}(1-y)^{\alpha-1}(1-y)^{2\eta-2\alpha-1}G^*(y)\mathrm{d}y  = \int_{0}^{z}(1-y)^{2\eta-\alpha-2}G^*(y)\mathrm{d}y.
\end{equation}
As for $J_2(z,t)$, using (\ref{G1absol}) one has the upper bound
\begin{align*}
|J_2(z,t)| & \le \frac{3}{2}\frac{1}{\eta-1}\int_{z}^{1}(1-y)^{\alpha-1}\Big(1-\frac{t}{\theta+t}y \Big)^{2\eta-2\alpha -1}\Big(1-\frac{t}{\theta+t}y \Big)^{-\eta+1}\mathrm{d}y,\\
& = \frac{3}{2}\frac{1}{\eta-1}\int_{z}^{1}(1-y)^{\alpha-1}\Big(1-\frac{t}{\theta+t}y \Big)^{\eta-2\alpha}\mathrm{d}y ~~\text{for all}~~t\ge t_0.
\end{align*}
Thus, it follows from Lemma \ref{L4} that
\[
\limsup_{t\rightarrow\infty} |J_2(z,t) |\le \frac{3}{2}\frac{1}{\eta-1}\int_{z}^{1}(1-y)^{\eta-\alpha-1}\mathrm{d}y,
\]
which is equivalent to
\begin{equation}\label{d*}
\limsup_{t\rightarrow\infty} \bigg\vert \frac{A^*(t)}{t^{\alpha-1}}\frac{1}{\theta^{\alpha-2\eta+1}}\frac{\theta+t}{t} -J_1(z,t) \bigg \vert \le \frac{3}{2}\frac{1}{\eta-1}\int_{z}^{1}(1-y)^{\eta-\alpha-1}\mathrm{d}y
\end{equation}
thanks to (\ref{J1J2}). Let us now observe that, by the triangular inequality,
\begin{align*}
&\bigg\vert \frac{A^*(t)}{t^{\alpha-1}}\frac{1}{\theta^{\alpha-2\eta+1}}\frac{\theta+t}{t} - \int_{0}^{1}(1-y)^{2\eta-\alpha-2} G^*(y)\mathrm{d}y\bigg\vert \\
& \le \bigg\vert  \frac{A^*(t)}{t^{\alpha-1}}\frac{1}{\theta^{\alpha-2\eta+1}}\frac{\theta+t}{t} -J_1(z,t) \bigg\vert  +\bigg\vert  J_1(z,t) -\int_{0}^{1}(1-y)^{2\eta-\alpha-2} G^*(y)\mathrm{d}y \bigg\vert,
\end{align*}
and therefore, with the help of (\ref{d*}) and (\ref{J1ztasy}) we arrive at
\begin{align}\label{limsupJ2}
&\limsup_{t\rightarrow\infty} \bigg\vert \frac{A^*(t)}{t^{\alpha-1}}\frac{1}{\theta^{\alpha-2\eta+1}}\frac{\theta+t}{t}- \int_{0}^{1}(1-y)^{2\eta-\alpha-2} G^*(y)\mathrm{d}y \bigg\vert \nonumber\\
&\qquad \le \frac{3}{2}\frac{1}{\eta-1}\int_{z}^{1}(1-y)^{\eta-\alpha-1}\mathrm{d}y + \int_{z}^{1}(1-y)^{2\eta-\alpha-2}G^*(y)\mathrm{d}y.
\end{align}
While the first integral on the right-hand side is clearly finite, the second integral is also convergent thanks to the estimate (\ref{G1yasym}). Since $z\in(x_0,1)$ is arbitrary, we can let $z\rightarrow 1^-$ in the above inequality, ending with the right-hand side converging to $0$. We thus obtain (\ref{Atasym}) for $\eta>1$.

Finally, we remark that the procedure to derive (\ref{Atasym}) in the case $\eta=1$ is similar. The splitting of \eqref{J1J2} into $J_1(z,t)$ and $J_2(z,t)$ can still be adopted but the major difference is that \eqref{G1absol} needs to be modified using the second case of (\ref{G1yasym}). The limit \eqref{J1ztasy} is still valid, and one can obtain a similar upper bound to \eqref{limsupJ2} which converges to $0$. The details are omitted.\vspace{-.1in}

\subsection{Proof of Proposition \ref{CovPareto}}\label{appendix_proof_cov_Pareto}
Recall that $s$ is fixed such that $0<s\le t$. As in the proof of Proposition \ref{CovExpPro} for the case of exponential delays, to analyze the asymptotic behaviour of $\mathbb{C}\mathrm{ov}[Z(s),Z(t)]$ as $t\to\infty$ we shall look at the terms appearing in the general expression (\ref{Covd0}). We can write $\mathbb{C}\mathrm{ov}[Z(s), Z(t)]  =J_1^*(t) + J_2^*(t) + J_3^*(t)$ for $0<s\leq t$,
where
\begin{align}
J_1^*(t):=&~\mu_2\int^s_0 \overline{W}(t-x)\mathrm{d}m(x),\label{Jstar1}\\
J_2^*(t):=&~\mu_1\int^s_0 \mathbb{E}[Z(s-x)] \overline{W}(t-x)\mathrm{d}m(x),\label{Jstar2}\\
J_3^*(t):=&~\mu_1\int^s_0  \overline{W}(s-x) \mathbb{E}[Z(t-x)]\mathrm{d}m(x) -\mathbb{E}[Z(s)] \mathbb{E}[Z(t)].\label{Jstar3}
\end{align}
The following is dedicated to obtaining the asymptotics for $J_1^*(t)$, $J_2^*(t)$ and $J_3^*(t)$ as $t\to \infty$.

{\bf Term $J_1^*(t)$.} Substituting the renewal density (\ref{FP}) and the Pareto survival function (\ref{Ptail}) into \eqref{Jstar1} followed by some simple manipulations, we arrive at the two sided bounds
\begin{align*}
&\frac{\mu_2\lambda s^\alpha}{\Gamma(\alpha+1)} \Big(\frac{\theta}{\theta + t }\Big)^\eta =\frac{\mu_2\lambda}{\Gamma(\alpha)} \Big(\frac{\theta}{\theta + t }\Big)^\eta \int_0^s x^{\alpha -1}\mathrm{d} x \\
&\le J_1^*(t)=\frac{\mu_2\lambda}{\Gamma(\alpha)} \int_0^s \Big(\frac{\theta}{\theta + t - x}\Big)^\eta x^{\alpha -1}\mathrm{d} x \le \frac{\mu_2\lambda}{\Gamma(\alpha)} \Big(\frac{\theta}{\theta + t -s}\Big)^\eta \int_0^s x^{\alpha -1}\mathrm{d} x =\frac{\mu_2\lambda s^\alpha}{\Gamma(\alpha+1)} \Big(\frac{\theta}{\theta + t -s}\Big)^\eta.
\end{align*}
By squeezing principle one obtains the asymptotic behaviour
\begin{equation}\label{J_1_Pareto}
J_1^*(t)\sim \frac{\mu_2\lambda \theta^\eta s^\alpha}{\Gamma(\alpha+1)} t^{-\eta}~~\text{as}~~t\rightarrow\infty.
\end{equation}

{\bf Term $J_2^*(t)$.} Again, putting (\ref{FP}) and (\ref{Ptail}) into \eqref{Jstar2} yields, in the same spirit as (\ref{J_1_Pareto}),
\begin{equation}\label{J_2_Pareto}
J_2^*(t)\sim  \frac{\mu_1\lambda \theta^\eta}{\Gamma(\alpha)} \left\{\int_0^s \mathbb{E}[Z(s-x)]x^{\alpha-1}\mathrm{d} x\right\}t^{-\eta}~~\text{as}~~t\rightarrow\infty.
\end{equation}

{\bf Term $J_3^*(t)$.} Utilizing \eqref{EZdtdelta0}, one sees that $J_3^*(t)$ in \eqref{Jstar3} can be written as
\begin{equation}\label{expression_J_3}
J_3^*(t) = \mu_1\int^s_0  \overline{W}(s-x) \big\{\mathbb{E}[Z(t-x)]- \mathbb{E}[Z(t)]\big\} \mathrm{d}m(x).
\end{equation}
In order to derive the asymptotics for $J_3^*(t)$ as $t\to\infty$, we proceed by studying the asymptotic behaviour of $\mathbb{E}[Z(t-x)]- \mathbb{E}[Z(t)]$ for all $x\in [0,s]$ and then apply dominated convergence. Let us first define, for $v\ge0$,
\[
\varphi(v):=\theta^\eta\int_0^1 \frac{1}{(v+1-z)^\eta} z^{\alpha -1}\mathrm{d} z,
\]
which is decreasing in $v$. By changing variable $y=z/(v+1)$, we get that
\begin{align}\label{expression_phi}
\varphi(v)=&~\frac{\theta^\eta}{(v+1)^\eta}\int_0^1 \frac{1}{[1-z/(v+1)]^\eta} z^{\alpha -1}\mathrm{d} z =\theta^\eta(v+1)^{\alpha-\eta} \int_0^{\frac{1}{v+1}} (1-y)^{-\eta}y^{\alpha-1}\mathrm{d} y \nonumber\\
=&~\theta^\eta (v+1)^{\alpha-\eta} G^*\left(\frac{1}{v+1}\right),
\end{align}
where $G^*(\cdot)$ is defined in (\ref{G1ydef}) (with the definition extended to allow for $0<\eta<1$ as well). Hence, one has $\mathbb{E}[Z(t)]=\frac{\mu_1\lambda }{\Gamma(\alpha)}t^{\alpha-\eta}\varphi\left(\frac{\theta}{t}\right)$ thanks to (\ref{expr_EZ_phi}). Defining
\begin{equation}\label{def_psi}
\psi(t,x):=\varphi\left(\frac{\theta}{t-x}\right)-\varphi\left(\frac{\theta}{t}\right)
\end{equation}
\noindent for $x\in[0,s]$, we can write
\[
\mathbb{E}[Z(t-x)]- \mathbb{E}[Z(t)] =\frac{\mu_1\lambda }{\Gamma(\alpha)}\left\{(t-x)^{\alpha-\eta}\psi(t,x)+[(t-x)^{\alpha-\eta}-t^{\alpha-\eta}]\varphi\left(\frac{\theta}{t}\right)\right\},
\]
and then $J_3^*(t)$ in \eqref{expression_J_3} can be decomposed as $J_3^*(t)=J_{31}^*(t)+J_{32}^*(t)$ where
\begin{align}
J_{31}^*(t):=&~\frac{\mu_1^2\lambda }{\Gamma(\alpha)} \int^s_0  \overline{W}(s-x) (t-x)^{\alpha-\eta}\psi(t,x) \mathrm{d} m(x), \label{def_S1}\\
J_{32}^*(t):=&~\frac{\mu_1^2\lambda }{\Gamma(\alpha)} \varphi\left(\frac{\theta}{t}\right) \int^s_0  \overline{W}(s-x)[(t-x)^{\alpha-\eta}-t^{\alpha-\eta}] \mathrm{d}m(x).  \label{def_S2}
\end{align}
The terms $J_{31}^*(t)$ and $J_{32}^*(t)$ are studied separately as follows.

{\bf Term $J_{31}^*(t)$ of $J_3^*(t)$.} We shall first proceed to prove the auxiliary result
\begin{equation}\label{lim_psi_prime}
\lim_{t\to \infty}\frac{\psi(t,x)}{t^{\eta-2}}=-\theta x~~\text{for all}~~x\in[0,s].
\end{equation}
Note that the limiting behaviour of $t^{\eta-2}$ as $t\rightarrow\infty$ depends on the value of $\eta$. It will be seen that L'H\^opital's rule can be applied when $\eta\ne 2$ where the limit on the left-hand side of \eqref{lim_psi_prime} is of the indeterminate form $0/0$ or $\infty/\infty$. This involves differentiating $\psi(x,t)$ in \eqref{def_psi} with respect to $t$, yielding
\begin{equation}\label{compute_psi_prime}
\frac{\partial}{\partial t}\psi(t,x) = - \frac{\theta}{(t-x)^2}\varphi ' \left( \frac{\theta}{t-x}\right) + \frac{\theta}{t^2} \varphi ' \left( \frac{\theta}{t}\right)=\left[ - \frac{\theta}{(t-x)^2}+\frac{\theta}{t^2}\right]\varphi ' \left( \frac{\theta}{t-x}\right)+ \frac{\theta}{t^2} \left[ \varphi ' \left( \frac{\theta}{t}\right) - \varphi ' \left( \frac{\theta}{t-x}\right)\right],
\end{equation}
where the derivative of $\varphi(\cdot)$ in \eqref{expression_phi} is given by
\begin{equation*}
\varphi ' (v) = - \theta^\eta v^{-\eta}(v+1)^{-1}+\theta^\eta (\alpha-\eta)(v+1)^{\alpha-\eta-1}G^*\left(\frac{1}{v+1}\right).
\end{equation*}
Four cases are considered as follows.

\noindent\textbf{Case 1.} $\eta\in (1,2)\cup(2,\infty)$: 
Utilizing \eqref{G1yasym} in Lemma \ref{L3}, it is not difficult to see that $G^*\big(\frac{1}{v+1}\big)= \frac{1}{\eta-1}v^{-\eta+1} + o(v^{-\eta+1})$ as $v\to 0^+$, and consequently one has the expansion
\begin{align}\label{equivalent_phi_prime}
\varphi ' (v)=& -\theta^\eta v^{-\eta} +\theta^\eta v^{-\eta}[1-(v+1)^{-1}] +\theta^\eta (\alpha-\eta)(v+1)^{\alpha-\eta-1}\left[\frac{1}{\eta-1}v^{-\eta+1} + o(v^{-\eta+1})\right]\nonumber\\
=&-\theta^\eta v^{-\eta} +\theta^\eta v^{-\eta+1} +\theta^\eta \frac{\alpha-\eta}{\eta-1} v^{-\eta+1} + o(v^{-\eta+1}) \nonumber\\
=&-\theta^\eta v^{-\eta} -\theta^\eta \frac{1-\alpha}{\eta-1} v^{-\eta+1} + o(v^{-\eta+1}) ~~\text{as}~~v\to 0^+.
\end{align}
The first term on the right-hand side of \eqref{compute_psi_prime} thus satisfies the asymptotic formula
\begin{align*}
\left[ - \frac{\theta}{(t-x)^2}+\frac{\theta}{t^2}\right]\varphi ' \left( \frac{\theta}{t-x}\right) \sim &  - \left[ - \frac{\theta}{(t-x)^2}+\frac{\theta}{t^2}\right] \theta^\eta \left(\frac{\theta}{t}\right)^{-\eta} = \theta t^{\eta-2}\left[ (1-x/t)^{-2}-1\right]\\
\sim &~ 2\theta x t^{\eta-3}~~\text{as}~~t\to\infty.
\end{align*}
Since $\frac{1}{t^2} \left[t^\eta - (t-x)^\eta\right]\sim\eta x t^{\eta-3}$ and $\frac{1}{t^2} \left[t^{\eta-1} - (t-x)^{\eta-1}\right]\sim(\eta-1) x t^{\eta-4}=o(t^{\eta-3})$ as $t\to\infty$, the second term on the right-hand side of (\ref{compute_psi_prime}) can again be estimated using (\ref{equivalent_phi_prime}) via
\begin{align*}
\frac{\theta}{t^2} \left[ \varphi ' \left( \frac{\theta}{t}\right) - \varphi ' \left( \frac{\theta}{t-x}\right)\right]= & - \frac{\theta}{t^2} \left[ t^\eta - (t-x)^\eta \right]  - \frac{\theta^2}{t^2} \frac{1-\alpha}{\eta-1} \left[ t^{\eta -1} - (t-x)^{\eta -1} \right]+o(t^{\eta -3})\\
=& -\eta \theta x t^{\eta-3} + o(t^{\eta-3})~~\text{as}~~t\to\infty.
\end{align*}
Plugging the above two estimates into (\ref{compute_psi_prime}) yields, for all $x\in[0,s]$,
\begin{equation}\label{lhopital}
\frac{\partial}{\partial t}\psi(t,x) \sim - (\eta-2) \theta xt^{\eta-3}~~\text{as}~~t\to\infty.
\end{equation}
This is further split into two cases as follows.

\noindent\textbf{Case 1a.} $\eta\in (2,\infty)$: In view of the expression on the right-hand side of \eqref{lhopital}, it is noted that $\int(\eta-2) \theta xt^{\eta-3}\mathrm{d}t=\theta x t^{\eta-2}\to\infty$ as $t\to\infty$, asserting that $\psi(t,x)\to\infty$ as $t\to\infty$. As a result, the left-hand side of \eqref{lim_psi_prime} is in the form $\infty/\infty$, and application of L'H\^opital's rule yields
\begin{equation}\label{lhopital2}
\lim_{t\to \infty}\frac{\psi(t,x)}{t^{\eta-2}}=\lim_{t\to \infty}\frac{\frac{\partial}{\partial t}\psi(t,x) }{(\eta-2)t^{\eta-3}}~~\text{for all}~~x\in[0,s],
\end{equation}
which proves \eqref{lim_psi_prime} thanks to \eqref{lhopital}.

\noindent\textbf{Case 1b.} $\eta\in (1,2)$: Using the finer asymptotic result \eqref{G1yasymfine}, one finds that \eqref{expression_phi} with $v=\theta/t$ can be expressed as
\begin{equation*}
\varphi\left(\frac{\theta}{t}\right)=\theta^\eta \left(\frac{\theta}{t}+1\right)^{\alpha-\eta} \bigg[ \frac{1}{\eta-1}\left(\frac{\theta/t}{\theta/t+1}\right)^{-\eta+1}+C^*+o(1)\bigg]~~\text{as}~~t\to\infty
\end{equation*}
for some constant $C^*$. Further expanding $\big(\frac{\theta}{t}+1\big)^a=1+a\frac{\theta}{t}+o\big(\frac{1}{t}\big)$ at $a=\alpha-\eta$ and $a=\eta-1$, we obtain
\begin{align}\label{prop7case1bvarphid}
\varphi\left(\frac{\theta}{t}\right)&=\theta^\eta \left[1+(\alpha-\eta)\frac{\theta}{t}+o\left(\frac{1}{t}\right)\right] \bigg\{ \frac{1}{\eta-1}\left(\frac{\theta}{t}\right)^{-\eta+1} \left[1+(\eta-1)\frac{\theta}{t}+o\left(\frac{1}{t}\right)\right] +C^*+o(1)\bigg\}\nonumber\\
&=\theta^\eta \left[1+(\alpha-\eta)\frac{\theta}{t}+o\left(\frac{1}{t}\right)\right] \left[ \frac{1}{\eta-1}\left(\frac{\theta}{t}\right)^{-\eta+1} +\left(\frac{\theta}{t}\right)^{-\eta+2} +o\left(\frac{1}{t^{-\eta+2}}\right) +C^*+o(1)\right]\nonumber\\
&=\theta^\eta \bigg[\frac{1}{\eta-1}\left(\frac{\theta}{t}\right)^{-\eta+1} +C^*+o(1)\bigg] ~~\text{as}~~t\to\infty,
\end{align}
where the last line follows from the assumption $-\eta+2>0$. By noting that (again because of $-\eta+2>0$)
\begin{equation*}
\left(\frac{\theta}{t-x}\right)^{-\eta+1} =\left(\frac{\theta}{t}\right)^{-\eta+1} \left(1-\frac{x}{t}\right)^{\eta-1} =\left(\frac{\theta}{t}\right)^{-\eta+1} \left[1-(\eta-1)\frac{x}{t}+o\left(\frac{1}{t}\right)\right] =\left(\frac{\theta}{t}\right)^{-\eta+1}+o(1)~~\text{as}~~t\to\infty,
\end{equation*}
replacing $t$ by $t-x$ in \eqref{prop7case1bvarphid} gives rise to
\begin{equation*}
\varphi\left(\frac{\theta}{t-x}\right) =\theta^\eta \bigg[\frac{1}{\eta-1}\left(\frac{\theta}{t}\right)^{-\eta+1} +C^*+o(1)\bigg] ~~\text{as}~~t\to\infty.
\end{equation*}
In view of the definition \eqref{def_psi}, taking difference of the above equation with \eqref{prop7case1bvarphid} results in $\psi(t,x)=o(1)$ as $t\to\infty$, i.e. the left-hand side of \eqref{lim_psi_prime} is in the form $0/0$. Consequently, the L'H\^opital's rule can be applied as in \eqref{lhopital2}, and further use of \eqref{lhopital} leads to \eqref{lim_psi_prime}.

\noindent\textbf{Case 2.} $\eta\in(0,1)$:  In this case, it is clear that $G^*\big(\frac{1}{v+1}\big)=\mathrm{B}\big(\alpha,1-\eta,\frac{1}{v+1}\big)=\mathrm{B}(\alpha,1-\eta)+o(1)$ as $v\to 0^+$. Therefore, it is clear from \eqref{expression_phi} and \eqref{def_psi} that $\psi(t,x)\rightarrow0$ as $t\to\infty$, implying that the left-hand side of \eqref{lim_psi_prime} is $0/0$. Note that (\ref{equivalent_phi_prime}) is now replaced by $\varphi'(v)= -\theta^\eta v^{-\eta}+\theta^\eta (\alpha-\eta) \mathrm{B}(\alpha,1-\eta) + o(1)$, and one can check that the arguments leading to \eqref{lhopital} and \eqref{lim_psi_prime} still hold.

\noindent\textbf{Case 3.} $\eta=1$: According to \eqref{G1yasym} in Lemma \ref{L3}, one has $G^*\big(\frac{1}{v+1}\big)=-\ln v + C+o(1)$ as $v\to 0^+$ for some constant $C$. With some simple algebra it can be shown that \eqref{equivalent_phi_prime} is replaced by $\varphi'(v)= - \theta v^{-1}+\theta(1-\alpha) \ln v +C^* + o(1)$ as $v\to 0^+$ where $C^*$ is a constant. Note that the first term is the dominant one, and the procedure leading to \eqref{lhopital} is still valid. It remains to show that $\psi(t,x)\rightarrow0$ as $t\to\infty$ so that L'H\^opital's rule in \eqref{lhopital2} can be applied and \eqref{lim_psi_prime} holds true. To see this, we substitute the above estimation of $G^*\big(\frac{1}{v+1}\big)$ into \eqref{expression_phi} (under $v=\theta/t$) so that, similar to \eqref{prop7case1bvarphid},
\begin{equation*}
\varphi\left(\frac{\theta}{t}\right)=\theta \left(\frac{\theta}{t}+1\right)^{\alpha-1} \left[-\ln \left(\frac{\theta}{t}\right) + C+o(1)\right] =\theta [\ln t + C^*+o(1)]~~\text{as}~~t\to\infty
\end{equation*}
for some constant $C^*$. One also has the same equation but with $t$ replaced by $t-x$, and therefore taking difference results in
\[
\psi(t,x) =\theta \left[\ln \left(\frac{t-x}{t}\right) + o(1)\right]=o(1)~~\text{as}~~t\to\infty.
\]

\noindent\textbf{Case 4.} $\eta=2$: This final case is different from the previous ones in the sense that $t^{\eta-2}$ in the denominator on the left-hand side of \eqref{lim_psi_prime} equals to one. We shall proceed directly to calculate the limit $\lim_{t\to \infty}\psi(t,x)$ as follows. First, application of the finer asymptotic result \eqref{G1yasymfine} to \eqref{expression_phi} with $v=\theta/t$ yields, for some constant $C^{**}$,
\begin{align}\label{varphicase4}
\varphi\left(\frac{\theta}{t}\right)=&~\theta^2 \left(\frac{\theta}{t}+1\right)^{\alpha-2} \bigg[ \left(\frac{\theta/t}{\theta/t+1}\right)^{-1}-(1-\alpha)\ln\left(\frac{\theta/t}{\theta/t+1}\right)+C^{**}+o(1)\bigg]\nonumber\\
=&~\theta^2 \bigg[ \left(\frac{\theta}{t}\right)^{-1} \left(\frac{\theta}{t}+1\right)^{\alpha-1} -(1-\alpha)\left(\frac{\theta}{t}+1\right)^{\alpha-2} \ln\left(\frac{\theta}{t}\right) +(1-\alpha)\left(\frac{\theta}{t}+1\right)^{\alpha-2} \ln\left(\frac{\theta}{t}+1\right)\nonumber\\
&+C^{**}\left(\frac{\theta}{t}+1\right)^{\alpha-2}+o(1)\bigg]~~\text{as}~~t\to\infty.
\end{align}
We further expand each term in the square bracket as
\begin{align*}
\left(\frac{\theta}{t}\right)^{-1} \left(\frac{\theta}{t}+1\right)^{\alpha-1} =&\left(\frac{\theta}{t}\right)^{-1}+(\alpha-1)+o(1),\\
\left(\frac{\theta}{t}+1\right)^{\alpha-2} \ln\left(\frac{\theta}{t}\right) =&~\ln\left(\frac{\theta}{t}\right) +(\alpha-2)\frac{\theta}{t}\ln\left(\frac{\theta}{t}\right) +o\left(\frac{\theta}{t}\ln\left(\frac{\theta}{t}\right)\right) =\ln\left(\frac{\theta}{t}\right) +o(1),\\
\left(\frac{\theta}{t}+1\right)^{\alpha-2} \ln\left(\frac{\theta}{t}+1\right)=&~o(1),\\
C^{**}\left(\frac{\theta}{t}+1\right)^{\alpha-2}=&~C^{**}+o(1).
\end{align*}
Therefore, \eqref{varphicase4} is reduced to
\[
\varphi\left(\frac{\theta}{t}\right) =\theta^2 \bigg[ \left(\frac{\theta}{t}\right)^{-1} -(1-\alpha)\ln\left(\frac{\theta}{t}\right) +C^{***}+o(1)\bigg]=\theta^2 \bigg[ \frac{t}{\theta} +(1-\alpha)\ln t +C^{****}+o(1)\bigg]~~\text{as}~~t\to\infty,
\]
where $C^{***}$ and $C^{****}$ are constants. The equation is also valid when $t$ is replaced by $t-x$, and therefore taking difference results in
\[
\psi(t,x) =\theta^2 \left[-\frac{x}{\theta} +(1-\alpha)\ln\left(\frac{t-x}{t}\right) + o(1)\right]=-\theta x+o(1)~~\text{as}~~t\to\infty,
\]
which is \eqref{lim_psi_prime}.

Having proved \eqref{lim_psi_prime} for $\eta>0$, we are ready to conclude on $J_{31}^*(t)$. For all $x\in[0,s]$ and $t\ge 2s$, two observations are made as follows. First, it is clear that $|(1-x/t)^{\alpha-\eta}|\le \max(1,1/2^{\alpha-\eta})$. Second, the fact that $\varphi(v)$ is decreasing in $v$ implies $|\psi(t,x)/t^{\eta-2}|\le |\psi(t,s)/t^{\eta-2}|$, which is upper bounded as it is continuous and convergent to a finite limit as $t\to \infty$ thanks to \eqref{lim_psi_prime}. Using the definition \eqref{def_S1}, we can thus apply the dominated convergence theorem to obtain
\begin{equation}\label{equivalent_S1_t}
\frac{J_{31}^*(t)}{t^{\alpha-2}} =\frac{\mu_1^2\lambda }{\Gamma(\alpha)}\int^s_0  \overline{W}(s-x) \left(1-\frac{x}{t}\right)^{\alpha-\eta}\frac{\psi(t,x)}{t^{\eta-2}} \mathrm{d} m(x) \longrightarrow -\frac{\mu_1^2\lambda \theta}{\Gamma(\alpha)}\int^s_0  \overline{W}(s-x) x \mathrm{d} m(x)~~\text{as}~~t\to\infty.
\end{equation}

{\bf Term $J_{32}^*(t)$ of $J_3^*(t)$.} Note that $J_{32}^*(t)$ defined in (\ref{def_S2}) is identical to zero if $\eta=\alpha$, and therefore it is sufficient to focus on the case $\eta\ne\alpha$. For all $x\in[0,s]$ and $t\ge s$, one obtains easily that $\frac{1}{t^{\alpha-\eta-1}}[(t-x)^{\alpha-\eta}-t^{\alpha-\eta}]\rightarrow (\eta-\alpha)x$ as $t\to\infty$. Moreover, one has the inequality $|\frac{1}{t^{\alpha-\eta-1}}[(t-x)^{\alpha-\eta}-t^{\alpha-\eta}]|\le |\frac{1}{t^{\alpha-\eta-1}}[(t-s)^{\alpha-\eta}-t^{\alpha-\eta}]|$ which is an upper bounded quantity. The dominated convergence theorem thus yields from (\ref{def_S2}) the equivalent
\begin{equation}\label{equivalent1_S2_t}
J_{32}^*(t) \sim \frac{\mu_1^2\lambda(\eta-\alpha)}{\Gamma(\alpha)} \varphi\left(\frac{\theta}{t}\right) \left[\int^s_0  \overline{W}(s-x)x \mathrm{d}m(x)\right] t^{\alpha-\eta-1}~~\text{as}~~t\to\infty.
\end{equation}
An asymptotic formula of $J_{32}^*(t)$ will follow once we find one for $\varphi\big(\frac{\theta}{t}\big)$. Recall that (\ref{expression_phi}) implies $\varphi\big(\frac{\theta}{t}\big) \sim \theta^\eta G^*\big(\frac{1}{\theta/t+1}\big)$ as $t\to\infty$, where $G^*\big(\frac{1}{\theta/t+1}\big) $ is asymptotically equivalent to $\frac{1}{\eta-1}\big(\frac{\theta}{t}\big)^{-\eta+1}$ (resp. $\ln t$) when $\eta>1$ (resp. $\eta=1$) according to \eqref{G1yasym} in Lemma \ref{L3}, and tends to $\mathrm{B}(\alpha,1-\eta)$ when $0<\eta<1$. Consolidating these results with \eqref{equivalent1_S2_t}, we arrive at the asymptotic behaviour for $J_{32}^*(t)$ given by, as $t\to\infty$,
\begin{equation}\label{equivalent_S2_t}
J_{32}^*(t) \sim
\begin{cases}
\displaystyle \frac{\mu_1^2\lambda\theta^\eta(\eta-\alpha)}{\Gamma(\alpha)} \mathrm{B}(\alpha,1-\eta) \left[\int^s_0  \overline{W}(s-x)x \mathrm{d}m(x)\right] t^{\alpha-\eta-1}, & \mathrm{if}~~0<\eta<1,\vspace{.07in}\\
\displaystyle \frac{\mu_1^2\lambda\theta (1-\alpha)}{\Gamma(\alpha)} \left[\int^s_0  \overline{W}(s-x)x \mathrm{d}m(x)\right] t^{\alpha-2}\ln t, & \mathrm{if}~~\eta=1,\vspace{.07in}\\
\displaystyle \frac{\mu_1^2\lambda\theta (\eta-\alpha)}{(\eta -1)\Gamma(\alpha)} \left[\int^s_0 \overline{W}(s-x)x \mathrm{d}m(x)\right] t^{\alpha-2}, & \mathrm{if}~~\eta>1.
\end{cases}
\end{equation}

{\bf Asymptotics for $J_3^*(t)$.} We can now gather the asymptotic results (\ref{equivalent_S1_t}) and (\ref{equivalent_S2_t}) to give the asymptotics for $J_3^*(t)$ as $t\to\infty$. First, if $\eta=\alpha$, then only $J_{31}^*(t)$ matters as $J_{32}^*(t)$ equals zero. Second, if $\eta\in(0,\alpha)\cup(\alpha,1]$, then $J_{32}^*(t)$ asymptotically dominates $J_{31}^*(t)$. Third, if $\eta>1$, then both $J_{31}^*(t)$ and $J_{32}^*(t)$ asymptotically behave like $t^{\alpha-2}$ and therefore $J_3^*(t)\sim J_{31}^*(t)+J_{32}^*(t)$. These are all summed up to yield
\begin{equation}\label{equivalent_J3_t}
J_3^*(t) \sim
\begin{cases}
\displaystyle -\frac{\mu_1^2\lambda \theta}{\Gamma(\alpha)} \left[\int^s_0  \overline{W}(s-x) x \mathrm{d} m(x)\right]t^{\alpha-2}, & \mathrm{if}~~\eta=\alpha,\vspace{.07in}\\
\displaystyle \frac{\mu_1^2\lambda\theta^\eta(\eta-\alpha)}{\Gamma(\alpha)} \mathrm{B}(\alpha,1-\eta) \left[\int^s_0  \overline{W}(s-x)x \mathrm{d}m(x)\right] t^{\alpha-\eta-1}, & \mathrm{if}~~\eta\in(0,\alpha)\cup(\alpha,1),\vspace{.07in}\\
\displaystyle \frac{\mu_1^2\lambda\theta (1-\alpha)}{\Gamma(\alpha)} \left[\int^s_0  \overline{W}(s-x)x \mathrm{d}m(x)\right] t^{\alpha-2}\ln t, & \mathrm{if}~~\eta=1,\vspace{.07in}\\
\displaystyle \frac{\mu_1^2\lambda\theta (1-\alpha)}{(\eta -1)\Gamma(\alpha)} \left[\int^s_0 \overline{W}(s-x)x \mathrm{d}m(x)\right] t^{\alpha-2}, & \mathrm{if}~~\eta>1.
\end{cases}
\end{equation}

{\bf The result.} From (\ref{J_1_Pareto}) and (\ref{J_2_Pareto}), it is clear that $J_1^*(t)$ and $J_2^*(t)$ are both asymptotically proportional to $t^{-\eta}$ as $t\to\infty$. By comparing with (\ref{equivalent_J3_t}), one sees that $J_3^*(t)$ is dominated by $J_1^*(t)+J_2^*(t)$ when $\eta< 2-\alpha$, so that $\mathbb{C}\mathrm{ov}[Z(s), Z(t)]\sim J_1^*(t)+J_2^*(t)$. If $\eta=2-\alpha$, then $J_1^*(t)$, $J_2^*(t)$ and $J_3^*(t)$ all asymptotically behave like $t^{\alpha-2}$, and one has $\mathbb{C}\mathrm{ov}[Z(s), Z(t)]\sim J_1^*(t)+J_2^*(t)+J_3^*(t)$. Finally, $J_3^*(t)$ is the dominant term when $\eta>2-\alpha$. All these results are summarized in (\ref{asymp_cov_Pareto}).

\section*{Supplementary materials}

\subsection*{Part 1: Proving the validity of interchanging limit and infinite summation in Equation (4.23)
in Proposition 3 
}

In the proof of Proposition 3 
, we have determined the constant $C$ in Equation (4.23), namely
\begin{equation}\label{limitsum}
\lim_{t\rightarrow\infty} \frac{\sum^\infty_{n=0} A_n t^{2\alpha+n} e^{-2\beta t} {}_1\mathrm{F}_1 (\alpha,2\alpha+n+1,2\beta t)}{t^{\alpha-1}}=C
\end{equation}
(where $A_n$ is defined in (4.18)
) by assuming that the limit and the infinite summation on the left-hand side can be interchanged. Here we shall prove that this is valid. Setting
\begin{equation}\label{chinDef}
\chi_n(t) =A_n t^{\alpha+n+1} e^{-2\beta t} {}_1\mathrm{F}_1 (\alpha,2\alpha+n+1,2\beta t),
\end{equation}
the left-hand side of \eqref{limitsum} equals $\lim_{t\rightarrow\infty}\sum^\infty_{n=0}\chi_n(t)$. To interchange the order of limit and infinite summation, a sufficient condition is that the series $\sum^\infty_{n=0}\chi_n(t)$ converges uniformly on $t\in [0,\infty)$. According to the Weierstrass M-test, uniform convergence of $\sum^\infty_{n=0}\chi_n(t)$ can be proved by finding a sequence $\{V_n\}_{n=0}^\infty$ such that
\begin{equation}\label{chinbound}
|\chi_n(t)|\le V_n~~\text{for all}~~n\in\mathbb{N}\cup \{ 0\}~~\text{and}~~t\ge0,
\end{equation}
and
\begin{equation}\label{Vnbound}
\sum^\infty_{n=0}V_n<\infty.
\end{equation}

Upper bounding $|\chi_n(t)|$ by a well defined $V_n$ is, in view of \eqref{chinDef}, only possible if one is able to provide some fine asymptotics for the Kummer's function $_1\mathrm{F}_1(a,b,z)$ as $z\to \infty$ (see (4.10)
). We need some precise information on the term $O(|z|^{-1})$, and in particular it is important to analyze how $O(|z|^{-1})$ depends on the parameter $b$ since $b=2\alpha+n+1$ varies with $n$. Refined expansions are only available in the literature for the Tricomi's functions (to the best of our knowledge) but not for $_1\mathrm{F}_1(a,b,z)$. Hence, we will use those expansions coupled with the identity (see Olver et al. (2006, Equation (13.2.41)))
\begin{equation}\label{Identity_Kummer_Tricomi}
_1\mathrm{F}_1 (a,b,z) =\frac{\Gamma(b)}{\Gamma(b-a)}\,e^{\mp\mathrm{i}\pi a}\, \mathrm{U}(a,b,z) +\frac{\Gamma(b)}{\Gamma(a)}\,e^{\pm\mathrm{i}\pi(b-a)}\,e^z\, \mathrm{U}(b-a,b,e^{\pm\mathrm{i}\pi}z),
\end{equation}
where $\mathrm{U}$ is the Tricomi's function (or the confluent hypergeometric function of the second kind). We recall that $\mathrm{U}(a,b,z)$ is a multivalued complex function, and hence the term $e^{\pm\mathrm{i}\pi}$, multiplied by $z$, changes its argument by $\pm\mathrm{i}\pi $ and then changes the value of the term $\mathrm{U}(b-a,b,e^{\pm\mathrm{i}\pi}z)$.

We first decompose \eqref{chinDef} into two parts as
\begin{equation}\label{chindecompose}
\chi_n(t) =\chi_{n,1}(t)+\chi_{n,2}(t),
\end{equation}
where
\begin{align}
\chi_{n,1}(t) &:=A_n t^{\alpha+n+1} e^{-2\beta t} {}_1\mathrm{F}_1 (\alpha,2\alpha+n+1,2\beta t) I_{\{2\beta t\le \kappa(n+1)\}}\label{chin1Def},\\
\chi_{n,2}(t) &:=A_n t^{\alpha+n+1} e^{-2\beta t} {}_1\mathrm{F}_1 (\alpha,2\alpha+n+1,2\beta t) I_{\{2\beta t>\kappa(n+1)\}}.\label{chin2Def}
\end{align}
In the above definitions, $\kappa>0$ is a constant that will be chosen `large enough' later on. In other words, we break down $\chi_n(t)$ into two terms $\chi_{n,1}(t)$ and $\chi_{n,2}(t)$ according to whether $t$ grows slower or faster than $O(n)$ towards $\infty$, and the terms $ \chi_{n,1}(t)$ and $\chi_{n,2}(t)$ will be studied separately.

{\bf Term $\chi_{n,1}(t)$.} First, we look at the term $\chi_{n,1}(t)$. Using (4.18) 
and (4.8)
, \eqref{chin1Def} can be rewritten as
\begin{equation}\label{chin1Eq2}
\chi_{n,1}(t) =\frac{\alpha}{\alpha+n}\,\frac{1}{n!}\,\frac{1}{(2\beta)^\alpha}\,\frac{1}{2^n} \left[\int_0^1 e^{-2\beta ty} (2\beta ty)^{\alpha+n}(1-y)^{\alpha-1}\,\mathrm{d}y\right] tI_{\{2\beta t\le \kappa(n+1)\}}.
\end{equation}
It can be readily checked that $e^{-x}x^{\alpha+n}$ (as a function of $x$ for $x\ge0$) achieves maximum at $x=\alpha+n$. Therefore, one has that
$
e^{-2\beta ty} (2\beta ty)^{\alpha+n} \le e^{-(\alpha+n)} (\alpha+n)^{\alpha+n}.
$
Further noting that
$
tI_{\{2\beta t\le \kappa(n+1)\}} \le\frac{\kappa}{2\beta} (n+1),
$
\eqref{chin1Eq2} can be upper bounded as
\begin{align}\label{chin1bound}
|\chi_{n,1}(t)| \le &~\frac{\alpha}{\alpha+n}\,\frac{1}{n!}\,\frac{1}{(2\beta)^\alpha}\,\frac{1}{2^n} \left[\int_0^1 (1-y)^{\alpha-1}\,\mathrm{d}y\right] e^{-(\alpha+n)} (\alpha+n)^{\alpha+n} \frac{\kappa}{2\beta} (n+1)\nonumber\\
=&~\frac{\kappa}{(2\beta)^{\alpha+1}}\,\frac{n+1}{\alpha+n}\,\frac{1}{2^n}\,\frac{(\alpha+n)^{\alpha+n}}{n!}\,e^{-(\alpha+n)}:=~V_{n,1},
\end{align}
which does not depend on $t$. Utilizing the Stirling's formula
\begin{equation}\label{stirling1}
n!\sim\sqrt{2\pi n}\left(\frac{n}{e}\right)^n~~\text{as}~~n\rightarrow\infty,
\end{equation}
one observes that
\begin{equation}\label{asymp1}
\frac{(\alpha+n)^{\alpha+n}}{n!}\,e^{-(\alpha+n)} \sim \frac{1}{\sqrt{2\pi n}}\left(\frac{\alpha+n}{n}\right)^n (\alpha+n)^\alpha e^{-\alpha} \sim \frac{1}{\sqrt{2\pi n}}\,n^\alpha = \frac{1}{\sqrt{2\pi}}\,n^{\alpha-\frac{1}{2}} ~~\text{as}~~n\rightarrow\infty,
\end{equation}
and therefore
\begin{equation*}
V_{n,1} \sim \frac{\kappa}{(2\beta)^{\alpha+1}\sqrt{2\pi}}\,\frac{n+1}{\alpha+n}\,\frac{n^{\alpha-\frac{1}{2}}}{2^n} \sim \frac{\kappa}{(2\beta)^{\alpha+1}\sqrt{2\pi}}\,\frac{n^{\alpha-\frac{1}{2}}}{2^n}:=V_{n,1}^*~~\text{as}~~n\rightarrow\infty.
\end{equation*}
Then, we apply the ratio test to check that
\begin{equation*}
\lim_{n\rightarrow\infty}\left|\frac{V_{n+1,1}}{V_{n,1}}\right| =\lim_{n\rightarrow\infty}\frac{V_{n+1,1}^*}{V_{n,1}^*} =\lim_{n\rightarrow\infty}\frac{1}{2}\left(\frac{n+1}{n}\right)^{\alpha-\frac{1}{2}} =\frac{1}{2}<1,
\end{equation*}
and thus
\begin{equation}\label{Vn1bound}
\sum^\infty_{n=0}V_{n,1}<\infty.
\end{equation}

{\bf Term $\chi_{n,2}(t)$.} Next, applying (\ref{Identity_Kummer_Tricomi}) and (4.18) 
in \eqref{chin2Def}, we further decompose $\chi_{n,2}(t)$ as
\begin{equation}\label{chin2decompose}
\chi_{n,2}(t) =\chi_{n,21}(t)+\chi_{n,22}(t),
\end{equation}
where
\begin{equation}\label{chin21Def}
\chi_{n,21}(t) :=I_{\{2\beta t>\kappa(n+1)\}} t^{\alpha+n+1} e^{-2\beta t}\frac{\alpha}{\alpha+n}\,\frac{\beta^n }{n!}\,\Gamma(\alpha)\,e^{\mp\mathrm{i}\pi \alpha}\, \mathrm{U}(\alpha, 2\alpha+n+1,2\beta t)
\end{equation}
and
\begin{equation}\label{chin22Def}
\chi_{n,22}(t) :=I_{\{2\beta t>\kappa(n+1)\}} t^{\alpha+n+1} \frac{\alpha}{\alpha+n}\,\frac{\beta^n }{n!}\,\Gamma(\alpha+n+1) \,e^{\pm\mathrm{i}\pi(\alpha+n+1)}\, \mathrm{U}(\alpha+n+1,2\alpha+n+1,e^{\pm\mathrm{i}\pi}(2\beta t)).
\end{equation}
To upper bound the above two functions, we can make use of Olver et al. (2006, Chapter 13.7(ii)) regarding bounds on $\mathrm{U}(a,b,z)$. In particular, we shall put $n=1$, $a=\alpha$ or $a=\alpha+n+1$, and $b= 2\alpha+n+1$ into their Equation (13.7.4). It is noted that the bounds are of slightly different forms depending on  the way $|z|$ tends to infinity. In particular, we will use their result by putting $z=\pm 2\beta t$ (i.e. $z$ is real), which tends to $\pm \infty$ as $t\to\infty$.

{\bf Term $\chi_{n,21}(t)$ of $\chi_{n,2}(t)$.} For the term $\chi_{n,21}(t)$, applying Olver et al. (2006, Chapter 13.7(ii)) with $n=1$, $a=\alpha$ and $b= 2\alpha+n+1$, we have that
\begin{equation}\label{Uexpression1}
\mathrm{U}(\alpha,2\alpha+n+1,2\beta t) =\frac{1}{(2\beta t)^\alpha} +\varepsilon_{n,21}(2\beta t),
\end{equation}
where
\begin{equation}\label{epsilon1Def}
|\varepsilon_{n,21}(2\beta t)|\le 2\omega_{n,21}(t)\,\frac{\alpha(\alpha+n)}{(2\beta t)^{\alpha+1}}\,\exp\left(\frac{2\omega_{n,21}(t)\rho_{n,21}(t)}{2\beta t}\right)~~\text{for}~~ 2\beta t>n+1 .
\end{equation}
We note that the above inequality corresponds to the term $2\beta t$ that belongs to region $R_1$ in Olver et al. (2006, Chapter 13.7(ii)), and therefore it is valid only when $2\beta t>n+1$. The intermediate functions $\sigma_{n,21}(t)$, $ \omega_{n,21}(t)$ and $ \rho_{n,21}(t)$ are given by their Equations (13.7.8) and (13.7.9) as $\sigma_{n,21}(t):=\frac{n+1}{2\beta t}$, $\omega_{n,21}(t):=\frac{1}{1-\sigma_{n,21}(t)}$,
and
\[
\rho_{n,21}(t) :=\frac{1}{2}|2\alpha^2-2\alpha(2\alpha+n+1)+2\alpha+n+1| +\frac{\sigma_{n,21}(t)\left[1+\frac{\sigma_{n,21}(t)}{4}\right]}{[1-\sigma_{n,21}(t)]^2}.
\]
Recall that $\kappa$ is to be chosen `large enough'. Here we shall choose $2\beta t>\kappa(n+1)$ for some $\kappa$ that we, at this point, impose strictly larger than $1$ (so that the constraint in \eqref{epsilon1Def} is satisfied) to arrive at the inequalities $\sigma_{n,21}(t)\le\frac{1}{\kappa}$,
\begin{equation}\label{omegan21bound}
\omega_{n,21}(t)\le\frac{1}{1-\frac{1}{\kappa}}=\frac{\kappa}{\kappa-1},
\end{equation}
and
\begin{equation}\label{rhon21bound}
\rho_{n,21}(t)\le \frac{1}{2}|2\alpha(1-\alpha)+(n+1)(1-2\alpha)| +\frac{\frac{1}{\kappa}\left(1+\frac{1}{4\kappa}\right)}{\left(1-\frac{1}{\kappa}\right)^2}\le B_1(n+1)+C_1,
\end{equation}
where $B_1:=\frac{|1-2\alpha|}{2}$ and $C_1:=\alpha(1-\alpha) +\frac{\kappa+\frac{1}{4}}{(\kappa-1)^2}$ are non-negative constants. A crucial point concerning these two constants is that they only depend on $\kappa$ (that is to be fixed later on) and are always finite under $\kappa>1$. Incorporating the inequalities \eqref{omegan21bound} and \eqref{rhon21bound} into \eqref{epsilon1Def} leads to
\begin{align}\label{epsilon1bound}
|\varepsilon_{n,21}(2\beta t)|\le &~\frac{2\kappa}{\kappa-1}\,\frac{\alpha(\alpha+n)}{(2\beta t)^{\alpha+1}}\,\exp\left(\frac{\frac{2\kappa}{\kappa-1} [B_1(n+1)+C_1]}{2\beta t}\right)\nonumber\\
\le &~\frac{2\kappa\alpha}{\kappa-1}\,\frac{1}{(2\beta t)^\alpha} \,\frac{\alpha+n}{2\beta t} \,\exp\left(\frac{\frac{2}{\kappa-1} [B_1(n+1)+C_1]}{n+1}\right)\nonumber\\
\le &~\frac{M_1^*}{(2\beta t)^\alpha}~~\text{for}~~ 2\beta t>\kappa(n+1) ~~\text{with}~~\kappa>1,
\end{align}
where $M_1^* :=\frac{2\kappa\alpha}{\kappa-1}\,\exp\left(\frac{2(B_1+C_1)}{\kappa-1}\right)$. With \eqref{Uexpression1} and \eqref{epsilon1bound}, we obtain the upper bound
\begin{equation}\label{U1bound}
|\mathrm{U}(\alpha,2\alpha+n+1,2\beta t)| \le \frac{1}{(2\beta t)^\alpha} +|\varepsilon_{n,21}(2\beta t)|\le \frac{M_1^{**}}{(2\beta t)^\alpha}~~\text{for}~~ 2\beta t>\kappa(n+1) ~~\text{with}~~\kappa>1,
\end{equation}
where $M_1^{**}:=M_1^*+1$.

Next, noting that $e^{-x}x^{n+1}$ (as a function of $x$) is decreasing on $[n+1,\infty)$, we have
\begin{equation*}
e^{-2\beta t} (2\beta t)^{n+1} \le e^{-\kappa(n+1)} [\kappa(n+1)]^{n+1} ~~\text{for}~~ 2\beta t>\kappa(n+1) ~~\text{with}~~\kappa>1.
\end{equation*}
Consequently,
\begin{equation*}
t^{n+1} e^{-2\beta t} \beta^n =\frac{1}{\beta}\, \frac{1}{2^{n+1}}\, e^{-2\beta t} (2\beta t)^{n+1} \le \frac{1}{\beta}\, \frac{1}{2^{n+1}}\, e^{-\kappa(n+1)} [\kappa(n+1)]^{n+1} ~~\text{for}~~ 2\beta t>\kappa(n+1) ~~\text{with}~~\kappa>1.
\end{equation*}
Then, with the help of \eqref{U1bound}, $\chi_{n,21}(t)$ in \eqref{chin21Def} can be upper bounded as, for $\kappa>1$,
\begin{align}\label{chin21bound}
|\chi_{n,21}(t)| \le &~I_{\{2\beta t>\kappa(n+1)\}} t^\alpha \frac{1}{\beta}\, \frac{1}{2^{n+1}}\, e^{-\kappa(n+1)} [\kappa(n+1)]^{n+1}\frac{\alpha}{\alpha+n}\,\frac{1}{n!}\,\Gamma(\alpha)\frac{M_1^{**}}{(2\beta t)^\alpha}\nonumber\\
\le &~\frac{\alpha\Gamma(\alpha)M_1^{**}}{\beta(2\beta)^\alpha} \left(\frac{\kappa e^{-\kappa}}{2}\right)^{n+1} \frac{1}{\alpha+n}\,\frac{(n+1)^{n+1}}{n!}
:=~V_{n,21}.
\end{align}
Using \eqref{asymp1} with $\alpha=1$, it is noted that
\begin{equation*}
\frac{(n+1)^{n+1}}{n!} \sim \frac{1}{\sqrt{2\pi}}\,n^{\frac{1}{2}}e^{n+1} ~~\text{as}~~n\rightarrow\infty,
\end{equation*}
and thus for $\kappa>1$,
\begin{equation*}
V_{n,21} \sim \frac{\alpha\Gamma(\alpha)M_1^{**}}{\beta(2\beta)^\alpha\sqrt{2\pi}} \left(\frac{\kappa e^{-\kappa+1}}{2}\right)^{n+1} \frac{n^{\frac{1}{2}}}{\alpha+n} \sim \frac{\alpha\Gamma(\alpha)M_1^{**}}{\beta(2\beta)^\alpha\sqrt{2\pi}} \left(\frac{\kappa e^{-\kappa+1}}{2}\right)^{n+1} \frac{1}{\sqrt{n}} :=V_{n,21}^*~~\text{as}~~n\rightarrow\infty.
\end{equation*}
As a result, application of ratio test gives, for $\kappa>1$,
\begin{equation*}
\lim_{n\rightarrow\infty}\left|\frac{V_{n+1,21}}{V_{n,21}}\right| =\lim_{n\rightarrow\infty}\frac{V_{n+1,21}^*}{V_{n,21}^*} =\lim_{n\rightarrow\infty}\frac{\kappa e^{-\kappa+1}}{2}\sqrt{\frac{n}{n+1}} =\frac{\kappa e^{-\kappa+1}}{2}\le\frac{1}{2}<1,
\end{equation*}
where the second last inequality follows from the fact that $xe^{-x+1}/2$ (as a function of $x$) is decreasing on $[1,\infty)$. Hence,
\begin{equation}\label{Vn21bound}
\sum^\infty_{n=0}V_{n,21}<\infty.
\end{equation}

{\bf Term $\chi_{n,22}(t)$ of $\chi_{n,2}(t)$.} Lastly, we consider the term $\chi_{n,22}(t)$ in (\ref{chin22Def}). Applying again Olver et al. (2006, Chapter 13.7(ii)) with $n=1$, $a=\alpha+n+1$ and $b= 2\alpha+n+1$ implies
\begin{equation}\label{Uexpression2}
\mathrm{U}(\alpha+n+1,2\alpha+n+1,e^{\pm\mathrm{i}\pi}(2\beta t)) =\frac{1}{[e^{\pm\mathrm{i}\pi}(2\beta t)]^{\alpha+n+1}} +\varepsilon_{n,22}(e^{\pm\mathrm{i}\pi}(2\beta t)).
\end{equation}
The term $e^{\pm\mathrm{i}\pi}(2\beta t)$ in the above Tricomi's function this time tends to $-\infty$ regardless of whether one chooses to multiply $2\beta t$ by $e^{\mathrm{i}\pi}$ or $e^{-\mathrm{i}\pi}$. Even though the value of $ \mathrm{U}(\alpha+n+1,2\alpha+n+1,e^{\pm\mathrm{i}\pi}(2\beta t))$ is different depending on whether we choose the argument $e^{\mathrm{i}\pi}(2\beta t)$ or $e^{-\mathrm{i}\pi}(2\beta t)$, the upcoming bounds \eqref{epsilon2bound} and \eqref{U2bound} are the same as they only involve the modulus of $e^{\pm\mathrm{i}\pi}(2\beta t)$. Hence, we choose to keep the notation $e^{\pm\mathrm{i}\pi}$ for presentation purpose. Thus, the error term $\varepsilon_{n,22}(e^{\pm\mathrm{i}\pi}(2\beta t))$ in \eqref{Uexpression2} satisfies the inequality
\begin{equation}\label{epsilon2Def}
|\varepsilon_{n,22}(e^{\pm\mathrm{i}\pi}(2\beta t))|\le 2\omega_{n,22}(t)\zeta_{n,22}(t) \,\frac{(\alpha+n+1)(1-\alpha)}{(2\beta t)^{\alpha+n+2}}\,\exp\left(\frac{2\omega_{n,22}(t)\rho_{n,22}(t)\zeta_{n,22}(t)}{2\beta t}\right)
\end{equation}
for $2\beta t>2(n+1)$. Note that, contrary to (\ref{epsilon1Def}), the above inequality corresponds this time to $e^{\pm\mathrm{i}\pi}(2\beta t)$ which belongs to region $R_3\cup\overline{R}_3$ in Olver et al. (2006, Chapter 13.7(ii)), and hence it is valid for $2\beta t>2(n+1)$. The functions $\sigma_{n,22}(t)$, $ \omega_{n,22}(t)$ and $ \rho_{n,22}(t)$ are given by their Equations (13.7.8) and (13.7.9) as $\sigma_{n,22}(t):=\frac{n+1}{2\beta t}\le\frac{1}{2}$,
\begin{equation*}
\nu_{n,22}(t):=\left(\frac{1}{2}+\frac{1}{2}\sqrt{1-4[\sigma_{n,22}(t)]^2}\right)^{-\frac{1}{2}}\le\sqrt{2},
\end{equation*}
\begin{equation*}
\zeta_{n,22}(t):=\left(\frac{\pi}{2} +\sigma_{n,22}(t)[\nu_{n,22}(t)]^2\right)\nu_{n,22}(t) \le \left(\frac{\pi}{2}+1\right)\sqrt{2},
\end{equation*}
\begin{equation*}
\omega_{n,22}(t):=\frac{1}{1-\nu_{n,22}(t)\sigma_{n,22}(t)}\le \frac{1}{1-\frac{\sqrt{2}}{2}} =2+\sqrt{2},
\end{equation*}
and
\begin{align*}
\rho_{n,22}(t) :=&~\frac{1}{2}|2(\alpha+n+1)^2-2(\alpha+n+1)(2\alpha+n+1)+2\alpha+n+1| +\frac{\nu_{n,22}(t)\sigma_{n,22}(t)\left[1+\frac{\nu_{n,22}(t)\sigma_{n,22}(t)}{4}\right]}{[1-\nu_{n,22}(t)\sigma_{n,22}(t)]^2}\\
\le&~\frac{1}{2}|2\alpha(1-\alpha)+(n+1)(1-2\alpha)| +\frac{19+14\sqrt{2}}{4}\le~B_2(n+1)+C_2,
\end{align*}
where $B_2:=\frac{|1-2\alpha|}{2}$ and $C_2:=\alpha(1-\alpha) +\frac{19+14\sqrt{2}}{4}$ are non-negative constants. By substituting the above results into \eqref{epsilon2Def}, we arrive at
\begin{align}\label{epsilon2bound}
|\varepsilon_{n,22}(e^{\pm\mathrm{i}\pi}(2\beta t))|\le &~2(1+\sqrt{2})(\pi+2) \frac{(\alpha+n+1)(1-\alpha)}{(2\beta t)^{\alpha+n+2}}\,\exp\left(\frac{2(1+\sqrt{2})(\pi+2)[B_2(n+1)+C_2]}{2\beta t}\right)\nonumber\\
\le &~2(1+\sqrt{2})(\pi+2) (1-\alpha) \frac{1}{(2\beta t)^{\alpha+n+1}}\,\frac{\alpha+n+1}{2\beta t} \,\exp\left(\frac{(1+\sqrt{2})(\pi+2)[B_2(n+1)+C_2]}{n+1}\right)\nonumber\\
\le &~\frac{M_2^*}{(2\beta t)^{\alpha+n+1}}~~\text{for}~~ 2\beta t>2(n+1),
\end{align}
where $M_2^* :=2(1+\sqrt{2})(\pi+2) (1-\alpha) e^{(1+\sqrt{2})(\pi+2)(B_2+C_2)}$. Using \eqref{epsilon2bound}, an upper bound for \eqref{Uexpression2} is given by
\begin{equation}\label{U2bound}
|\mathrm{U}(\alpha+n+1,2\alpha+n+1,e^{\pm\mathrm{i}\pi}(2\beta t))| \le \frac{1}{(2\beta t)^{\alpha+n+1}} +|\varepsilon_{n,22}(e^{\pm\mathrm{i}\pi}(2\beta t))| \le \frac{M_2^{**}}{(2\beta t)^{\alpha+n+1}}~~\text{for}~~ 2\beta t>2(n+1),
\end{equation}
where $M_2^{**}:=M_2^*+1$.

Now, we can get from \eqref{chin22Def} and \eqref{U2bound} the upper bound, for $\kappa\ge 2$,
\begin{align}\label{chin22bound}
|\chi_{n,22}(t)| \le &~I_{\{2\beta t>\kappa(n+1)\}} t^{\alpha+n+1} \frac{\alpha}{\alpha+n}\,\frac{\beta^n }{n!}\,\Gamma(\alpha+n+1)\, \frac{M_2^{**}}{(2\beta t)^{\alpha+n+1}}\nonumber\\
\le &~\frac{\alpha M_2^{**}}{(2\beta)^{\alpha+1}}\,\frac{1}{2^n}\,\frac{1}{\alpha+n}\,\frac{\Gamma(\alpha+n+1)}{n!}:=~V_{n,22}.
\end{align}
Utilizing the Stirling's formula \eqref{stirling1} and the equivalent for Gamma function (see Abramowitz and Stegun (1972, Equation (6.1.37), p.257), namely
\begin{equation*}
\Gamma(y+1)\sim\sqrt{2\pi y}\left(\frac{y}{e}\right)^y~~\text{as}~~y\rightarrow\infty,
\end{equation*}
we obtain the asymptotic relationship
\begin{equation*}
\frac{\Gamma(\alpha+n+1)}{n!} \sim \frac{\sqrt{\alpha+n}\left(\frac{\alpha+n}{e}\right)^{\alpha+n}}{\sqrt{n}\left(\frac{n}{e}\right)^n} =\sqrt{\frac{\alpha+n}{n}}\,e^{-\alpha} \left(\frac{\alpha+n}{n}\right)^n (\alpha+n)^\alpha \sim n^\alpha ~~\text{as}~~n\rightarrow\infty.
\end{equation*}
Thus, for $\kappa\ge 2$,
\begin{equation*}
V_{n,22} \sim \frac{\alpha M_2^{**}}{(2\beta)^{\alpha+1}}\,\frac{1}{2^n}\,\frac{1}{\alpha+n}\,n^\alpha  \sim \frac{\alpha M_2^{**}}{(2\beta)^{\alpha+1}}\,\frac{n^{\alpha-1}}{2^n} :=V_{n,22}^*~~\text{as}~~n\rightarrow\infty.
\end{equation*}
Consequently, the ratio test verifies that, for $\kappa\ge 2$,
\begin{equation*}
\lim_{n\rightarrow\infty}\left|\frac{V_{n+1,22}}{V_{n,22}}\right| =\lim_{n\rightarrow\infty}\frac{V_{n+1,22}^*}{V_{n,22}^*} =\lim_{n\rightarrow\infty}\frac{1}{2}\left(\frac{n+1}{n}\right)^{\alpha-1} =\frac{1}{2}<1,
\end{equation*}
and hence,
\begin{equation}\label{Vn22bound}
\sum^\infty_{n=0}V_{n,22}<\infty.
\end{equation}

{\bf Verifying \eqref{chinbound} and \eqref{Vnbound}.} Combining \eqref{chindecompose}, \eqref{chin2decompose}, \eqref{chin1bound}, \eqref{chin21bound} and \eqref{chin22bound}, it is clear that by choosing $\kappa\ge 2$ in \eqref{chin1Def} and \eqref{chin2Def}, we have
\begin{equation*}
|\chi_n(t)| = |\chi_{n,1}(t)+\chi_{n,21}(t)+\chi_{n,22}(t)| \le |\chi_{n,1}(t)|+|\chi_{n,21}(t)|+|\chi_{n,22}(t)| \le V_{n,1}+V_{n,21}+V_{n,22},
\end{equation*}
i.e. the upper bound \eqref{chinbound} is obtained by defining the sequence $\{V_n\}_{n=0}^\infty$ via $V_n:=V_{n,1}+V_{n,21}+V_{n,22}$. Then, the condition \eqref{Vnbound} is satisfied thanks to \eqref{Vn1bound}, \eqref{Vn21bound} and \eqref{Vn22bound}.

\subsection*{Part 2: Proof of Lemma 1
}

First, we split $G^*(y)$ defined in (4.28) 
into two parts as
\begin{equation}\label{G1ya}
G^*(y) = \int_{0}^{1/2}(1-x)^{-\eta} x^{\alpha -1}\mathrm{d}x +\int_{1/2}^{y}(1-x)^{-\eta}x^{\alpha -1}\mathrm{d}x,
\end{equation}
where the first integral is always finite. Therefore, we focus on the divergent second integral and distinguish between the cases $\eta>1$ and $\eta=1$ for ease of presentation.

\noindent\textbf{Case 1.} $\eta>1$: Performing integration by parts, one obtains
\begin{equation}\label{integral2IBP}
\int_{1/2}^{y}(1-x)^{-\eta}x^{\alpha -1}\mathrm{d}x = \frac{1}{\eta-1}y^{\alpha-1}(1-y)^{-\eta+1} - \frac{1}{2^{-\eta+\alpha}}\frac{1}{\eta-1} + \frac{1-\alpha}{\eta-1}\int_{1/2}^{y}(1-x)^{-\eta+1}x^{\alpha-2}\mathrm{d}x.
\end{equation}
The analysis is further separated into three cases as follows.

\noindent\textbf{Case 1a.} $1<\eta<2$: In this case, $-\eta+1>-1$ and therefore the integral $\int_{1/2}^{1}(1-x)^{-\eta+1}x^{\alpha-2}\mathrm{d}x$ is convergent as $x^{\alpha -2}$ is bounded on $x\in[\frac{1}{2},1]$. Incorporating such an observation into \eqref{integral2IBP}, it is clear that
\begin{equation}\label{lem2case1astep1}
\int_{1/2}^{y}(1-x)^{-\eta}x^{\alpha -1}\mathrm{d}x = \frac{1}{\eta-1}y^{\alpha-1}(1-y)^{-\eta+1} + C +o(1)\text{~~as~~}y\rightarrow 1^-
\end{equation}
for some constant $C$. Further expanding $y^{\alpha-1}=[1+(y-1)]^{\alpha-1}=1+(\alpha-1)(y-1)+o(y-1)$, we obtain
\begin{align}
y^{\alpha-1}(1-y)^{-\eta+1}=&~(1-y)^{-\eta+1}+(1-\alpha)(1-y)^{-\eta+2}+o((1-y)^{-\eta+2})\label{lem2case1astep20} \\
=&~(1-y)^{-\eta+1}+o(1)\text{~~as~~}y\rightarrow 1^-,\label{lem2case1astep2}
\end{align}
where the last line is due to $-\eta+2>0$. Combining \eqref{G1ya}, \eqref{lem2case1astep1} and \eqref{lem2case1astep2} yields the first asymptotic result in (4.30). 

\noindent\textbf{Case 1b.} $\eta=2$: 
Equation \eqref{integral2IBP} is still valid, but the integral $\int_{1/2}^{y}(1-x)^{-1}x^{\alpha-2}\mathrm{d}x$ on the right-hand side can easily shown to be divergent at $y=1$. Integration by parts gives
\begin{equation}\label{lem2case1bstep1}
\int_{1/2}^{y}(1-x)^{-1}x^{\alpha-2}\mathrm{d}x =-y^{\alpha -2}\ln(1-y)-2^{2-\alpha}\ln 2 + (\alpha -2) \int_{1/2}^{y}x^{\alpha- 3}\ln(1-x) \mathrm{d}x,
\end{equation}
where the integral $\int_{1/2}^{y}x^{\alpha- 3}\ln(1-x) \mathrm{d}x$ converges as $y\rightarrow1^-$ because $\int_{1/2}^{1}\ln(1-x) \mathrm{d}x$ is finite. Incorporating \eqref{integral2IBP} and \eqref{lem2case1bstep1} into \eqref{G1ya} and utilizing the expansion $y^{a}=[1+(y-1)]^{a-1}=1+a(y-1)+o(y-1)$ at $a=\alpha-1$ and $a=\alpha-2$, one finds (for some constant $C$)
\begin{align*}
G^*(y) =&~y^{\alpha-1}(1-y)^{-1} -(1-\alpha)y^{\alpha -2}\ln(1-y) +C +o(1)\\
=&~(1-y)^{-1} -(\alpha-1)-(1-\alpha)\ln(1-y)-(1-\alpha)(\alpha-2)(y-1)\ln(1-y)\\
&-(1-\alpha)(\alpha-2)o(y-1)\ln(1-y)+C +o(1)\text{~~as~~}y\rightarrow 1^-,
\end{align*}
which simplifies to the second asymptotic formula in (4.30) 
because $(y-1)\ln(1-y)=o(1)$.

\noindent\textbf{Case 1c.} $\eta>2$: The integral on the right-hand side of \eqref{integral2IBP} satisfies
\begin{equation*}
\int_{1/2}^{y}(1-x)^{-\eta+1}x^{\alpha-2}\mathrm{d}x =O\Big(\int_{1/2}^{y}(1-x)^{-\eta +1} \mathrm{d}x\Big) =O((1-y)^{-\eta+2})=o((1-y)^{-\eta+1})\text{~~as~~}y\rightarrow 1^-.
\end{equation*}
By applying the expansion \eqref{lem2case1astep20} (which is also valid for $\eta>2$) and the above result to \eqref{integral2IBP}, one proves the first result in (4.29) 
when $\eta>2$.

It is instructive to note that the proven result (4.30) 
also asserts that the first asymptotic formula in (4.29) 
is also valid when $1<\eta\le2$.


\noindent\textbf{Case 2.} $\eta=1$: Similar to \eqref{lem2case1bstep1}, integrating by parts on the second integral in (\ref{G1ya}) results in
\begin{equation*}
\int_{1/2}^{y}(1-x)^{-1}x^{\alpha -1}\mathrm{d}x =-y^{\alpha -1}\ln(1-y)-2^{1-\alpha}\ln 2 + (\alpha -1) \int_{1/2}^{y}x^{\alpha- 2}\ln(1-x) \mathrm{d}x.
\end{equation*}
Since $\int_{1/2}^{1}\ln(1-x) \mathrm{d}x$ is finite, the integral on the right-hand side above is convergent as $y\rightarrow1^-$. The second result in (4.29) 
follows by further noting that $\lim_{y\rightarrow1^-}(1-y^{\alpha -1})\ln(1-y)=0$.

\subsection*{Part 3: Proof of Lemma 2
}

Performing a change of variable $x=(v+s)(1-y)$ (i.e. $y=1-\frac{x}{v+s}$) in (4.33) 
yields
\begin{equation}\label{Js0}
J(s) = (v+s)^{1-\gamma -\xi}\int_{\frac{v}{v+s}}^{1}y^{-\gamma}(1-y)^{-\xi}G\bigg(1-\frac{v}{(v+s)y}\bigg)\mathrm{d}y.
\end{equation}
To analyze the asymptotic behaviour of $J(s)$ as $s\rightarrow\infty$ (so that $\frac{v}{v+s}\rightarrow 0^+$), we shall study the integral, for small $h>0$,
\begin{equation}\label{Js1}
\int_{h}^{1}y^{-\gamma}(1-y)^{-\xi}G\bigg(1-\frac{h}{y}\bigg)\mathrm{d}y=U_1(h)+U_2(h),
\end{equation}
where
\begin{equation}\label{U12h}
U_1(h):= \int_{h}^{1/2}y^{-\gamma}(1-y)^{-\xi}G\bigg(1-\frac{h}{y}\bigg)\mathrm{d}y\text{~~and~~}
U_2(h):=\int_{1/2}^{1}y^{-\gamma}(1-y)^{-\xi}G\bigg(1-\frac{h}{y}\bigg)\mathrm{d}y.
\end{equation}
It is clear that $U_2(h)$ converges to $G(1)\int_{1/2}^{1}y^{-\gamma}(1-y)^{-\xi}\mathrm{d}y$ as $h\rightarrow0^+$ by dominated convergence (as $G(1-\frac{h}{y})$ is bounded on $h\in (0,\frac{1}{2}]$ and $y\in[\frac{1}{2},1]$ and $\int_{1/2}^{1}y^{-\gamma}(1-y)^{-\xi}\mathrm{d}y$ is finite). For $U_1(h)$, a change of variable $z=1-\frac{h}{y}$ (i.e. $y=\frac{h}{1-z}$ with $\mathrm{d}y =\frac{h}{(1-z)^2}\mathrm{d}z$) results in
\begin{equation}\label{U1h}
U_1(h)= h^{1-\gamma}\int_{0}^{1-2h}(1-z)^{\gamma - 2}\bigg(1-\frac{h}{1-z}\bigg)^{-\xi}G(z)\mathrm{d}z.
\end{equation}
We need to consider two cases as follows.

\noindent\textbf{Case 1.} $\gamma-2>-1 \Leftrightarrow \gamma >1$: We rewrite \eqref{U1h} as
\begin{equation}\label{U1h2}
U_1(h)= h^{1-\gamma}\int_{0}^{1}(1-z)^{\gamma - 2}\omega(z,h)\mathrm{d}z,
\end{equation}
where $\omega(z,h) := I_{\{0\le z\le 1-2h\}}(1-\frac{h}{1-z})^{-\xi}G(z)$. It can be seen that $\omega(z,h)$ is bounded by some constant for all $z\in [0,1]$ and a small enough $h$. Besides, one has that $\int_{0}^{1}(1-z)^{\gamma-2}\mathrm{d}z$ is finite. Hence, by the dominated convergence theorem we obtain
\begin{equation*}
\lim_{h\rightarrow 0^+}\int_{0}^{1}(1-z)^{\gamma-2}\omega(z,h)\mathrm{d}z = \int_{0}^{1}(1-z)^{\gamma-2}w(z,0)\mathrm{d}z =\int_{0}^{1}(1-z)^{\gamma-2}G(z)\mathrm{d}z.
\end{equation*}
Consequently, \eqref{U1h2} (or (\ref{U1h})) behaves asymptotically as
\begin{equation}\label{*1}
U_1(h) \sim  h^{1-\gamma} \int_{0}^{1}(1-z)^{\gamma-2}G(z)\mathrm{d}z \text{~~as~~}h\rightarrow 0^+.
\end{equation}

\noindent\textbf{Case 2.} $\gamma-2 = -1 \Leftrightarrow \gamma=1$: By the change of variable $t=\frac{\ln(1-z)}{\ln2h}$ (i.e. $1-z =(2h)^t$ with $\mathrm{d}t = \frac{1}{\ln2h}\frac{-1}{1-z}\mathrm{d}z$), (\ref{U1h}) becomes
\begin{equation}\label{U1ha}
U_1(h)=(-\ln2h) \int_{0}^{1}(1-2^{-t}h^{1-t})^{-\xi}G(1-(2h)^t)\mathrm{d}t.
\end{equation}
For $t\in[0,1]$ and $h$ small enough (less than $1/2$), it is observed that both $(1-2^{-t}h^{1-t})^{-\xi}$ and $G(1-(2h)^t)$ are upper bounded by some constant. Therefore, one checks easily that
\[
\lim_{h\rightarrow 0^+}(1-2^{-t}h^{1-t})^{-\xi}G(1-(2h)^t) =
\begin{cases}
G(0), &\mathrm{if}~~t=0,
\\
G(1), &\mathrm{if}~~0<t<1,\\
2^{\xi}G(1), & \mathrm{if}~~t=1.
\end{cases}
\]
Hence, by the dominated convergence theorem together with the fact that $-\ln 2h = -\ln 2 -\ln h \sim -\ln h$ as $h\rightarrow 0^+$, one finds that the asymptotic behaviour of (\ref{U1ha}) (or (\ref{U1h}) with $\gamma=1$) is given by
\begin{equation}\label{*2}
U_1(h)  \sim (-\ln h)G(1)\text{~~as~~}h\rightarrow 0^+.
\end{equation}	

Therefore, combining the two cases from (\ref{*1}) and (\ref{*2}), we have that $U_1(h)$ is $O(h^{1-\gamma})$ for $\gamma>1$ but $O(-\ln h)$ for $\gamma=1$. Since $U_2(h)$ in (\ref{U12h}) is $O(1)$, it can be concluded that (\ref{Js1}) asymptotically behaves like, as $s\rightarrow\infty$,
\[
\int_{h}^{1}y^{-\gamma}(1-y)^{-\xi}G\bigg(1-\frac{h}{y}\bigg)\mathrm{d}y \sim
\begin{cases}
h^{1-\gamma}\int_{0}^{1}(1-z)^{\gamma-2}G(z)\mathrm{d}z,& \text{if }\gamma>1, \\
  (-\ln h)G(1),& \text{if }\gamma =1.
\end{cases}
\]
Finally, plugging $h = \frac{v}{v+s} $ in (\ref{Js0}) results in (4.34). 

\end{appendices}	

\end{document}